\input ./style/arxiv-general.cfg
\documentclass[labelcite,MSNbibl,number,citesort,dvips]{arxbj}
\makeatletter
   \@ifpackageloaded{graphicx}{}{\usepackage{graphicx}}
\makeatother
\usepackage{upgreek}

\volume{22}
\issue{2}
\pubyear{2016}
\firstpage{1227}
\lastpage{1277}
\doi{10.3150/14-BEJ691}
\docsubty{FLA}

\makeatletter
\newcommand{\error}{\operatorname{error}}
\newcommand{\algo}{\operatorname{algo}}
\newcommand{\diag}{\operatorname{diag}}

\newcommand{\rrvert}{\vert}

\newcommand{\llvert}{\vert}
\newtheorem{theorem}{Theorem}
\newproclaim{condition}{Condition}
\newtheorem{corollary}{Corollary}
\newremark{example}{Example}
\newtheorem{lemma}{Lemma}
\newremark{remark}{Remark}
\makeatother

\begin{document}
\begin{frontmatter}

\title{Greedy algorithms for prediction}
\runtitle{Greedy algorithms for prediction}

\begin{aug}
\author{\inits{A.}\fnms{Alessio}~\snm{Sancetta}\corref{}\ead[label=e1]{asancetta@gmail.com}\ead[label=u1,url]{https://sites.google.com/site/wwwsancetta/}}
\address{Department of Economics, Royal Holloway University of London,
Egham Hill, Egham TW20 0EX, UK.
\printead{e1}; \printead{u1}}
\end{aug}

%
\received{\smonth{8} \syear{2013}}
%
\revised{\smonth{12} \syear{2014}}

%
\begin{abstract}
In many prediction problems, it is not uncommon that the number of
variables used to construct a forecast is of the same order of magnitude
as the sample size, if not larger. We then face the problem of constructing
a prediction in the presence of potentially large estimation error.
Control of the estimation error is either achieved by selecting variables
or combining all the variables in some special way. This paper considers
greedy algorithms to solve this problem. It is shown that the resulting
estimators are consistent under weak conditions. In particular, the
derived rates of convergence are either minimax or improve on the
ones given in the literature allowing for dependence and unbounded
regressors. Some versions of the algorithms provide fast solution
to problems such as Lasso.
\end{abstract}

%
\end{frontmatter}

\section{Introduction}\label{sec1}

The goal of this paper is to address the problem of forecasting in
the presence of many explanatory variables or individual forecasts.
Throughout the paper, the explanatory variables will be referred to
as regressors even when they are individual forecasts that we wish
to combine or basis functions, or in general elements in some dictionary.

The framework is the one where the number of regressors is often large
relatively to the sample size. This is quite common in many fields,
for example, in macroeconomic predictions (e.g., Stock and Watson \cite
{StoWat99,StoWat03,StoWat04}). Moreover, when there is evidence of
structural breaks,
it is not always possible to use the full sample without making further
assumptions. Indeed, it is often suggested to forecast using different
sample sizes in an effort to mitigate the problem (e.g., Pesaran \textit{et~al.}
\cite{PesPetTim06}, Pesaran and Picks
\cite{PesPic11}). When doing so, we still need
to make sure that the forecasts built using smaller sample sizes are
not too noisy.

For these reasons, it is critical to consider procedures that allow
us to select and/or combine variables in an optimal way when the data
are dependent. It is clear that in large-dimensional problems, variable
selection via information criteria is not feasible, as it would require
the estimation of a huge number of models. For example, if we are
considering 100 regressors, naive model selection of a model with
only 10 variables (i.e., an order of magnitude lower) would require
estimation and comparison of ${100\choose 10}$ models, which is in the order of billions.

This paper considers greedy algorithms to do automatic variable selection.
There are many references related to the algorithms considered here
(e.g., B{\"u}hlmann \cite{Buh06}, Barron \textit{et~al.}
\cite{Baretal08}, Huang, Cheang and Barron \cite{Hua08},
B{\"u}hlmann and van~de Geer \cite{Buhvan11}). These existing results
are not applicable
to standard prediction problems, as they assume i.i.d. random variable
with bounded regressors and in some case bounded error terms.

Greedy algorithms have been applied to time series problems both in
a linear and non-linear context (e.g., Audrino and B\"uhlmann \cite{AudBuh06,AudBuh09},
Audrino and Barone-Adesi \cite{AudBar06}, amongst
others). However, to the author's knowledge, in the linear case, only
Lutz and B{\"u}hlmann \cite{LutBuh06} derive consistency under
strong mixing. There,
no rates of convergence are given. (See Audrino and B{\"u}hlmann
\cite{AudBuh09},
for the non-linear case, again where no rates are given.) The above
references only consider Boosting. It is known that other greedy algorithms
possess better convergence rates (e.g., Barron \textit{et~al.}
\cite{Baretal08}). Here,
only linear predictions are considered. Of course, when the regressors
are a basis for some function space, the results directly apply to
series estimators, hence, non-linear prediction (e.g., Mallat and Zhang
\cite{MalZha93}, Daubechies, Defrise and De Mol \cite{DauDefDeM04}, Barron \textit
{et~al.} \cite{Baretal08},
B{\"u}hlmann and van~de Geer \cite{Buhvan11}, Sancetta
\cite{San14}, for more details along these lines).

To be precise, this paper shall consider greedy algorithms and provide
rates of convergence which are best possible for the given set up
or considerably improve on the existing ones, even under dependence
conditions. The first algorithm is the $L_{2}$-Boosting studied by
B{\"u}hlmann \cite{Buh06}, also known as Projection Pursuit
in signal processing
(e.g., Mallat and Zhang \cite{MalZha93}) and Pure Greedy Algorithm in approximation
theory (e.g., DeVore and Temlyakov \cite{DeVTem96}). As mentioned
above, it is
routinely used in many applications, even in time series problems.
The second algorithm is known as Orthogonal Greedy Algorithm (OGA)
in approximation theory (e.g., DeVore and Temlyakov \cite
{DeVTem96}, Temlyakov \cite{Tem00}), and has
also been studied in the statistical literature (Barron \textit{et~al.}
\cite{Baretal08}). It is the one the most resembles
OLS estimation. The
OGA is also reviewed in B{\"u}hlmann and van~de Geer \cite{Buhvan11},
where it
is called Orthogonal Matching Pursuit (see also Zhang
\cite{Zha09}, Cai and Wang \cite{CaiWan11}, for
recent results). The third algorithm is a version
of the Hilbert space projection algorithm studied by Jones \cite{Jon92}
and Barron \cite{Bar93} with the version studied in
this paper taken from
Barron \textit{et~al.} \cite{Baretal08}, and called the Relaxed Greedy Algorithm (RGA).
Adding a natural restriction to the RGA, the algorithm leads to the
solution of the Lasso problem, which appears to be relatively new
(see Sancetta \cite{San14}). This constrained version
will be called Constrained
Greedy Algorithm (CGA). Finally, closely related to the CGA is the
Frank--Wolfe Algorithm (FWA) (see Frank and Wolfe \cite
{FraWol56}, and
Clarkson \cite{Cla10}, Jaggi \cite
{Jag13}, and Freund, Grigas and Mazumder \cite{FreGriMaz13}, for recent
results). This
selection seems to span the majority of known algorithms used in applied
work.

The general problem of variable selection is often addressed relying
on penalized estimation with an $l_{1}$ penalty. Greedy algorithms
can be related to Lasso as they both lead to automatic variable selection.
Algorithms that use a penalty in the estimation will not be discussed
here. It is well known (Friedman \textit{et~al.} \cite{Frietal07})
that the Lasso solution
can be recovered via Pathwise Coordinate Optimization (a stagewise
recursive algorithm), using the results of Tseng
\cite{Tse01} (see also
Daubechies, Defrise and De Mol \cite{DauDefDeM04}, for related results). On the other hand,
Huang, Cheang and Barron \cite{Hua08} have extended the RGA to the case of a
Lasso penalty.
(For recent advances on asymptotics for Lasso, the reader may consult
Greenshtein and Ritov \cite{vanvan04}, Bunea, Tsybakov and Wegkamp \cite{BunTsyWeg07N1},
van~de Geer \cite{van08},
Huang, Cheang and Barron \cite{Hua08}, Zhang \cite{Zha10}, Belloni and
Chernozhukov \cite{BelChe11}, Belloni \textit{et~al.} \cite
{Beletal12}, amongst others.) Another related approach for variable
selection under sparsity and design matrix constraints is via linear
programming (e.g., Candes and Tao \cite{CanTao07}).

One related question which is also considered here is the one of persistence
as defined by Greenshtein and Ritov \cite{vanvan04} and explored
by other authors
(e.g., Greenshtein \cite{Gre06}, B{\"u}hlmann and van~de Geer \cite{Buhvan11}, Bartlett, Mendelson and Neeman \cite{BarMenNee12}). This problems is of interest in a prediction context
and relates to the idea of pseudo true value. Loosely speaking, one
is interested in finding the largest class of linear models relative
to which the estimator is still optimal in some sense. Here, it is
shown that for mixing data, persistence holds for the class of linear
models as large as the ones considered in Greenshtein and Ritov
\cite{vanvan04}
and Bartlett, Mendelson and Neeman~\cite{BarMenNee12}.

The focus of the paper is on prediction and consistency of the forecasts.
Asymptotic normality of the estimators is not derived due to the weak
conditions used (e.g., see B{\"u}hlmann \cite{Buh13},
Nickl and van~de Geer \cite{Nicvan13}, van~de Geer \textit{et~al.}
\cite{vanetal14}, Zhang and Zhang \cite{ZhaZha14} for results
on statistical significance for high-dimensional, sparse models, under
different estimation procedures and assumptions).

The paper is structured as follows. The remainder of this section
defines the estimation set-up, the objectives and the conditions to
be used. Two different sets of dependence conditions are used: beta
mixing, which gives the best convergence rates, and more general conditions
allowing for non-mixing data and possibly long memory data.
Section~\ref{SectionAlgorithms} starts with a summary of existing results
comparing them with some of the ones derived here. The actual statement
of all the results follows afterward. With the exception of the PGA,
it is shown that the algorithms can achieve the minimax rate under
beta mixing. However, for the PGA, the rates derived here considerably
improve on the ones previously obtained. The algorithms are only reviewed
later on in Section~\ref{SectionreviewAlgorithms}. The reader unfamiliar
with these algorithms can browse through Section~\ref{SectionAlgorithms}
right after Section~\ref{sec1} if needed. A discussion of the conditions and
examples and applications of the results are given in Section~\ref
{Sectiondiscussion}.
In particular, Section~\ref{Sectionexamples} gives examples of applications
to long memory achieving convergence rates as good or better than
the ones derived by other authors under i.i.d. observations, though
requiring the population Gram matrix of the regressors to have full
rank. In Section~\ref{Sectionimplementation}, details on implementation
are given. Section~\ref{Sectionimplementation} contains remarks
of practical nature including vectorized versions of the algorithms,
which are useful when implemented in scripting languages such as R
and Matlab. This section also gives details on finite sample performance
via simulation examples to complement the theoretical results. For
example, the simulations in Section~\ref{SectionnumericalResults}
show that~-- despite the slower rates of convergence~-- the PGA seems
to perform particularly well when the signal to noise is low (see
also B{\"u}hlmann and van~de Geer \cite{Buhvan11}, Section~12.7.1.1).
The proofs
are all in Section~\ref{SectionProofs}. Section~\ref{SectionProofs}
contains results on the approximation properties of the algorithms
that can be of interest in their own. For example, a simple extension
of the result in DeVore and Temlyakov \cite{DeVTem96} to
statistical least
square estimation is given in order to bound the approximation error
of the PGA ($L_{2}$-Boosting). Moreover, it is also shown that the
complexity of the PGA grows sub-linearly with the number of iteration,
hence compensating this way for the higher approximation error (Lemma
\ref{LemmaEstimatorL1Bound} in Section~\ref{SectionProofs}). This
observation appears to be new and it is exploited when considering
convergence under non-mixing data.

\subsection{Estimation setup}

There are possibly many more regressors than the sample size. However,
most of the regressors are not needed or useful for prediction, for example,
they may either be zero or have a progressively decreasing importance.
This means that most of the regressors are redundant. Redundancy is
formally defined in terms of a bound on the absolute sum of the regression
coefficients. In particular, let $\mathcal{X}$ be a set of regressors
of cardinality $K$, possibly much larger than the sample size $n$
and growing with $n$ if needed. Then the focus is on the linear
regression function $\mu(x )=\sum_{k=1}^{K}b_{k}x^{
(k )}$
where $\sum_{k=1}^{K}\llvert b_{k}\rrvert \leq B<\infty$, and
$x^{ (k )}$
is the $k$th element in $x$. As $B$ increases, the class of
functions representable by $\mu$ becomes larger (e.g., when $\mathcal{X}$
is a set of functions whose linear span is dense in some space of
functions). The same remark is valid when $K$ grows with $n$, as
for sieve estimators. The absolute summability of the regression coefficients
is standard (e.g., B{\"u}hlmann \cite{Buh06}, Barron \textit{et~al.}
\cite{Baretal08}). This
restriction
is also used in compressed sensing, where a signal with no noise admits
a sparse representation in terms of a dictionary (e.g., Temlyakov \cite{Tem11}, Chapter~5). Nevertheless, high-dimensional statistics also considers
the problem of consistency when $B\rightarrow\infty$ at the rate
$\mathrm{o} (\sqrt{n/\ln K} )$ (e.g., Greenshtein and Ritov
\cite{vanvan04}, Greenshtein \cite{Gre06},
B{\"u}hlmann and van~de Geer \cite{Buhvan11}, Bartlett, Mendelson and Neeman \cite{BarMenNee12}). Here, it is shown that Greedy algorithms are consistent in
this situation when the data are dependent and the regressors are
not necessarily bounded.

Notational details and conditions are introduced next. Given random
variables $Y$ and $X$, interest lies in approximating the conditional
regression\vspace*{1pt} function $\mathbb{E} [Y|X ]=\mu_{0}
(X )$,
with the linear regression $\mu(X ):=\sum_{k=1}^{K}b_{k}X^{ (k )}$,
where $\sum_{k=1}^{K}\llvert b_{k}\rrvert \leq B$, and $X^{
(k )}$
denotes the $k$th element of $X$. Hence, $\mu_{0}$ does not
need to be linear. (Most of the literature, essentially, considers
the case when the true regression function $\mu_{0}\in\mathcal
{L} (B )$
with Barron \textit{et~al.} \cite{Baretal08} being one of the few
exceptions.) Let
$ \{ Y_{i},X_{i}\dvt i=1,2,\ldots,n \} $
be possibly dependent copies of $Y,X$. Define the empirical inner
product
\[
\bigl\langle Y,X^{ (k )} \bigr\rangle_{n}:=\frac
{1}{n}\sum
_{i=1}^{n}Y_{i}X_{i}^{ (k )}
\quad\mbox{and}\quad\bigl\llvert X^{ (k )}\bigr\rrvert
_{n}^{2}:= \bigl\langle X^{
(k )},X^{ (k )} \bigr
\rangle_{n}.
\]
To make sure that the magnitude of the regression coefficients is
comparable, assume that $\llvert X^{ (k )}\rrvert _{n}^{2}=1$.
This is a standard condition that also simplifies the discussion throughout
(e.g., B{\"u}hlmann \cite{Buh06}, Barron \textit{et~al.}
\cite{Baretal08}). In practice, this
is achieved
by dividing the original variables by $\llvert X^{ (k
)}\rrvert _{n}$.
Throughout, it is assumed that the variables have unit $\llvert \cdot
\rrvert _{n}$
norm. This also implies that $\mathbb{E}\llvert X^{ (k
)}\rrvert _{n}^{2}=1$.
Denote by
\[
\mathcal{L} (B ):= \Biggl\{ \mu\dvt \mu(X )=\sum_{k=1}^{K}b_{k}X^{ (k )},
\sum_{k=1}^{K}\llvert b_{k}\rrvert
\leq B,X\in\mathcal{X} \Biggr\},
\]
the space of linear functions on $\mathcal{X}$ with $l_{1}$ coefficients
bounded by $B$. It follows that $\mathcal{L} (B )$ is a
Hilbert space under the inner product $ \langle X^{ (k
)},X^{ (l )} \rangle=\mathbb{E}X^{ (k
)}X^{ (l )}$
as well as the empirical inner product $ \langle X^{ (k
)},X^{ (l )} \rangle_{n}$.
Also, let $\mathcal{L}:=\bigcup_{B<\infty}\mathcal{L} (B )$
be the union of the above spaces for finite~$B$. The goal it to estimate
the regression function when the true expectation is replaced by the
empirical one, that is, when we use a finite sample of $n$ observations
$ \{ Y_{i},X_{i}\dvt i=1,2,\ldots,n \} $. As already mentioned,
$B$ is only known to be finite, and this is a standard set up used
elsewhere (e.g., B{\"u}hlmann \cite{Buh06}, Barron \textit{et~al.} \cite
{Baretal08}).
Moreover, $\mu_{0}$
does not need to be an element of $\mathcal{L} (B )$ for
any finite $B$.

Results are sometimes derived using some restricted eigenvalue condition
on the empirical Gram matrix of the regressors also called compatibility
condition (e.g., B{\"u}hlmann and van~de Geer \cite{Buhvan11}, for a list
and discussion).
For example, the minimum eigenvalue of the empirical Gram matrix of
any possible $m$ regressors out of the $K$ possible ones, is given
by
%
%
\begin{equation}
\rho_{m,n}:=\inf\Biggl\{ \frac{\llvert \sum_{k=1}^{K}X^{
(k )}b_{k}\rrvert _{n}^{2}}{\sum_{k=1}^{K}\llvert b_{k}\rrvert
^{2}}\dvt \sum
_{k=1}^{K} \{ b_{k}\neq0 \} =m \Biggr\},\label
{EQrestrictedEigenvalueCondition}
\end{equation}
where $ \{ b_{k}\neq0 \} $ is the indicator function of
a set (e.g., Zhang \cite{Zha09}, and many of the
references on Lasso cited
above; see also the isometry condition in Candes and Tao
\cite{CanTao07}). The
above condition means that the regressors are approximately orthogonal,
and typical examples are frames (e.g., Daubechies, Defrise and De Mol \cite{DauDefDeM04}). This
condition is usually avoided in the analysis of convergence rates
of greedy algorithms. Note that unless one uses a fixed design for
the regressors, (\ref{EQrestrictedEigenvalueCondition}) is random.
In this paper, $m$ usually refers to the number of iterations or
greedy steps at which the algorithm is stopped. The population counterpart
of (\ref{EQrestrictedEigenvalueCondition}) will be denoted by $\rho_{m}$,
that is,
%
%
\begin{equation}
\rho_{m}:=\inf\Biggl\{ \frac{\mathbb{E}\llvert \sum_{k=1}^{K}X^{ (k
)}b_{k}\rrvert _{n}^{2}}{\sum_{k=1}^{K}\llvert b_{k}\rrvert
^{2}}\dvt \sum
_{k=1}^{K} \{ b_{k}\neq0 \} =m \Biggr\}.\label
{EQeigenvaluePopulationRestriction}
\end{equation}
When $m$ is relatively small, $\rho_{m}$ plus an $\mathrm{o}_{p} (1 )$
term can be used to bound $\rho_{m,n}$ from below (e.g., Loh and
Wainwright \cite{LohWai12}; see also Nickl and van~de Geer \cite{Nicvan13}). Eigenvalue
restrictions will be avoided here under mixing dependent conditions.
However, under non-mixing and possibly long memory conditions, the
convergence rates can deteriorate quite quickly. Restricting attention
to the case $\rho_{m}>0$ allows one to derive more interesting results.

Throughout the paper, the following symbols are used: $\lesssim$
and $\gtrsim$ indicate inequality up to a multiplicative finite absolute
constant, $\asymp$ when the left-hand side and the right-hand side are
of the same
order of magnitude, $\wedge$ and $\vee$ are $\min$ and $\max$,
respectively, between the left-hand side and the right-hand side.

\subsection{Objective}

To ease notation, let $\llvert \cdot\rrvert _{2}= (\mathbb
{E}\llvert \cdot\rrvert ^{2} )^{1/2}$
and define
%
%
\begin{equation}
\gamma(B ):=\inf_{\mu\in\mathcal{L} (B
)}\llvert\mu-\mu_{0}\rrvert
_{2}\label{EQapproximationError}
\end{equation}
to be the approximation error of the best element in $\mathcal{L}
(B )$,
and let $\mu_{B}$ be the actual minimizer. Since for each $B<\infty$
the set $\mathcal{L} (B )$ is compact, one can replace the
$\inf$ with $\min$ in the above display. The approximation can improve
if $B$ increases. For simplicity, the notation does not make explicit
the dependence of the approximation error on $K$, as $K$ is the
same for all the algorithms, while $B$ can be different for the CGA
and FWA, as it will be shown in due course.

Let $X'$ be a random variable distributed like $X$ but independent
of the sample. Let $\mathbb{E}'$ be expectation w.r.t. $X'$ only.
The estimator from any of the greedy algorithms will be denoted by~$F_{m}$. The bounds are of the following kind:
%
%
\begin{equation}
\bigl(\mathbb{E}'\bigl\llvert\mu_{0}
\bigl(X' \bigr)-F_{m} \bigl(X' \bigr)\bigr
\rrvert^{2} \bigr)^{1/2}\lesssim \error (B,K,n,m )+\algo (B,m )+
\gamma(B )\label{EQerrorDecomposition}
\end{equation}
for any $B$ in some suitable range, where relates to the $B$ in
the approximation $\mu_{B}$ from (\ref{EQapproximationError}).
The possible values of $B$ depend on the algorithm. For the PGA,
OGA and RGA, $B<\infty$, that is, the algorithms allow to approximate
any function in $\mathcal{L}$, the union of $\mathcal{L}
(B )$
for any $B>0$. The CGA and FWA restrict $B\leq\bar{B}$ which is
a user specified parameter. This gives direct control of the estimation
error. The results for the CGA and FWA will be stated explicitly in
$\bar{B}$, so that $\bar{B}\rightarrow\infty$ is allowed. The term
$\gamma(B )$ is defined in (\ref{EQapproximationError}),
while
\[
\algo (B,m )^{2}\gtrsim\llvert Y-F_{m}\rrvert
_{n}^{2}-\inf_{\mu\in\mathcal{L} (B )}\llvert Y-\mu\rrvert
_{n}^{2}
\]
defines an upper bound for the error due to estimating using any of
the algorithms rather than performing a direct optimization. It could
be seen as part of the approximation error, but to clearly identify
the approximation properties of each algorithm, $\algo (B,m )$
is explicitly defined. Finally, the term $\error (B,K,n,m )$
is the estimation error.

\subsection{Approximation in function spaces}

When $\mu_{0}\notin\mathcal{L}$, the approximation can be large.
This is not to say that functions in $\mathcal{L}$ cannot represent
non-linear functions. For example, the set of regressors $\mathcal{X}$
could include functions that are dense in some set, or generally be
a subset of some dictionary (e.g., Mallat and Zhang \cite
{MalZha93}, Barron \textit{et~al.} \cite{Baretal08}, Sancetta \cite{San14}).

Consider the framework in Section~2.3 of Barron \textit{et~al.}
\cite{Baretal08}. Let
$\mu_{0}$ be a univariate function on $ [0,1 ]$, that is,
$\mu_{0}$
is the expectation of $Y$ conditional on a univariate variable with
values in $ [0,1 ]$. Suppose $\mathcal{D}$ is a dictionary
of functions on $ [0,1 ]$, and denote its elements by $g$.
Suppose that $\mu_{0}$ is in the closure of functions admitting the
representation\vspace*{1pt} $\mu(x )=\sum_{g\in\mathcal
{D}}b_{g}g (x )$,
where $\sum_{g\in\mathcal{D}}\llvert b_{g}\rrvert \leq B$; $b_{g}$
are coefficients that depend on the functions $g$. Examples include
sigmoid functions, polynomials, curvelet, frames, wavelets, trigonometric
polynomials, etc. Since $\mathcal{D}$ might be infinite or too large
for practical applications, one considers a subset $\mathcal{X}\subset
\mathcal{D}$,
which is a dictionary of $K$ functions on $ [0,1 ]$. Then
$\mu_{0} (x )=\sum_{g\in\mathcal{X}}b_{g}g (x
)+\sum_{g\in\mathcal{D}\setminus\mathcal{X}}b_{g}g (x )$.
Assuming that $\llvert \sum_{g\in\mathcal{D}\setminus\mathcal
{X}}b_{g}g (x )\rrvert _{2}\lesssim K^{-\alpha}$
for some $\alpha>0$, the approximation error decreases as one expands
the dictionary. Examples for non-orthogonal dictionaries are discussed
in Barron \textit{et~al.} \cite{Baretal08}. However, to aid
intuition, one can consider
Fourier basis for smooth enough functions to ensure that $\sum_{g\in
\mathcal{D}}\llvert b_{g}\rrvert <\infty$.
If $\mathcal{X}$ is large enough, one may expect the second summation
to have a marginal contribution.

Hence,\vspace*{1pt} with abuse of notation, the result of the present paper cover
the aforementioned problem, where the functions $g\in\mathcal{X}$
are then denoted by $ \{ x^{ (k )}\dvt k=1,2,\ldots,K \} $;
here $x\in[0,1 ]$, while each $g (x )$ is denoted
by $x^{ (k )} (x )$, so that $x^{ (k )}$
is not the $k$th entry in $x$ but a function of $x$ (the $k$th
element in a dictionary). As mentioned in the \hyperref[sec1]{Introduction}, this paper
does not make any distinction whether $\mathcal{X}$ is a set of explanatory
variables or functions (in general a dictionary), so it also covers
problems addressed in compress sensing with error noise.

\subsection{Conditions}

The theoretical properties of the algorithms are a function of the
dependence conditions used. At first, absolute regularity is used.
This allows to obtain results as good as if the data were independent
(e.g., Chen and Shen \cite{CheShe98}). However, for some
prediction problems,
absolute regularity might not be satisfied. Hence, more general dependence
conditions shall be used. Generality comes at a big cost in this case.

Some notation is needed to recall the definition of absolute regularity.
Suppose that $ (W_{i} )_{i\in\mathbb{Z}}$ is a stationary
sequence of random variables and, for any $d\geq0$, let $\sigma
(W_{i}\dvt i\leq0 )$,
$\sigma(W_{i}\dvt i\geq d )$ be the sigma algebra generated
by $ \{ W_{i}\dvt i\leq0 \} $ and $ \{ W_{i}\dvt i\geq d
\} $,
respectively. For any $d\geq0$, the beta mixing coefficient $\beta
(d )$
for $ (W_{i} )_{i\in\mathbb{Z}}$ is
\[
\beta(d ):=\mathbb{E}\sup_{A\in\sigma
(W_{i}\dvt i\geq d )}\bigl\llvert\Pr\bigl(A|
\sigma(W_{i}\dvt i\leq0 ) \bigr)-\Pr(A )\bigr\rrvert
\]
(see Rio \cite{Rio00}, Section~1.6, for other equivalent
definitions). The
sequence $ (W_{i} )_{i\in\mathbb{Z}}$ is absolutely regular
or beta mixing if $\beta(d )\rightarrow0$ for
$d\rightarrow\infty$.

Throughout, with slight abuse of notation, for any $p>0$, $\llvert
\cdot
\rrvert _{p}^{p}=\mathbb{E}\llvert \cdot\rrvert ^{p}$
is the $L_{p}$ norm (i.e., do not confuse $\llvert \cdot\rrvert _{n}$
with $\llvert \cdot\rrvert _{p}$). Moreover, $\mu_{0} (X
):=\mathbb{E} [Y|X ]$
is the true regression function, $Z:=Y-\mu_{0} (X )$ is
the error term, $\Delta(X )=\mu_{B} (X )-\mu
_{0} (X )$
is the approximation residual (recall that $\mu_{B}$ is the best
$L_{2}$ approximation to $\mu_{0}$ in $\mathcal{L} (B )$).

The asymptotics of the greedy algorithms are studied under the following
conditions.

%
%
\begin{condition}\label{ConditionEY|X} $\max_{k}\llvert X^{
(k )}\rrvert _{n}^{2}=1$,
$\max_{k}\llvert X^{ (k )}\rrvert _{2}=1$.
\end{condition}

%
%
\begin{condition}\label{ConditionabsoluteRegularityBounded}The\vspace*{1pt}
sequence $ (X_{i},Z_{i} )_{\in\mathbb{Z}}$ is stationary
absolutely regular with beta mixing coefficients $\beta(i
)\lesssim\beta^{i}$
for some $\beta\in[0,1)$ and $\mathbb{E}\llvert Z\rrvert ^{p}<\infty$
for some $p>2$, $\max_{k\leq K}\llvert X^{ (k )}\rrvert $
is bounded, and the approximation residual $\Delta(X )=\mu
_{B} (X )-\mu_{0} (X )$
is also bounded. Moreover, $1<K\lesssim\exp\{ Cn^{a} \} $,
for some absolute constant $C$ and $a\in[0,1)$.
\end{condition}

Bounded regressors and sub-Gaussian errors are the common conditions
under which greedy algorithms are studied. Condition~\ref{ConditionabsoluteRegularityBounded}
already weakens this to the error terms only possessing a $p>2$ moment.
However, restricting attention to bounded regressors can be limiting.
The next condition replaces this with a moment condition.

%
%
\begin{condition}\label{ConditionabsoluteRegularity}The sequence
$ (X_{i},Z_{i} )_{\in\mathbb{Z}}$ is stationary absolutely
regular with beta mixing coefficients $\beta(i )\lesssim
\beta^{i}$
for some $\beta\in[0,1)$ and
%
%
\begin{equation}
\mathbb{E}\bigl\llvert ZX^{ (k )}\bigr\rrvert^{p}+\mathbb{E}
\bigl\llvert X^{ (k )}\bigr\rrvert^{2p}+\mathbb{E}\bigl\llvert
\Delta(X )X^{ (k )}\bigr\rrvert^{p}<\infty,\label{EQboundVariance}
\end{equation}
for some $p>2$. Moreover, $1<K\lesssim n^{\alpha}$ for some $\alpha
< (p-2 )/2$
(with $p$ as just defined).
\end{condition}

Note that in the case of independent random variables, one could relax
the moment condition to $p\geq2$. Recall that $\mu_{0}$ is not restricted
to be in $\mathcal{L} (B )$. Only the resulting estimator
will be. The expectation of $\Delta(X )$ is the bias.

There are examples of models that are not mixing (e.g., Andrews \cite
{And84}, Bradley \cite{Bra86}). For
example, the sieve bootstrap is not mixing
(Bickel and B{\"u}hlmann \cite{BicBuh99}). It is important to
extend the applicability of
the algorithms to such case. The gain in generality leads to a considerably
slower rate of convergence than the i.i.d. and beta mixing case. This
is mostly due to the method of proof. It is not known whether the
results can be improved in such cases. Dependence is now formalized
by the following.

%
%
\begin{condition}\label{Conditiondependence}Denote by $\mathbb{E}_{0}$
the expectation conditional at time $0$ (w.r.t. the natural filtration
of the random variables). Recall that $\llvert \cdot\rrvert _{p}:=
(\mathbb{E}\llvert \cdot\rrvert ^{p} )^{1/p}$.
The sequence $ (X_{i},Z_{i} )_{\in\mathbb{Z}}$ is stationary,
and for some $p\geq2$,
\[
d_{n,p}:=\max_{k}\sum
_{i=0}^{n}\frac{ (\llvert \mathbb
{E}_{0}Z_{i}X_{i}^{ (k )}\rrvert _{p}+\llvert \mathbb
{E}_{0} [ (1-\mathbb{E} )\llvert X_{i}^{ (k
)}\rrvert ^{2} ]\rrvert _{p}+\llvert \mathbb{E}_{0} [
(1-\mathbb{E} )\Delta(X_{i} )X_{i}^{ (k
)} ]\rrvert _{p} )}{ (i+1 )^{1/2}}<\infty
\]
for any $n$.
\end{condition}

Note that the dependence condition is in terms of mixingales and for
weakly dependent data, $\sup_{n}d_{n,p}<\infty$ when the $p$th
moment exists, under certain conditions. The general framework allows
us to consider data that might be strongly dependent (long memory),
when $d_{n,p}\rightarrow\infty$ (see Example~\ref{ExamplelongMemory}
for some details).

\section{Algorithms}\label{SectionAlgorithms}

The algorithms have already appeared elsewhere, and they will be reviewed
in Section~\ref{SectionreviewAlgorithms}. All the algorithms studied
here do achieve a global minimum of the empirical risk. This minimum
might not be unique if the number of variables are larger than the
sample size. Moreover, the convergence rates of the algorithms to
the global minimum can differ. The reader unfamiliar with them, can
skim through Section~\ref{SectionreviewAlgorithms} before reading
the following. In particular, the PGA has the slowest rate, while
all the others have a faster rate which is essentially optimal (see
Lemmas~\ref{LemmaL2Boosting},~\ref{LemmaOGA},~\ref{Lemmarelaxedgreedy}
and~\ref{LemmafrankWolfeApproximation}, for the exact rates used
here; see DeVore and Temlyakov \cite{DeVTem96}, and Barron \textit
{et~al.} \cite{Baretal08}, for
discussions on optimality of convergence rates). The optimal rate
toward the global minimum is $m^{-1/2}$ under the square root of
the empirical square error loss, where $m$ is the number of greedy
iterations. For the PGA the convergence rate of the approximation
error of the algorithm, $\algo (B,m )$, is only $m^{-1/6}$,
without requiring the target $Y$ to be itself an element of $\mathcal
{L} (B )$,
that is, a linear function with no noise (Lemma~\ref{LemmaL2Boosting}).
For functions in $\mathcal{L} (B )$, Konyagin and Temlyakov \cite{KonTem99} improved the rate to $m^{-11/62}$, while Livshitz and Temlyakov
\cite{LivTem03}
show a lower bound $m^{-0.27}$. Hence, the approximation rate
of the PGA is an open question. The slow rate of the PGA ($L_{2}$-Boosting)
has led Barron \textit{et~al.} \cite{Baretal08} to disregard it.
While the approximating
properties of the PGA are worse than the other algorithms, its finite
sample properties tend to be particularly good in many cases (e.g.,
Section~\ref{SectionnumericalResults}). An overview of how the present
results add to the literature and further details are summarized next.

\subsection{Comparison with existing results}\label{SectionsummaryAlgorithms}

\begin{table}
\tabcolsep=0pt
\tablewidth=250pt
\caption{Comparison of results}\label{tab1}
\begin{tabular*}{\tablewidth}{@{\extracolsep{\fill}}@{}ll@{}}
\hline
Algorithm/author/conditions & Rates\\
\hline
PGA & \\
\quad $\mathrm{M}(X;b)$, $\mathrm{M}(Z;p)$, $\mathrm{D}(X,Z;\beta^{n})$, $\mathrm{K}(E)$, $L_{2}$ & $(\frac{\ln K}{n} )^{1/8}$\\[3pt]
\quad $\mathrm{M}(X,Z;g)$, $\mathrm{D}(X,Z;\mathrm{NM})$, $\mathrm{K}(P)$, X, $L_{2}$ & $ (\frac{d_{n,\bar{p}}^{2}}{n} )^{ (1-\epsilon)/8}$\\[6pt]
B{\"u}hlmann and van~de Geer \cite{Buhvan11} & \\[3pt]
\quad $\mathrm{M}(X;b)$, $\mathrm{M}(Z;g)$, $\mathrm{D}(X,Z;\mathit{iid})$, $\mathrm{K}(E)$, $L_{2n}$ & $ (\frac{\ln K}{n} )^{ (1-\epsilon)/16}$\\[3pt]
Lutz and B{\"u}hlmann \cite{LutBuh06} & \\[3pt]
\quad $\mathrm{M}(X,Z;p)$, $\mathrm{D}(X,Z;n^{\alpha})$, $\mathrm{K}(E)$, $L_{2}$ & $\mathrm{o} (1)$\\[9pt]
OGA & \\[3pt]
\quad $\mathrm{M}(X;b)$, $\mathrm{M}(Z;p)$, $\mathrm{D}(X,Z;\beta^{n})$, $\mathrm{K}(E)$, $L_{2}$ & $(\frac{\ln K}{n} )^{1/4}$\\[3pt]
\quad $\mathrm{M}(X,Z;g)$, $\mathrm{D}(X,Z;\mathrm{NM})$, $\mathrm{K}(P)$, X, $L_{2}$ & $ (\frac{d_{n,\bar{p}}^{2}}{n} )^{ (1-\epsilon)/6}$\\[6pt]
B{\"u}hlmann and van~de Geer \cite{Buhvan11} & \\[3pt]
\quad $\mathrm{M}(X;b)$, $\mathrm{M}(Z;g)$, $\mathrm{D}(X,Z;\mathit{iid})$, $\mathrm{K}(E)$, $L_{2}$ & $ (\frac{1}{n} )^{1/6}\vee(\frac{\ln K}{n} )^{1/4}$\\[3pt]
Barron \textit{et~al.} \cite{Baretal08} & \\[3pt]
\quad $\mathrm{M}(X,Z;b)$, $\mathrm{D}(X,Z;\mathit{iid})$, $\mathrm{K}(P)$, $EL_{2}$ & $ (\frac{\ln K}{n} )^{1/4}$\\[3pt]
Zhang \cite{Zha09} & \\[3pt]
\quad $\mathrm{M}(X,Z;b)$, $\mathrm{D}(X,Z;\mathit{iid})$, X, $L_{2n}$, $+$ & $ (\frac{K_{0}}{n} )^{1/2}$ \\[9pt]
RGA & \\[3pt]
\quad $\mathrm{M}(X;b)$, $\mathrm{M}(Z;p)$, $\mathrm{D}(X,Z;\beta^{n})$, $\mathrm{K}(E)$, $L_{2}$ & $(\frac{\ln K}{n} )^{1/4}$\\[3pt]
\quad $\mathrm{M}(X,Z;g)$, $\mathrm{D}(X,Z;\mathrm{NM})$, $\mathrm{K}(P)$, X, $L_{2}$ & $ (\frac{d_{n,\bar{p}}^{2}}{n} )^{ (1-\epsilon)/6}$\\[3pt]
Barron \textit{et~al.} \cite{Baretal08} & \\[3pt]
\quad $\mathrm{M}(X,Z;b)$, $\mathrm{D}(X,Z;\mathit{iid})$, $\mathrm{K}(P)$, $EL_{2}$ & $ (\frac{\ln K}{n} )^{1/4}$\\[9pt]
CGA and FWA & \\[3pt]
\quad $\mathrm{M}(X;b)$, $\mathrm{M}(Z;p)$, $\mathrm{D}(X,Z;\beta^{n})$, $\mathrm{K}(E)$, $L_{2}$, $+$ & $(\frac{\ln K}{n} )^{1/4}$\\[3pt]
\quad $\mathrm{M}(X,Z;g)$, $\mathrm{D}(X,Z;\mathrm{NM})$, $\mathrm{K}(P)$, $L_{2}$, $+$ & $ (\frac{d_{n,\bar{p}}^{2}}{n} )^{ (1-\epsilon)/4}$\\
\hline
\end{tabular*}\vspace*{-2pt}
\end{table}

There are many results on greedy algorithms under different conditions.
Table~\ref{tab1} summarizes and compares some of these results. For each algorithm
the most interesting results from the present paper are presented
first. The symbols used to describe the conditions are defined in
the glossary of symbols at the end of this section.

\subsubsection{Glossary for Table~\texorpdfstring{\protect\ref{tab1}}{1}}

\paragraph{Moments}

M (variable; moment type); moment types: $p={}$moments, refer to paper
for exact $p$, $g={}$sub-Gaussian tails, $b={}$bounded random variables;
for example, $\mathrm{M}(X,Z;b)$ means that both $X$ and $Z$ are bounded.

\paragraph{Dependence}

D (variable, dependence type); dependence types are all stationary:
$\mathit{iid}={}$i.i.d. or just independence, $\alpha^{n}$/$\beta^{n}={}$geometric
alpha/beta mixing, $n^{\alpha}$/$n^{\beta}={}$polynomial alpha/beta
mixing; $\mathrm{NM}={}$non-mixing; see paper for details on the polynomial
rate and how it relates to moments.

\paragraph{K}

K (growth rate); number of regressors $K$: $P=n^{a}$ for any $a<\infty$,
$E=\exp\{ Cn^{a} \} $ for $a\in[0,1)$, $C<\infty$.

\paragraph{Design matrix}

$X$ if conditions are imposed on the design matrix, for example, compatibility
conditions, otherwise, no symbol is reported.

\paragraph{Loss function}

$L_{2}=L_{2}$ loss as in the l.h.s. of (\ref{EQerrorDecomposition})
and results holding in probability, $EL_{2}={}$same as $L_{2}$ but
results holding in $L_{1}$ (i.e., take a second expectation w.r.t.
to the sample data), $L_{2n}={}$empirical $L_{2}$ loss.

\paragraph{Additional remarks on glossary}

The true function $\mu_{0}$ is assumed to be in $\mathcal{L}
(B )$
for some finite $B$. When rates are included, $\epsilon$ is understood
to be a positive arbitrarily small constant. Also, $\bar{p}$ in
$d_{n,\bar{p}}$
refers to a large $p$ depending on $\epsilon$ and $K$, with exact
details given in Corollary~\ref{CorollarydependenceFiniteMoments}.
In some cases, conditions may not fit exactly within the classification
given above due to minor differences, in which case they may still
be classified within one group. The symbol $+$ is used to denote additional
conditions which can be found in the cited paper. For Zhang \cite{Zha09},
$K_{0}$ represents the true number of non-zero coefficients and it
is supposed to be small. For the CGA and the FWA, the symbol${}+{}$refers
to the fact that the user pre-specifies a $\bar{B}<\infty$ and constrains
estimation in $\mathcal{L} (B )$ with $B\leq\bar{B}$, and
it also assumes that $\mu_{0}\in\mathcal{L} (\bar{B} )$.
The results in the paper are more general, and the restrictions in
Table~\ref{tab1} are for the sake of concise exposition and comparison.

\subsubsection{Comments}

Table~\ref{tab1} only provides upper bounds. Interest would also lie in deriving
lower bound estimates (e.g., Donoho and Johnstone \cite
{DonJoh98}, Birg\'e and Massart \cite{BirMas01}, Tsybakov \cite{Tsy03}, and Bunea, Tsybakov and Wegkamp \cite{BunTsyWeg07N2}, for such rates
for certain nonparametric parametric problems; see also Tsybakov \cite{Tsy09}, Chapter~2, for a general discussion on lower bounds). The results
in Tsybakov \cite{Tsy03} and Bunea, Tsybakov and Wegkamp \cite{BunTsyWeg07N2}
provide minimax rates
and explicit estimators for certain function classes which exactly
apply in the present context. Suppose that the error term $Z$ is
Gaussian, the regressors $X$ are bounded and an i.i.d. sample is
available. Let $\mu_{n}$ be any estimator in $\mathcal{L}
(B )$.
From Theorem~2 in Tsybakov \cite{Tsy03}, one can deduce that
\[
\sup_{\mu\in\mathcal{L} (B )}\llvert\mu-\mu_{n}\rrvert\gtrsim\cases{
\displaystyle B\sqrt{\frac{K}{n}}, &\quad if $K\lesssim\sqrt{n}$,
\vspace*{5pt}\cr
\displaystyle B \biggl(\frac{\ln K}{n} \biggr)^{1/4}, &\quad if $K
\gtrsim\sqrt{n}$.}
\]
This results is also useful to understand the difference between the
result derived by Zhang \cite{Zha09} for the OGA and
usual results for Lasso
under sparsity. In these cases, the target function is in a much smaller
class than $\mathcal{L} (B )$, that is, $\mu_{0}$ is a linear
function with a small number of $K_{0}$ non-zero regression coefficients.
Within this context, one can infer that the result from Zhang \cite{Zha09}
is the best possible (e.g., use Theorem~3 in Tsybakov
\cite{Tsy03}).

Under mixing conditions, the convergence rates for the OGA, RGA, CGA
and FWA are optimal. Table~\ref{tab1} shows that the results in Barron \textit
{et~al.} \cite{Baretal08} for the OGA and RGA are also optimal,
but require i.i.d. bounded
regressors and noise. The convergence rates for the PGA are not optimal,
but considerably improve the ones of B{\"u}hlmann and van~de Geer \cite{Buhvan11}
also allowing for unbounded regressors and dependence.

\subsection{Statement of results}

\subsubsection{Mixing data}

In the following, when some relation is said to hold in probability,
it means it holds with probability going to one as $n\rightarrow\infty$.
Also, note that the linear projection of $\mu_{0}$ onto the space
spanned by the regressors is in $\mathcal{L}$ (the union of the
$\mathcal{L} (B )$
spaces) because the number of regressors $K$ is finite. Hence, let
%
%
\begin{equation}
B_{0}:=\arg\inf_{B>0}\gamma(B )\label{EQOLSB}
\end{equation}
be the absolute sum of the coefficients in the unconstrained linear
projection of $\mu_{0}$ onto the space spanned by the regressors
($\gamma(B )$ as in (\ref{EQapproximationError})). Of
course, $K$ is allowed to diverge to infinity with $n$, if needed,
which in consequence may also imply $B_{0}$ in (\ref{EQOLSB}) can
go to infinity.

%
%
\begin{theorem}\label{TheoremabsoluteRegularity}Under Condition
\ref{ConditionEY|X} and either Conditions~\ref{ConditionabsoluteRegularityBounded}
or~\ref{ConditionabsoluteRegularity},
%
%
\begin{equation}
\bigl(\mathbb{E}'\bigl\llvert\mu_{0}
\bigl(X' \bigr)-F_{m} \bigl(X' \bigr)\bigr
\rrvert^{2} \bigr)^{1/2}\lesssim \error (B,K,n,m )+\algo (B )+
\gamma(B )\label{EQmainObjective}
\end{equation}
in probability, where
%
%
\begin{equation}
B\in\cases{ [B_{0},\infty), &\quad for the PGA, OGA, RGA,
\vspace*{3pt}\cr
(0,
\bar{B}], &\quad for the CGA and FWA,}\label{EQsizeTrueParameter}
\end{equation}
where\vspace*{-4pt}
%
%
\begin{eqnarray}
\error (B,K,n,m ) &=&\cases{ \displaystyle\sqrt{\frac{m\ln K}{n}}, &\quad
for the
PGA, OGA, RGA,
\cr
\displaystyle\bar{B} \biggl(\frac{\ln K}{n}
\biggr)^{1/4}, &\quad for the CGA and FWA,} \label{EQerrorAbsoluteReg}
\\
\algo (B,m ) &=& \cases{ B^{1/3}m^{-1/6}, &\quad for the PGA,
\vspace*{3pt}\cr
Bm^{-1/2}, &\quad for the OGA and RGA,
\vspace*{3pt}\cr
\bar{B}m^{-1/2}, &\quad
for the CGA and FWA.}\label{EQalgoError}
\end{eqnarray}
\end{theorem}

%
%
\begin{remark}When $B_{0}\leq\bar{B}$, asymptotically, the CGA and
FWA impose no constraint on the regression coefficients. In this case,
these algorithms also satisfy (\ref{EQmainObjective}) with (\ref
{EQsizeTrueParameter})
as for the OGA and RGA. While $B_{0}$ is unknown, this observation
will be used to deduce Corollary~\ref{Corollarypersistence}. Also
note that~(\ref{EQmainObjective}) for the PGA, OGA and RGA is minimized
by $B=B_{0}$.
\end{remark}

Theorem~\ref{TheoremabsoluteRegularity} allows one to answer several
questions of interest about the algorithms. Note that $\error (B,K,n,m )$
in (\ref{EQerrorAbsoluteReg}) does not depend on $B$, as a consequence
of the method of proof; it will depend on $B$ for some of the other
results. The next two results will focus on two related important
questions. One concerns the overall convergence rates of the estimator
when the true function $\mu_{0}\in\mathcal{L}$, that is, $\mu_{0}$ is
linear with absolutely summable coefficients. The other concerns the
largest linear model in reference of which the estimator is optimal
in a square error sense (i.e., persistence in the terminology of
Greenshtein and Ritov \cite{vanvan04}, or traditionally, this is
termed consistency for
the linear pseudo true value). Rates of convergence are next. These
rates directly follow from Theorem~\ref{TheoremabsoluteRegularity},
using the fact that $B<\infty$ and equating $\error (B,K,n,m )$
with $\algo (B,m )$ and solving for $m$.

%
%
\begin{corollary}Under the conditions of Theorem~\ref{TheoremabsoluteRegularity},
if
\[
m\mbox{ satisfies }\cases{ \displaystyle\asymp\biggl(\frac{n}{\ln K}
\biggr)^{3/4}, &\quad for the PGA,
\vspace*{3pt}\cr
\displaystyle\asymp\sqrt{
\frac{n}{\ln K}}, &\quad for the OGA and RGA,
\vspace*{3pt}\cr
\displaystyle\gtrsim\sqrt{
\frac{n}{\ln K}}, &\quad for the CGA and FWA}
\]
then, in probability
\[
\bigl(\mathbb{E}'\bigl\llvert\mu_{0}
\bigl(X' \bigr)-F_{m} \bigl(X' \bigr)\bigr
\rrvert^{2} \bigr)^{1/2}\lesssim\cases{ \displaystyle\biggl(
\frac{\ln K}{n} \biggr)^{1/8}, &\quad for the PGA if $
\mu_{0}\in\mathcal{L}$,
\vspace*{3pt}\cr
\displaystyle\biggl(\frac{\ln K}{n}
\biggr)^{1/4}, &\quad for the OGA and RGA if $\mu_{0}\in
\mathcal{L}$,
\vspace*{3pt}\cr
\displaystyle\bar{B} \biggl(\frac{\ln K}{n}
\biggr)^{1/4}, &\quad for the CGA and FWA if $\mu_{0}\in
\mathcal{L} (\bar{B} )$.}
\]
\end{corollary}

The CGA and FWA achieve the minimax rate under either Conditions~\ref{ConditionabsoluteRegularityBounded}
or~\ref{ConditionabsoluteRegularity} if $\mu_{0}\in\mathcal{L}
(\bar{B} )$
as long as the number of iterations $m$ is large enough. The drawback
in fixing $\bar{B}$ is that if $\mu_{0}\in\mathcal{L} (B )$
with $\bar{B}<B$, there can be an increase in bias. This can be avoided
by letting $\bar{B}\rightarrow\infty$ with the sample size. The following
can then be used to bound the error when $\bar{B}<B$, if $\mu_{0}\in
\mathcal{L} (B )$
(Sancetta \cite{San14}).

%
%
\begin{lemma}\label{LemmaBbarapproximation}Let $\mu\in\mathcal
{L} (B )$
for some $B<\infty$. Then
\[
\inf_{\mu'\in\mathcal{L} (B' )}\bigl\llvert\mu-\mu'\bigr\rrvert
_{2}\leq\max\bigl\{ B-B',0 \bigr\}.
\]
\end{lemma}

The bounds are explicit in $\bar{B}$ so that one can let $\bar
{B}\rightarrow\infty$
if needed and apply Lemma~\ref{LemmaBbarapproximation} to show
that the approximation error goes to zero if $\mu_{0}\in\mathcal
{L} (B )$
for some bounded $B$.

Next, one can look at the idea of persistence, which is also related
to consistency of an estimator for the pseudo true value in the class
of linear functions. Adapting the definition of persistence to the
set up of this paper, the estimator $F_{m}$ is persistent at the
rate $B\rightarrow\infty$ if
%
%
\begin{equation}
\mathbb{E}'\bigl\llvert Y'-F_{m}
\bigl(X' \bigr)\bigr\rrvert^{2}-\inf_{\mu\in
\mathcal{L} (B )}
\mathbb{E}'\bigl\llvert Y'-\mu\bigl(X'
\bigr)\bigr\rrvert^{2}=\mathrm{o}_{p} (1 ),\label{EQpersistence}
\end{equation}
where $X'$ and $Y'$ are defined to have same marginal distribution
as the $X_{i}$'s and $Y_{i}$'s, but independent of them. Directly
from Theorem~\ref{TheoremabsoluteRegularity} deduce the following.

%
%
\begin{corollary}\label{Corollarypersistence}Let $\bar{B}=B$ for
the CGA and FWA. Under the conditions of Theorem~\ref{TheoremabsoluteRegularity},
(\ref{EQpersistence}) holds if $m\rightarrow\infty$ such that
$m=\mathrm{o} (n/\ln K )$
and $B=\mathrm{o} (\sqrt{m} )$ for all algorithms.
\end{corollary}

\subsubsection{Non-mixing and strongly dependent data}

In the non-mixing case, the rates of convergence of the estimation
error can quickly deteriorate. Improvements can then be obtained by
restricting the population Gram matrix of the regressors to be full
rank. The next result does not restrict $\rho_{m}$.

%
%
\begin{theorem}\label{TheoremDependenceCrude}Under Conditions~\ref{ConditionEY|X}
and~\ref{Conditiondependence}, (\ref{EQmainObjective}) holds in
probability, with
\[
\error (B,K,n,m )= \biggl(\frac{d_{n,p}^{2}K^{4/p}}{n} \biggr
)^{1/4}\times\cases{
\bigl(B+m^{1/2} \bigr), &\quad for the PGA,
\vspace*{3pt}\cr
(B+m ), &\quad for the
OGA, RGA,
\vspace*{3pt}\cr
\bar{B}, &\quad for the CGA and FWA}
\]
and
\[
B\in\cases{ (0,\infty), &\quad for the PGA, OGA, RGA,
\vspace*{3pt}\cr
(0,\bar{B}], &\quad
for the CGA and FWA}
\]
and $\algo (B,m )$ as in (\ref{EQalgoError}).
\end{theorem}

Unlike $\error (B,K,n,m )$ in (\ref{EQerrorAbsoluteReg})
which did not depend on $B$, the above is derived using a different
method of proof and does depend on $B$. Also note the different restriction
on $B$. Letting $\mu_{0}\in\mathcal{L}$, one obtains the following
explicit convergence rates.

%
%
\begin{corollary}Suppose that
%
%
\begin{equation}
m\mbox{ satisfies }\cases{ \asymp \displaystyle\biggl(
\frac{n}{d_{n,p}^{2}K^{4/p}} \biggr)^{3/8}, &\quad for the PGA,
\vspace*{3pt}\cr
\asymp \displaystyle\biggl(\frac
{n}{d_{n,p}^{2}K^{4/p}} \biggr)^{1/12}, &\quad for the OGA and
RGA,
\vspace*{3pt}\cr
\gtrsim \displaystyle\biggl(\frac
{n}{d_{n,p}^{2}K^{4/p}}
\biggr)^{1/8}, &\quad for the CGA and FWA.}\label
{EQoptimalmnoMixingnoEigenvalue}
\end{equation}
Under the conditions of Theorem~\ref{TheoremDependenceCrude}, in
probability,
\[
\bigl(\mathbb{E}'\bigl\llvert\mu_{0}
\bigl(X' \bigr)-F_{m} \bigl(X' \bigr)\bigr
\rrvert^{2} \bigr)^{1/2}\lesssim\cases{ \displaystyle\biggl(
\frac{d_{n,p}^{2}K^{4/p}}{n} \biggr)^{1/16}, &\quad for the PGA if $
\mu_{0}\in\mathcal{L}$,
\vspace*{3pt}\cr
\displaystyle\biggl(\frac{d_{n,p}^{2}K^{4/p}}{n}
\biggr)^{1/12}, &\quad for the OGA and RGA if $\mu_{0}\in
\mathcal{L}$,
\vspace*{3pt}\cr
\displaystyle\biggl(\frac{d_{n,p}^{2}K^{4/p}}{n} \biggr)^{1/4}, &
\quad for the CGA and FWA if $\mu_{0}\in\mathcal{L} (\bar{B} )$.}
\]
\end{corollary}

The results are close to the lower bound $\mathrm{O} (n^{-1/4} )$
only for the CGA and FWA under weak dependence, not necessarily mixing
data (i.e., $\sup_{n}d_{n,p}<\infty$) and variables with moments of
all orders (i.e., $p$ arbitrary large). Now, restrict attention to
$\rho_{K}>0$, that is, $\rho_{m}$ in (\ref
{EQeigenvaluePopulationRestriction})
with $m=K$. This is equivalent to say that the population Gram matrix
of the regressors has full rank. In this case, the results for the
PGA, OGA and RGA can be improved. By following the proofs in
Section~\ref{SectionProofs}, it is easy to consider $\rho_{m}$ going to
zero as $m\rightarrow\infty$, but at the cost of extra details, hence
these case will not be reported here.

%
%
\begin{theorem}\label{TheoremDependenceEigenvalue}Suppose that
$\rho_{K}>0$. Under Conditions~\ref{ConditionEY|X} and~\ref{Conditiondependence},
for the PGA, OGA and RGA, (\ref{EQmainObjective})~holds in probability
with
\[
\error (B,K,n,m )= \bigl(m+m^{1/2}B \bigr) \biggl(\frac
{d_{n,p}^{2}K^{4/p}}{n}
\biggr)^{1/2}
\]
for any positive $B$, and $\algo (B,m )$ as in (\ref{EQalgoError}),
as long as $\error (B,K,n,m )+\algo (B,m )=\mathrm{o}
(1 )$.
\end{theorem}

The above theorem leads to much better convergence rates.

%
%
\begin{corollary}Suppose that
%
%
\begin{equation}
m\asymp\cases{ \displaystyle\biggl(\frac{n}{d_{n,p}^{2}K^{4/p}} \biggr
)^{3/4}, &
\quad for the PGA,
\vspace*{3pt}\cr
\displaystyle\biggl(\frac{n}{d_{n,p}^{2}K^{4/p}}
\biggr)^{1/3}, &\quad for the OGA and RGA.} \label{EQoptimalmNonMixing}
\end{equation}
Under the conditions of Theorem~\ref{TheoremDependenceEigenvalue},
\[
\bigl(\mathbb{E}'\bigl\llvert\mu_{0}
\bigl(X' \bigr)-F_{m} \bigl(X' \bigr)\bigr
\rrvert^{2} \bigr)^{1/2}\lesssim\cases{ \displaystyle\biggl(
\frac{d_{n,p}^{2}K^{4/p}}{n} \biggr)^{1/8}, &\quad for the PGA if $
\mu_{0}\in\mathcal{L}$,
\vspace*{3pt}\cr
\displaystyle\biggl(\frac{d_{n,p}^{2}K^{4/p}}{n}
\biggr)^{1/6}, &\quad for the OGA and RGA if $\mu_{0}\in
\mathcal{L}$.}
\]
\end{corollary}

Under non-mixing dependence, deterioration in the rate of convergence
due to $K$ becomes polynomial rather than the logarithmic one of
Theorem~\ref{TheoremabsoluteRegularity}. On the positive side, the
dependence condition used is very simple and can be checked for many
models (e.g., Doukhan and Louhichi \cite{DouLou99}, Section~3.5,
Dedecker and Doukhan \cite{DedDou03}, for examples and
calculations). Interesting results
can be deduced when the regressors have a moment generating function.
Then the rates of convergence can be almost as good if not better
than the ones derived by other authors assuming i.i.d. data, though
only when $\rho_{K}>0$ holds.

%
%
\begin{corollary}\label{CorollarydependenceFiniteMoments}Suppose
that $X$ and $Z$ have moments of all order and $K\lesssim n^{\alpha}$
for some $\alpha\in\mathbb{N}$. Under Conditions~\ref{ConditionEY|X}
and~\ref{Conditiondependence}, choosing $m$ as in (\ref{EQoptimalmNonMixing})
for the PGA, OGA and RGA and as in (\ref{EQoptimalmnoMixingnoEigenvalue})
for the CGA and FWA, for any $\epsilon\in(0,1 )$, and
$p=4\alpha/\epsilon$,
\begin{eqnarray*}
&& \bigl(\mathbb{E}'\bigl\llvert\mu_{0}
\bigl(X' \bigr)-F_{m} \bigl(X' \bigr)\bigr
\rrvert^{2} \bigr)^{1/2}
\\
&&\quad \lesssim\cases{ \displaystyle\biggl(
\frac{d_{n,p}^{2}}{n} \biggr)^{ (1-\epsilon
)/8}, &\quad for the PGA if $
\mu_{0}\in\mathcal{L}$ and $\rho_{K}>0$,
\vspace*{3pt}\cr
\displaystyle
\biggl(\frac{d_{n,p}^{2}}{n} \biggr)^{ (1-\epsilon
)/6}, &\quad for the OGA and RGA if $
\mu_{0}\in\mathcal{L}$ and $\rho_{K}>0$,
\vspace*{3pt}\cr
\displaystyle
\biggl(\frac{d_{n,p}^{2}}{n} \biggr)^{ (1-\epsilon
)/4}, &\quad for the CGA and FWA if $
\mu_{0}\in\mathcal{L} (\bar{B} )$}
\end{eqnarray*}
in probability.
\end{corollary}

\subsection{Review of the algorithms}\label{SectionreviewAlgorithms}

The algorithms have been described in several places in the literature.
The following sections review them. The first two algorithms are boosting
algorithms and they are reviewed in B{\"u}hlmann and van~de Geer \cite{Buhvan11}.
The third algorithm has received less attention in statistics despite
the fact that it has desirable asymptotic properties (Barron \textit
{et~al.} \cite{Baretal08}). The fourth algorithm is a
constrained version of the third
one and further improves on it in certain cases. The fifth and last
algorithm is the basic version of the Frank--Wolfe \cite{FraWol56} algorithm.

\begin{figure}[b]
\fbox{
\begin{tabular}{l}
Set: \\
$m\in\mathbb{N}$\\
$F_{0}:=0$\\
$\nu\in(0,1]$\\
For: $j=1,2,\ldots,m$\\
$s (j ):=\arg\max_{k}\llvert \langle
Y-F_{j-1},X^{ (k )} \rangle_{n}\rrvert $\\[2pt]
$g_{j} (X ):= \langle Y-F_{j-1},X^{s (j)} \rangle_{n}X^{s (j )}$\\
$F_{j}:=F_{j-1}+\nu g_{j} (X )$
\end{tabular}}
\caption{PGA ($L_{2}$-Boosting).}\vspace*{-6pt}\label{fig1}
\end{figure}

\subsubsection{Pure Greedy Algorithm (a.k.a. $L_{2}$-Boosting)}

Boosting using the $L_{2}$ norm is usually called $L_{2}$-Boosting,
though some authors also call it Pure Greedy Algorithm (PGA) in order
to stress its origin in the approximation theory literature (e.g.,
Barron \textit{et~al.} \cite{Baretal08}), and this is how it will
be called here. The term
matching pursuit is also used by engineers (e.g., Mallat and Zhang
\cite{MalZha93}). Figure~\ref{fig1} recalls the algorithm. There, $\nu\in
(0,1]$ is the
shrinkage parameter and it controls the degree of greediness in the
algorithm. For example, as $\nu\rightarrow0$ the algorithm in
Figure~\ref{fig1} converges to Stagewise Linear Regression, a variant of the LARS
algorithm that has striking resemblance to Lasso (Efron \textit{et~al.}
\cite{Efretal04},
for details). In order to avoid ruling out good regressors that are
correlated to $X^{s (m )}$ ($s (m )$ as defined
in Figure~\ref{fig1} and $X^{s (m )}=X^{ (s (m
) )}$
throughout to ease notation) one chooses $\nu$ smaller then $1$,
usually $0.1$ (B{\"u}hlmann \cite{Buh06}).

\begin{figure}
\fbox{
\begin{tabular}{l}
Set: \\
$m\in\mathbb{N}$\\
$F_{0}:=0$\\
For: $j=1,2,\ldots,m$\\
$s (j ):=\arg\max_{k}\llvert \langle
Y-F_{j-1},X^{ (k )} \rangle_{n}\rrvert $\\[2pt]
$P_{X}^{j}:={}$OLS operator on $\operatorname{span} \{ X^{s (1)},X^{s (2 )},\ldots,X^{s (j )} \} $ \\[2pt]
$F_{j}:=P_{X}^{j}Y$
\end{tabular}}
\caption{OGA (Orthogonal Matching Pursuit).}\vspace*{-8pt}\label{fig2}
\end{figure}

The PGA recursively fits the residuals from the previous regression
to the univariate regressor that reduces the most the residual sum
of the squares. At each step $j$, the algorithm solves $\min
_{k,b}\llvert Y-F_{j-1}-X^{ (k )}b\rrvert _{n}^{2}$.
However, the coefficient can then be shrunk by an amount $\nu\in
(0,1 )$
in order to reduce the degree of greediness. The resulting function
$F_{m}$ is an element of $\mathcal{L} (B_{m} )$ for some
$B_{m}=\mathrm{O} (m^{1/2} )$ (Lemma~\ref{LemmaEstimatorL1Bound}).
The algorithm is known not to possess as good approximation properties
as the other algorithms considered in this paper. However, this is
compensated by $B_{m}$ not growing too fast, hence, also the estimation
error does not grow too fast.

\subsubsection{Orthogonal Greedy Algorithm (a.k.a. Orthogonal Matching Pursuit)}
\label{SectionOGA}

The Orthogonal Greedy Algorithm (OGA) (e.g., Barron \textit{et~al.}
\cite{Baretal08}) is
also known as Orthogonal Matching Pursuit. Figure~\ref{fig2} recalls that the
OGA finds the next regressor to be included based on the same criterion
as for PGA, but at each $m$ iteration, it re-estimates the regression
coefficients by OLS using the selected regressors. For convenience,
the OLS projection operator is defined by $P_{X}^{m}$ where the $m$
stresses that one is only using the regressors included up to iteration
$m$, that is, $P_{X}^{m}Y=\sum_{k=1}^{m}b_{kn}X^{s (k )}$ for
OLS coefficients $b_{kn}$'s. Hence, in some circumstances, the OGA
is too time consuming, and may require the use of generalized inverses
when regressors are highly correlated. However, Pati, Rezaiifar and Krishnaprasad \cite{PatRezKri93}
give a faster algorithm for its estimation.

\subsubsection{Relaxed Greedy Algorithm}\label{SectionRGA}

The Relaxed Greedy Algorithm (RGA) is a less popular method, which
however has the same estimation complexity of the PGA. It is reviewed
in Figure~\ref{fig3}. The RGA updates taking a convex combination of the existing
regression function with the new predictor. The RGA does not shrink
the estimated coefficient at each step, but does shrink the regression
from the previous iteration $j-1$ by an amount $1-w_{j}$, where
$w_{j}=j^{-1}$. Other weighting schemes such that $w_{j}\in
(0,1 )$
and $w_{j}=\mathrm{O} (j^{-1} )$ can be used and the results hold
as they are (see Remark~2.5 in Barron \textit{et~al.}~\cite{Baretal08}). The weight sequence
$w_{j}=j^{-1}$ produces an estimator that has the simple average
structure $F_{m}=\sum_{j=1}^{m} (\frac{j}{m} )g_{j}
(X )$.\vspace*{1pt}

\begin{figure}
\fbox{
\begin{tabular}{l}
Set: \\
$m\in\mathbb{N}$\\
$F_{0}:=0$\\
For: $j=1,2,\ldots,m$\\
$w_{j}=1/j$\\
$s (j ):=\arg\max_{k}\llvert \langle Y-(1-w_{j} )F_{j-1},X^{ (k )} \rangle_{n}\rrvert $\\[2pt]
$g_{j} (X ):= \langle Y- (1-w_{j})F_{j-1},X^{s (j )} \rangle_{n}X^{s (j)}$\\[2pt]
$F_{j}:= (1-w_{j} )F_{j-1}+g_{j} (X )$
\end{tabular}}
\caption{RGA.}\vspace*{-8pt}\label{fig3}
\end{figure}

The RGA is advocated by Barron \textit{et~al.} \cite{Baretal08},
as it possesses better
theoretical properties than PGA ($L_{2}$-Boosting) and it is simpler
to implement\vspace*{1pt} than the OGA. At each stage $j$, the algorithm solves
$\min_{k,b}\llvert Y- (1-w_{j} )F_{j-1}+w_{j}X^{
(k )}b\rrvert _{n}^{2}$.
It is possible to also consider the case where $w_{j}$ is not fixed
in advance, but estimated at each iteration. In Figure~\ref{fig3}, one just
replaces the line defining $s (j )$ with
%
%
\begin{equation}
\bigl[s (j ),w_{j} \bigr]:=\arg\max_{k\leq K,w\in
[0,1 ]}\bigl
\llvert\bigl\langle Y- (1-w )F_{j-1},X^{
(k )} \bigr
\rangle_{n}\bigr\rrvert.\label{EQRGALineSearch}
\end{equation}
The asymptotic results hold as they are, as in this case, the extra
optimization can only reduce $\algo (m,B )$, the error in
the algorithm. The same remark holds for the next algorithms.

\subsubsection{Constrained greedy and Frank--Wolfe Algorithms}\label{SectionCGA&FWA}

The Constrained Greedy Algorithm (CGA) is a variation of the RGA.
It is used in Sancetta \cite{San14} in a slightly
different context. The
Frank--Wolfe Algorithm (FWA) (Frank and Wolfe \cite
{FraWol56}; see Clarkson \cite{Cla10}, Jaggi \cite{Jag13}, Freund, Grigas and Mazumder \cite{FreGriMaz13}, for recent results on its
convergence) is a well-known algorithm for the optimization of functions
under convex constraints. Figure~\ref{fig4} review the algorithms. The two
algorithms are similar, though some notable differences are present.
The FWA chooses at each iteration the regressor that best fits the
residuals from the previous iteration model. Moreover, the regression
coefficient is chosen as the value of the constraint times the sign
of the correlation of the residuals with the chosen regressor. On
the other hand, the difference of the CGA from the RGA is that at
each step the estimated regression coefficient is constrained to be
smaller in absolute value than a pre-specified value $\bar{B}$. When
the function one wants to estimate is known to lie in $\mathcal
{L} (1 )$,
the algorithm is just a simplified version of the Hilbert Space Projection
algorithm of Jones \cite{Jon92} and Barron \cite{Bar93} and have been studied
by several authors for estimation of mixture of densities (Li and
Barron \cite{LiBar00}, Rakhlin, Panchenko and Mukherjee \cite{RakPanMuk05}, Klemel{\"{a}}
\cite{Kle07}, Sancetta \cite{San13}).

\begin{figure}
\fbox{
\begin{tabular}{l@{\hspace*{29pt}}l}
CGA & FWA \\
Set: & Set:\\
$m\in\mathbb{N}$ & $m\in\mathbb{N}$\\
$F_{0}:=0$ & $F_{0}:=0$\\
$\bar{B}<\infty$ & $\bar{B}<\infty$\\
For: $j=1,2,\ldots,m$ & For: $j=1,2,\ldots,m$\\
$w_{j}=1/j$ & $w_{j}:=2/ (1+j )$\\
$s (j ):=\arg\max_{k}\llvert \langle Y-
(1-w_{j} )F_{j-1},X^{ (k )} \rangle_{n}\rrvert $ & $s (j ):=\arg\max _{k}\llvert \langle Y-F_{j-1},X^{ (k )} \rangle_{n}\rrvert $\\[2pt]
$b_{j}:=\frac{1}{w_{j}} \langle Y- (1-w_{j})F_{j-1},X^{s (j )} \rangle_{n}$ & $b_{j}:=\bar{B}\operatorname{sign} ( \langle Y-F_{j-1},X^{s (j )}\rangle_{n} )$\\[2pt]
$g_{j} (X ):=\operatorname{sign} (b_{j} ) (\llvert b_{j}\rrvert \wedge\bar{B} )X^{s
(j )}$ & $g_{j} (X ):=b_{j}X^{s (j )}$\\[2pt]
$F_{j}:= (1-w_{j} )F_{j-1}+w_{j}g_{j} (X )$ & $F_{j}:= (1-w_{j} )F_{j-1}+w_{j}g_{j} (X )$
\end{tabular}}
\caption{CGA and FWA.}\label{fig4}
\end{figure}

At each step $j$, the CGA solves $\min_{k,\llvert b\rrvert \leq\bar
{B}}\llvert Y- (1-w_{j} )F_{j-1}+w_{j}X^{ (k
)}b\rrvert _{n}^{2}$.
Under the square loss with regression coefficients satisfying $\sum
_{k=1}^{K}\llvert b_{k}\rrvert \leq\bar{B}$,
the FWA reduces to minimization of the linear approximation\vspace*{1pt} of the
objective function, minimized over the simplex, that is, $\min
_{k,\llvert b\rrvert \leq\bar{B}} \langle bX^{ (k
)},F_{j-1}-Y \rangle_{n}$
with\vspace*{1pt} update of $F_{j}$ as in Figure~\ref{fig4}. Despite the differences, both
the CGA and the FWA lead to the solution of the Lasso problem. In
particular, the regression coefficients are the solution to the following
problem:
\[
\min_{b_{1},b_{2},\ldots,b_{K}}\Biggl\llvert Y-\sum_{k=1}^{K}b_{k}X^{
(k )}
\Biggr\rrvert_{n},\qquad\mbox{such that }\sum
_{k=1}^{K}\llvert b_{k}\rrvert\leq\bar{B}.
\]
The above is the standard Lasso problem due to Tibshirani
\cite{Tib96}.
In particular, CGA and FWA solve the above problem as $m\rightarrow
\infty$,
\[
\llvert Y-F_{m}\rrvert_{n}^{2}\leq\inf
_{\mu\in\mathcal{L}
(\bar{B} )}\bigl\llvert Y-\mu(X )\bigr\rrvert_{n}^{2}+
\frac
{\bar{B}^{2}}{m}
\]
(Lemma~\ref{Lemmarelaxedgreedy} and~\ref{LemmafrankWolfeApproximation},
in Section~\ref{SectionProofs}, where for simplicity, only the weighting
schemes as in Figure~\ref{fig4} are considered). The complexity of the estimation
procedure is controlled by $\bar{B}$. This parameter can be either
chosen based on a specific application, or estimated via cross-validation,
or splitting the sample into estimation and validation sample.

The CGA and FGA also allows one to consider the forecast combination
problem with weights in the unit simplex, by minor modification. To
this end, for the CGA let
%
%
\begin{equation}
g_{j} (X ):= \bigl[ (b_{j}\wedge1 )\vee0
\bigr]X^{s (j )},\label{EQRGASimplex}
\end{equation}
so that $\bar{B}=1$ and the estimated $b_{j}$'s parameters are bounded
below by zero. For the FWA change,
\[
s (j ):=\arg\max_{k} \bigl\langle Y-F_{j-1},X^{
(k )}
\bigr\rangle_{n}; b_{j}:=\bar{B}\operatorname{sign} \bigl( \bigl\langle
Y-F_{j-1},X^{s (j )} \bigr\rangle_{n} \bigr)\vee0,
\]
where one does not use the absolute value in the definition of $s
(j )$.
(This follows from the general definition of the Frank--Wolfe Algorithm,
which simplifies to the algorithm in Figure~\ref{fig4} when $\sum
_{k=1}^{K}\llvert b_{k}\rrvert \leq\bar{B}$.)
Hence, the resulting regression coefficients are restricted to lie
on the unit simplex.

As for the RGA, for the CGA and FWA it is possible to estimate $w_{j}$
at each greedy step. For the CGA, this requires to change the line
defining $s (j )$ with (\ref{EQRGALineSearch}). Similarly,
for the FWA, one adds the following line just before the definition
of $F_{j}$:
\[
w_{j}=\arg\min_{w\in[0,1 ]}\bigl\llvert Y- (1-w
)F_{j-1}+wg_{j} (X )\bigr\rrvert_{n}^{2}.
\]
These steps can only reduce the approximation error of the algorithm,
hence, the rates of convergence derived for the fixed sequence $w_{j}$
are an upper bound for the case when $w_{j}$ is estimated at each
step.

\subsection{Discussion}\label{Sectiondiscussion}

\subsubsection{Objective function}

The objective function is the same one used in B{\"u}hlmann
\cite{Buh06}, which
is the integrated square error (ISE), where integration is w.r.t.
the true distribution of the regressors (note that the expectation
is w.r.t. $X'$ only). This objective function is zero if the (out
of sample) prediction error is minimized (recall that $\mu_{0}
(X )=\mathbb{E} [Y|X ]$),
and for this reason it is used in the present study. Under this objective,
some authors derive consistency, but not explicit rates of convergence
(e.g., B{\"u}hlmann \cite{Buh06}, Lutz and B{\"u}hlmann
\cite{LutBuh06}). An exception is Barron \textit{et~al.}
\cite{Baretal08}, who derive rates of convergence for the mean integrated
square error. Rates of convergence of greedy algorithms are usually
derived under a weaker norm, namely the empirical $L_{2}$ norm and
the results hold in probability (e.g., B{\"u}hlmann and van~de Geer \cite{Buhvan11},
and references therein). This is essentially equivalent to assuming
a fixed design for the regressors. The empirical $L_{2}$ norm has
been used to show consistency of Lasso, hence deriving results under
this norm allows one to compare to Lasso in a more explicit way. Convergence
of the empirical $L_{2}$ norm does not necessarily guarantee that
the prediction error is minimized, asymptotically.

%
%
\begin{example}\label{ExamplelossFunctionControl}Let $F_{m}
(X )=\sum_{k=1}^{K}X^{ (k )}b_{kn}$
be the output of one of the algorithms, where the subscript $n$ is
used to stress that $b_{kn}$ depends on the sample. Also, let $Z:=Y-\mu
_{0} (X )$,
and $\mu_{0} (X )=\sum_{k=1}^{K}X^{ (k )}b_{k0}$,
where the $b_{k0}$'s are the true coefficients. Control of the empirical
$L_{2}$ norm only requires control of Control of $ \langle Z,\sum
_{k=1}^{K}X^{ (k )} (b_{kn}-b_{k0} )
\rangle_{n}$
(e.g., Lemma 6.1 in B{\"u}hlmann and van~de Geer \cite{Buhvan11}) and
this quantity
tends to be $\mathrm{O}_{p} (m\ln K/n )$ under regularity conditions.
On the other hand, control of the $L_{2}$ norm (i.e., ISE) also requires
control of $ (1-\mathbb{E} )\llvert \sum_{k=1}^{K}X^{
(k )} (b_{kn}-b_{k0} )\rrvert _{n}^{2}$. Sufficient conditions for this term to be $\mathrm{O}_{p} (m\ln K/n )$
are often used, but in important cases such as dependent non-mixing
random data, this does not seem to be the case anymore. Hence, this
term is more challenging to bound and requires extra care (see
van~de Geer \cite{van14}, for results on how to bound such
a term in an i.i.d.
case).
\end{example}

\subsubsection{Dependence conditions}

Absolute regularity is convenient, as it allows to use decoupling
inequalities. In consequence, the same rate of convergence under i.i.d.
observations holds under beta mixing when the mixing coefficients
decay fast enough. Many time series models are beta mixing. For example,
any finite order ARMA model with i.i.d. innovations and law absolutely
continuous w.r.t. the Lebesgue measure satisfies geometric mixing
rates (Mokkadem \cite{Mok88}). Similarly, GARCH models
and more generally
models that can be embedded in some stochastic recursive equations
are also beta mixing with geometric mixing rate for innovations possessing
a density w.r.t. the Lebesgue measure (e.g., Basrak, Davis and Mikosch \cite{BasDavMik02}, for
details: they derive the results for strong mixing, but the result
actually implies beta mixing). Many positive recurrent Markov chains
also satisfy geometric absolute regularity (e.g., Mokkadem \cite{Mok90}).
Hence, while restrictive, the geometric mixing rate of Conditions
\ref{ConditionabsoluteRegularityBounded} and~\ref{ConditionabsoluteRegularity}
is a convenient condition satisfied by common time series models.

In Condition~\ref{ConditionabsoluteRegularity}, (\ref{EQboundVariance})
is used to control the moments of the random variables. The geometric
mixing decay could be replaced with polynomial mixing at the cost
of complications linking the moments of the random variables (i.e.,
(\ref{EQboundVariance})) and their mixing coefficients (e.g.,
Rio \cite{Rio00}, for details).

Condition~\ref{Conditiondependence} only controls dependence in
terms of some conditional moments of the centered random variables.
Hence, if the dependence on the past decreases as we move towards
the future, the centered variables will have conditional moment closer
and closer to zero. On the other hand, Conditions~\ref{ConditionabsoluteRegularityBounded}
and~\ref{ConditionabsoluteRegularity} control dependence in terms
of the sigma algebra generated by the future and the past of the data.
This is much stronger than controlling conditional expectations, and
computation of the resulting mixing coefficients can be very complicated
unless some Markov assumptions are made as in Mokkadem
\cite{Mok88,Mok90}
or Basrak, Davis and Mikosch \cite{BasDavMik02} (see Doukhan and Louhichi \cite{DouLou99}, for
further discussion
and motivation).

\subsubsection{Examples for Conditions \texorpdfstring{\protect\ref{ConditionabsoluteRegularity}}{3} and 
\texorpdfstring{\protect\ref{Conditiondependence}}{4}}\label{Sectionexamples}

To highlight the scope of the conditions and how to establish them
in practice, consider a simple non-trivial example.

%
%
\begin{example}\label{ExampleshortMemory}Let $\mu_{0} (X
)=g (X^{ (k )};k\leq K )$,
where $g$ satisfies
\[
\bigl\llvert g \bigl(x^{ (k )};k\leq K \bigr)-g \bigl(z^{
(k )};k
\leq K \bigr)\bigr\rrvert\lesssim\sum_{k=1}^{K}
\lambda_{k}\bigl\llvert x^{ (k )}-z^{ (k )}\bigr\rrvert
\]
for $\sum_{k=1}^{K}\lambda_{k}\leq1$, $\lambda_{k}\geq0$ and
$g (x^{ (k )};k\leq K )=0$
when $x^{ (k )}=0$ for all $k\leq K$. Since $K\rightarrow
\infty$
with $n$, it is natural to impose this condition which is of the
same flavor as $\sum_{k=1}^{K}\llvert b_{k}\rrvert \leq B$ in the linear
model. Suppose that $ (Z_{i} )_{i\in\mathbb{Z}}$ is a sequence
of independent random variables ($Z=Y-\mathbb{E} [Y|X ]$)
with finite $p$ moments, and independent of the regressors $
(X_{i} )_{i\in\mathbb{Z}}$.
The regressors admit the following vector autoregressive representation,
$X_{i}=HW_{i}$, where $H$ is a $K\times L$ matrix with positive
entries and rows summing to one; $W_{i}=AW_{i-l}+\varepsilon_{i}$,
$A$ is a diagonal $L\times L$ matrix with entries less than one
in absolute values, and $ (\varepsilon_{i} )_{i\in\mathbb{Z}}$
is a\vspace*{1pt} sequence of i.i.d. $L$ dimensional random variables with finite
$2p$ moments, that is, $\mathbb{E}\llvert \varepsilon_{i,k}\rrvert
^{2p}<\infty$,
where $\varepsilon_{i,k}$ is the $k$th entry in~$\varepsilon_{i}$.
Throughout, the $K$ dimensional vectors are column vectors.
\end{example}

If one takes $L=K$ and $H$ to be diagonal, $X_{i}=W_{i}$. As
$K\rightarrow\infty$,
the process is not necessarily mixing. Hence, one is essentially required
to either keep $L$ fixed or impose very restrictive structure on
the innovations in order to derive mixing coefficients. Luz and
B{\"u}hlmann \cite{Buh06} consider vector autoregressive
models (VAR) with
the dimension
of the variables increasing to infinity. They then assume that the
model is strongly mixing. However, it is unclear that a VAR of increasing
dimensionality can be strongly mixing. The mixing coefficients of
functions of independent random variables are bounded above by the
sum of the mixing coefficients of the individual variables (e.g.,
Theorem 5.1 in Bradley \cite{Bra05}). If the number of
terms in the sum goes
to infinity (i.e., $K$ in the present context, $q$ in Luz and
B{\"u}hlmann \cite{Buh06}), such VAR may not be strongly
mixing. Even using a
known results
on Markov chain, it is not possible to show that VAR models with increasing
dimension are mixing without very restrictive conditions on the innovations
(e.g., condition iii in Theorem~1$'$ in Mokkadem \cite{Mok88}).

Restrictions such as $A$ being diagonal or $ (X_{i} )_{i\in
\mathbb{Z}}$
and $ (Z_{i} )_{i\in\mathbb{Z}}$ being independent are only
used to simplify the discussion, so that one can focus on the standard
steps required to establish the validity of the conditions in Example
\ref{ExampleshortMemory}. The above model can be used to show how
to check Conditions~\ref{ConditionabsoluteRegularity} and~\ref{Conditiondependence}
and how Condition~\ref{ConditionabsoluteRegularity} can fail.

%
%
\begin{lemma}\label{LemmaconditionAbsoluteRegularity}Consider the
model in Example~\ref{ExampleshortMemory}. Suppose that $\varepsilon_{i}$
has a density w.r.t. the Lebesgue measure and $L$ is bounded. Then
Condition~\ref{ConditionabsoluteRegularity} is satisfied.
\end{lemma}

%
%
\begin{lemma}\label{LemmaconditionDependence}Consider the model
in Example~\ref{ExampleshortMemory}. Suppose that $\varepsilon_{i,k}$
only takes values in $ \{ -1,1 \} $ with equal probability
for each $k$, $L=K$ and $H$ is the identity matrix (i.e., $X_{i}=W_{i}$),
while all the rest is as in Example~\ref{ExampleshortMemory}. Then
Condition~\ref{ConditionabsoluteRegularity} is not satisfied, but
Condition~\ref{Conditiondependence} is satisfied.
\end{lemma}

The proof of these two lemmas~-- postponed to Section~\ref
{SectionproofLemmataConditions}~-- shows how the conditions can be verified.

The next examples provides details on the applicability of Condition
\ref{Conditiondependence} to possibly long memory processes. In
particular, the goal is to show that Corollary~\ref{CorollarydependenceFiniteMoments}
can be applied. In consequence, new non-trivial models and conditions
are allowed. In these examples, the rates of convergence implied by
Corollary~\ref{CorollarydependenceFiniteMoments} are comparable
to, or better than the ones in B{\"u}hlmann and van~de Geer \cite{Buhvan11} which
require i.i.d. observations. However, one needs to restrict attention
to regressors whose population Gram matrix has full rank ($\rho_{K}>0$).
The following only requires stationarity and ergodicity of the error
terms.

%
%
\begin{example}\label{ExamplestationaryErrorsIIDXs}Let $
(Z_{i} )_{i\in\mathbb{Z}}$
be a stationary ergodic sequence with moments of all orders, and suppose
that $ (X_{i} )_{i\in\mathbb{Z}}$ is i.i.d., independent
of the $Z_{i}$'s, and with zero mean and moments of all orders and
such that $\rho_{K}>0$. Moreover, suppose that $\mu_{0}\in\mathcal{L}$.
By independence of $X_{i}$ and the $Z_{i}$'s, and the fact the that
$X_{i}$'s are i.i.d. mean zero, it follows that $\mathbb
{E}_{0}Z_{i}X_{i}^{ (k )}=0$.
Similarly, $\mathbb{E}_{0} (1-\mathbb{E} )\llvert X_{i}^{ (k
)}\rrvert ^{2}=0$
for $i>0$. Finally, given that $\mu_{0}\in\mathcal{L}$, $\Delta
(X )=\mu_{B} (X )-\mu_{0} (X )=0$
by choosing $B$ large enough so that $\mu_{B}=\mu_{0}$. Hence,\vspace*{1pt} this
implies that $\sup_{n}d_{n,p}<\infty$ in Corollary~\ref{CorollarydependenceFiniteMoments},
though for the CGA and FWA it is necessary to assume $\mu_{0}\in
\mathcal{L} (\bar{B} )$
and not just $\mu_{0}\in\mathcal{L}$, or just $\mu_{0}\in\mathcal{L}$
but $\bar{B}\rightarrow\infty$.
\end{example}

Remarkably, Example~\ref{ExamplestationaryErrorsIIDXs} shows that
if the regressors are i.i.d., it is possible to achieve results as
good as the ones derived in the literature only assuming ergodic stationary
noise. The next example restricts the noise to be i.i.d., but allows
for long memory Gaussian regressors and still derives convergence
rates as fast as the ones of Example~\ref{ExamplestationaryErrorsIIDXs}.

%
%
\begin{example}\label{ExamplelongMemory}
\textit{Let} $X_{i}^{ (k
)}=\sum_{l=0}^{\infty}a_{lk}\varepsilon_{i-l,k}$,
\textit{where} $ (\varepsilon_{i,k} )_{i\in\mathbb{Z}}$ \textit{is a sequence
of i.i.d. standard Gaussian random variables}, \textit{and} $a_{0k}=1$,
$a_{lk}=l^{- (1+\epsilon)/2}$
$\epsilon\in(0,1]$ \textit{for} $l>0$. \textit{Also}, \textit{suppose that} $ (X_{i}
)_{i\in\mathbb{Z}}$
\textit{is independent of} $ (Z_{i} )_{i\in\mathbb{Z}}$, \textit{which is
i.i.d. with moments of all orders. It is shown in Section}~\ref{SectionProofExample}
\textit{that for this} MA($\infty$) \textit{model with Gaussian errors},
\[
\bigl\llvert\mathbb{E}_{0} (1-\mathbb{E} )\bigl\llvert
X_{i}^{
(k )}\bigr\rrvert^{2}\bigr\rrvert
_{p}\lesssim i^{- (1+\epsilon
)}
\]
\textit{when} $i>0$. \textit{Hence}, \textit{in Condition}~\ref{Conditiondependence} $\sup
_{n}d_{n,p}<\infty$
\textit{for any} $p<\infty$, \textit{and in consequence}, \textit{one can apply Corollary}~\ref{CorollarydependenceFiniteMoments}
\textit{if} $\rho_{K}>0$ \textit{and the true function is in} $\mathcal{L}$ \textit{or}
$\mathcal{L} (\bar{B} )$
\textit{for the CGA and FWA. For the CGA and FWA}, $\rho_{K}=0$ \textit{is allowed}.
\end{example}

\section{Implementation and numerical comparison}\label{Sectionimplementation}

\subsection{Vectorized version}

Vectorized versions of the algorithms can be constructed. These versions
make quite clear the mechanics behind the algorithms. The vectorized
versions are useful when the algorithms are coded using scripting
languages or when $n$ and $K$ are very large, but $K=\mathrm{o} (n )$.
In this case, the time dimension $n$ could be about $\mathrm{O}
(10^{7} )$
or even $\mathrm{O} (10^{8} )$ and the cross-sectional dimension
$K=\mathrm{O} (10^{3} )$. The memory requirement to store a matrix
of doubles of size $10^{7}\times10^{3}$ is in excess of 70 gigabytes,
often too much to be stored in RAM on most desktops. On the other
hand, sufficient statistics such as $\mathbf{X}^{T}\mathbf{X}$ and
$\mathbf{X}^{T}\mathbf{Y}$ ($\mathbf{X}$ being the $n\times K$
matrix of regressors and $\mathbf{Y}$ the $n\times1$ vector of dependent
variables and the subscript $T$ stands for transpose) are manageable
and can be updated through summation.

Figure~\ref{fig5} shows vectorized versions of the algorithms. Of course, it
is always assumed that the regressors have been standardized, that is,
$\diag (\mathbf{X}^{T}\mathbf{X}/n )=I_{K}$, the identity
matrix, where $\diag (\cdot)$ stands for the diagonal
matrix constructed from the diagonal of its matrix argument. The symbol
$0_{K}$ is the $K$ dimensional vector of zeros, while for other
vector quantities, the subscript denotes the entry in the vector,
which are assumed to be column vectors.

\begin{figure}
\tabcolsep=9.6pt
\begin{tabular}{|lll|}
\multicolumn{1}{l}{PGA} & OGA & \multicolumn{1}{l}{RGA}\\[-6pt]
\multicolumn{3}{@{}c@{}}{\hrulefill}\\[-4pt]
\rule{0pt}{12pt}$\mbox{Set:}$ & & \\
$C=\mathbf{X}^{T}\mathbf{Y}/n$ & $C=\mathbf{X}^{T}\mathbf{Y}/n$ &
$C=\mathbf{X}^{T}\mathbf{Y}/n$\\
$D=\mathbf{X}^{T}\mathbf{X}/n$ & $D=\mathbf{X}^{T}\mathbf{X}/n$ &
$D=\mathbf{X}^{T}\mathbf{X}/n$\\
$b=0_{K}$ & $b=0_{K}$ & $b=0_{K}$\\
$\nu\in(0,1 )$ & & \\
$\mbox{For:}j=1,2,\ldots,m$ & & \\
$A=C-Db$ & $A=C-Db$ & $A=C- (1-\frac{1}{j} )Db$\\
$s (j )=\arg\max_{k\leq K}\llvert A_{k}\rrvert $ & $s(j )=\arg\max_{k\leq K}\llvert A_{k}\rrvert $ & $s (j)=\arg\max_{k\leq K}\llvert A_{k}\rrvert $\\[2pt]
$a=0_{K}$ & $P_{X}^{j}$ as in Figure~\ref{fig2} & $a=0_{K}$\\
$a_{s (j )}=A_{s (j )}$ & $b=P_{X}^{j}Y$ &
$a_{s (j )}=A_{s (j )}$\\
$b=b+\nu a$ & & $b= (1-\frac{1}{j} )b+\frac{1}{j}a$
\\[6pt]
CGA & FWA\phantom{\tsub{\rule{0pt}{10pt}}} & \\[-8pt]
\multicolumn{3}{@{}c@{}}{\hrulefill}\\[-4pt]
\rule{0pt}{12pt}$\mbox{Set:}$ & & \\
$C=\mathbf{X}^{T}\mathbf{Y}/n$ & $C=\mathbf{X}^{T}\mathbf{Y}/n$ & \\
$D=\mathbf{X}^{T}\mathbf{X}/n$ & $D=\mathbf{X}^{T}\mathbf{X}/n$ & \\
$b=0_{K}$ & $b=0_{K}$ & \\
$\bar{B}<\infty$ & $\bar{B}<\infty$ & \\
$\mbox{For:}j=1,2,\ldots,m$ & & \\
$A=C- (1-\frac{1}{j} )Db$ & $A=C-Db$ & \\
$s (j )=\arg\max_{k\leq K}\llvert A_{k}\rrvert $ & $s
(j )=\arg\max_{k\leq K}\llvert A_{k}\rrvert $ & \\
$a=0_{K}$ & $a=0_{K}$ & \\
$a_{s (j )}=\operatorname{sign} (A_{s (j )} ) (j\llvert A_{s (j )}\rrvert \wedge\bar{B} )$ &$a_{s (j )}=\operatorname{sign} (A_{s (j )} )\bar
{B}$ & \\[3pt]
$b= (1-\frac{1}{j} )b+\frac{1}{j}a$ & $b= (1-\frac
{2}{1+j} )b+\frac{2}{1+j}a$ & \\
\hline
\end{tabular}
\caption{Vectorized versions of the algorithms.}\label{fig5}
\end{figure}
%

\subsection{Choosing the number of iterations}

In order to achieve the bounds in the theorem, $m$ needs to be chosen
large enough for the algorithm to perform well in terms of approximation
error (see Lemmas~\ref{LemmaL2Boosting},~\ref{LemmaOGA} and~\ref{Lemmarelaxedgreedy}).
Nevertheless, an excessively large $m$ can produce poor results as
shown in the theorems with the exception of CGA and FWA. In consequence,
guidance on the number of iterations $m$ is needed. The number of
regressors can be left unconstrained in many situations, as long as
the dependence is not too strong. The number of iterations can be
chosen following results in the literature. Suppose the $F_{m}$ estimator
in the algorithm can be represented as $F_{m} (X )=P^{m}Y$
for some suitable projection operator $P^{m}$. Then one may choose
the number of iterations according the following AIC criterion:
\[
\ln\bigl(\bigl\llvert Y-F_{m} (X )\bigr\rrvert_{n}^{2}
\bigr)+2\mathrm{df} \bigl(P^{m} \bigr)/n,
\]
where $\mathrm{df} (P^{m} )$ are the degrees of freedom of the prediction
rule $P^{m}$, which are equal to the sum of the eigenvalues of $P^{m}$,
or equivalently they are equal to the trace of the operator.
B{\"u}hlmann \cite{Buh06} actually suggests using the
modified AIC based on
Hurvich, Simonoff and Tsai \cite{HurSimTsa98}:
\[
\ln\bigl(\bigl\llvert Y-F_{m} (X )\bigr\rrvert_{n}^{2}
\bigr)+\frac
{1+\mathrm{df} (P^{m} )/n}{1- (\mathrm{df} (P^{m} )+2 )/n}.
\]

For ease of exposition, let $\mathbf{X}_{m}$ be the $n\times m$
matrix of selected regressors and denote by $\mathbf{X}_{m}^{s
(j )}$
the $j$th column of $\mathbf{X}_{m}$. For the PGA, B{\"u}hlmann and Yu
\cite{BuhYu03} show that the degrees of freedom are
given by the trace
of
\[
\mathcal{B}_{m}:=I_{n}-\prod_{j=1}^{m}
\biggl(I_{n}-\nu\frac{\mathbf
{X}_{m}^{s (j )} (\mathbf{X}_{m}^{s (j
)} )'}{ (\mathbf{X}_{m}^{s (j )} )'\mathbf
{X}_{m}^{s (j )}} \biggr),
\]
where $I_{n}$ is the $n$ dimensional identity matrix.

The trace of the hat matrix $\mathcal{B}_{m}:=\mathbf{X}_{m}
(\mathbf{X}_{m}^{T}\mathbf{X}_{m} )^{-1}\mathbf{X}_{m}^{T}$
gives the degrees of freedom for the OGA, that is, $\operatorname{Trace}
(\mathcal{B}_{m} )=m$.

Unfortunately, the projection matrix of the RGA is complicated and
the author could not find a simple expression. Nevertheless, the degrees
of freedom could be estimated (e.g., Algorithm 1 in Jianming \cite{Ye98}).

Choice of $\bar{B}$ is equivalent to the choice of the penalty constant
in Lasso. Hence, under regularity conditions (Zou \textit{et~al.} \cite{ZouHasTib07},
Tibshirani and Taylor \cite{TibTay12}) the degrees of freedom of
the CGA and FWA are
approximated
by the number of non-zero coefficients or the rank of the population
Gram matrix of the selected variables. Alternatively, one has to rely
on cross-validation to choose $m$ for the PGA, OGA, RGA and $\bar{B}$
for the CGA and FWA.

\subsection{Numerical results}\label{SectionnumericalResults}

To assess the finite performance of the algorithms a comprehensive
set of simulations is carried out for all the algorithms. It is worth
mentioning that the CGA and FWA are equivalent to Lasso, hence, conclusions
also apply to the Lasso, even though the conditions used for consistency
are very different.

For each Monte Carlo set up, 100 simulations are run, where the sample
size is $n={}$20,100. Consider the model
\begin{eqnarray*}
Y_{i}&=&\sum_{k=1}^{K}X_{i}^{ (k )}b_{k}+Z_{i},\qquad
 X_{i}^{ (k )}=\sum_{s=0}^{S}
\theta_{s}\varepsilon_{i-s,k},\qquad  Z_{i}=
\frac{\kappa}{\sigma}\sum_{s=0}^{S}\theta
_{s}\varepsilon_{i-s,0},
\end{eqnarray*}
where\vspace*{1pt} $K=100$, $\kappa^{2}=\frac{\operatorname{Var} (\sum_{k=1}^{K}X_{i}^{ (k )}b_{k}
)}{\operatorname{Var} (\sum_{s=0}^{S}\theta_{s}\varepsilon_{i-s,0} )}$,
so that $\sigma^{2}\in\{ 8,0.25 \} $ is the signal to noise
ratio, corresponding roughly to an $R^{2}$ of $0.89,0.2$. The innovations
$ \{ (\varepsilon_{i,k} )_{i\in\mathbb
{Z}}\dvt k=0,1,\ldots,K \} $
are collections of i.i.d. standard normal random variables. For $k,l>0$,
$\mathbb{E}\varepsilon_{i,k}\varepsilon_{i,l}=\omega^{\llvert
k-l\rrvert }$
with $\omega= \{ 0,0.75 \} $ with convention $0^{0}=1$,
that is, a Toeplitz covariance matrix. Moreover, $\mathbb
{E}\varepsilon_{i,0}\varepsilon_{i,k}=0$
for any $k>0$. Finally, $ \{ \theta_{s}\dvt s=0,1,\ldots,S \} $
is as follows:

\begin{longlist}[Case WD:]
\item[Case ID:] $\theta_{0}=1$ and $\theta_{s}=0$ if $s>0$;

\item[Case WD:] $\theta_{s}= (0.95 )^{s}$ with $S=100+n$;

\item[Case SD:] $\theta_{s}= (s+1 )^{-1/2}$ with $S=1000+n$.
\end{longlist}

In other words, the above model allows for time dependent $Z_{i}$'s and
$X_{i}$'s as well for correlated regressors (when $\omega>0$). However,
the $X$ and the $Z$ are independent by construction. By different
choice of regression coefficients $b_{k}$'s, it is possible to define
different scenarios for the evaluation of the algorithms. These are
listed in the relevant subsections below. For each different scenario,
the mean integrated square error (MISE) from the simulations is computed:
that is, the Monte Carlo approximation of $\mathbb{E} [\mathbb
{E}'\llvert \mu_{0} (X' )-F_{m} (X' )\rrvert ^{2} ]$.
Standard errors were all relatively small, so they are not reported,
but available upon requests together with more detailed results.

The number of greedy steps $m$ or the bound $\bar{B}$ were chosen
by a cross-validation method for each of the algorithms (details are
available upon request). Hence, results also need to be interpreted
bearing this in mind, as cross-validation can be unstable at small
sample sizes (e.g., Efron \cite{Efr83}, see also
Sancetta \cite{San10},
for some
simulation evidence and alternatives, amongst many others). Moreover,
cross-validation is usually inappropriate for dependent data, often
leading to larger than optimal models (e.g., Burman and Nolan
\cite{BurNol92},
Burman, Chow and Nolan \cite{BurChoNol94}, for discussions and alternatives). Nevertheless,
this also allows one to assess how robust is the practical implementation
of the algorithms. Given the large amount of results, Section~\ref
{SectionremarksNumerical}
summarizes the main conclusions.

\subsection{Low-dimensional model}\vspace*{-1pt}

The true regression function has coefficients $b_{k}=1/3$ for $k=1,2,3$,
and $b_{k}=0$ for $k>3$.

\subsection{High-dimensional small equal coefficients}\vspace*{-1pt}

The true regression function has coefficients $b_{k}=1/K$, $k\leq K$.

\subsection{High-dimensional decaying coefficients}\vspace*{-1pt}

The true regression function has coefficients $b_{k}=k^{-1}$, $k\leq
K$.

\subsection{High-dimensional slowly decaying coefficients}\vspace*{-1pt}

The true regression function has coefficients $b_{k}=k^{-1/2}$, $k\leq
K$.

\subsection{Remarks on numerical results}\vspace*{-1pt}\label{SectionremarksNumerical}

Results from the simulations are reported in Tables \ref{tab2}--\ref{tab5}. These results show that the algorithms
are somehow comparable,
within a $\pm10\%$ relative performance. Overall, the PGA ($L_{2}$-Boosting)
is robust and often delivers the best results despite the theoretically
slower convergence rates.

On the other hand, the performance of the OGA is somehow disappointing
given the good theoretical performance. Table~\ref{tab2} shows that the OGA
can perform remarkably well under very\vadjust{\goodbreak} special circumstance, that is,
relatively large sample size $ (n=100 )$, time independent
and uncorrelated regressors and high signal to noise ratio. To some
extent, these are the conditions used by Zhang \cite
{Zha09} to show optimality
of the OGA.

The RGA, CGA and FWA provide good performance comparable to the PGA
and in some cases better, especially when the signal to noise ration
is higher. For example, Table~\ref{tab2} shows that these algorithms perform
well as long as the regressors are either uncorrelated or the time
dependence is low. Intuitively, time dependence leads to an implicit
reduction of information, hence it is somehow equivalent to estimation
with a smaller sample. This confirms the view that the PGA is usually
the most robust of the methods.

\begin{table}
\tabcolsep=0pt
\caption{MISE: low-dimensional, $K=100$}\label{tab2}
\begin{tabular*}{\tablewidth}{@{\extracolsep{\fill}}@{}lllllllllll@{}}
\hline
& \multicolumn{5}{l}{$n=20$} & \multicolumn{5}{l@{}}{$n=100$}\\[-6pt]
& \multicolumn{5}{l}{\hrulefill} & \multicolumn{5}{l@{}}{\hrulefill}\\
$(\omega,\sigma^{2} )$ & PGA & OGA & RGA & CGA & FWA & PGA & OGA & RGA & CGA & FWA\\
\hline
& \multicolumn{10}{c@{}}{Case ID}\\
$ (0,8 )$ & 0.40 & 0.51 & 0.36 & 0.36 & 0.40 & 0.08 & 0.03 & 0.09 & 0.09 & 0.09\\
$ (0,0.20 )$ & 0.59 & 0.87 & 0.93 & 0.75 & 0.77 & 0.47 & 0.52 & 0.49 & 0.44 & 0.44\\
$ (0.75,8 )$ & 0.25 & 0.39 & 0.26 & 0.36 & 0.35 & 0.09 & 0.15 & 0.07 & 0.13 & 0.13\\
$ (0.75,0.25 )$ & 0.86 & 1.20 & 1.29 & 1.00 & 1.14 & 0.50 & 0.45 & 0.49 & 0.48 & 0.47
\\[3pt]
& \multicolumn{10}{c@{}}{Case WD}\\
$ (0,8 )$ & 1.65 & 2.06 & 1.56 & 1.51 & 1.52 & 0.67 & 0.68 & 0.54 & 0.56 & 0.54\\
$ (0,0.20 )$ & 2.81 & 2.95 & 2.97 & 3.49 & 3.01 & 3.07 & 4.01 & 2.82 & 2.93 & 2.95\\
$ (0.75,8 )$ & 1.25 & 2.21 & 1.24 & 1.35 & 1.32 & 0.87 & 1.18 & 0.79 & 0.85 & 0.89\\
$ (0.75,0.25 )$ & 4.36 & 4.29 & 4.56 & 5.34 & 5.28 & 4.43 & 5.55 & 4.18 & 4.45 & 4.56
\\[3pt]
& \multicolumn{10}{c@{}}{Case SD}\\
$ (0,8 )$ & 1.26 & 1.63 & 1.26 & 1.24 & 1.25 & 0.50 & 0.50 & 0.43 & 0.42 & 0.41\\
$ (0,0.20 )$ & 2.31 & 2.36 & 2.55 & 2.61 & 2.53 & 2.20 & 2.72 & 2.16 & 2.15 & 2.14\\
$ (0.75,8 )$ & 0.88 & 1.82 & 0.91 & 0.98 & 1.00 & 0.63 & 0.86 & 0.58 & 0.58 & 0.58\\
$ (0.75,0.25 )$ & 3.28 & 3.37 & 3.58 & 3.88 & 4.14 & 3.13 & 3.74 & 3.05 & 3.11 & 3.12\\
\hline
\end{tabular*}
\end{table}

\begin{table}
\tabcolsep=0pt
\caption{MISE: high-dimensional small coefficients, $K=100$}\label{tab3}
\begin{tabular*}{\tablewidth}{@{\extracolsep{\fill}}@{}lllllllllll@{}}
\hline
& \multicolumn{5}{l}{$n=20$} & \multicolumn{5}{l@{}}{$n=100$}\\[-6pt]
& \multicolumn{5}{l}{\hrulefill} & \multicolumn{5}{l@{}}{\hrulefill}\\
$(\omega,\sigma^{2} )$ & PGA & OGA & RGA & CGA & FWA & PGA & OGA & RGA & CGA & FWA\\
\hline
& \multicolumn{10}{c@{}}{Case ID}\\
$ (0,8 )$ & 0.10 & 0.12 & 0.11 & 0.10 & 0.10 & 0.08 & 0.10 & 0.08 & 0.08 & 0.09\\
$ (0,0.20 )$ & 0.11 & 0.16 & 0.16 & 0.13 & 0.14 & 0.10 & 0.11 & 0.12 & 0.10 & 0.10\\
$ (0.75,8 )$ & 0.20 & 0.27 & 0.17 & 0.17 & 0.14 & 0.09 & 0.12 & 0.08 & 0.09 & 0.09\\
$ (0.75,0.25 )$ & 0.26 & 0.38 & 0.38 & 0.29 & 0.33 & 0.23 & 0.28 & 0.25 & 0.22 & 0.22
\\[3pt]
& \multicolumn{10}{c@{}}{Case WD}\\
$ (0,8 )$ & 0.35 & 0.40 & 0.35 & 0.37 & 0.33 & 0.27 & 0.36 &0.25 & 0.25 & 0.22\\
$ (0,0.20 )$ & 0.50 & 0.56 & 0.53 & 0.59 & 0.52 & 0.53 &0.68 & 0.51 & 0.54 & 0.56\\
$ (0.75,8 )$ & 0.65 & 0.88 & 0.65 & 0.63 & 0.50 & 0.34 &0.44 & 0.29 & 0.31 & 0.33\\
$ (0.75,0.25 )$ & 1.28 & 1.28 & 1.34 & 1.58 & 1.50 & 1.27 &1.62 & 1.22 & 1.32 & 1.37
\\[3pt]
& \multicolumn{10}{c@{}}{Case SD}\\
$ (0,8 )$ & 0.28 & 0.30 & 0.28 & 0.28 & 0.26 & 0.22 & 0.29 & 0.21 & 0.21 & 0.19\\
$ (0,0.20 )$ & 0.38 & 0.39 & 0.45 & 0.45 & 0.43 & 0.38 & 0.49 & 0.40 & 0.40 & 0.39\\
$ (0.75,8 )$ & 0.51 & 0.70 & 0.51 & 0.50 & 0.43 & 0.25 & 0.37 & 0.24 & 0.26 & 0.27\\
$ (0.75,0.25 )$ & 0.95 & 1.00 & 1.05 & 1.07 & 1.12 & 0.90 & 1.15 & 0.88 & 0.93 & 0.91\\
\hline
\end{tabular*}
\end{table}
%

\begin{table}
\tabcolsep=0pt
\caption{MISE: high-dimensional decaying coefficients, $K=100$}\label{tab4}
\begin{tabular*}{\tablewidth}{@{\extracolsep{\fill}}@{}lllllllllll@{}}
\hline
& \multicolumn{5}{l}{$n=20$} & \multicolumn{5}{l@{}}{$n=100$}\\[-6pt]
& \multicolumn{5}{l}{\hrulefill} & \multicolumn{5}{l@{}}{\hrulefill}\\
$(\omega,\sigma^{2} )$ & PGA & OGA & RGA & CGA & FWA & PGA & OGA & RGA & CGA & FWA\\
\hline
& \multicolumn{10}{c@{}}{Case ID}\\
$ (0,8 )$ & \phantom{0}2.28 & \phantom{0}2.60 & \phantom{0}2.33 & \phantom{0}2.26 & \phantom{0}2.13 & \phantom{0}1.44 & \phantom{0}2.03 & \phantom{0}1.39 & \phantom{0}1.59 & \phantom{0}1.78\\
$ (0,0.20 )$ & \phantom{0}2.42 & \phantom{0}3.61 & \phantom{0}3.72 & \phantom{0}3.02 & \phantom{0}3.12 & \phantom{0}2.23 & \phantom{0}2.48 & \phantom{0}2.56 & \phantom{0}2.22 & \phantom{0}2.25\\
$ (0.75,8 )$ & \phantom{0}3.98 & \phantom{0}5.32 & \phantom{0}3.25 & \phantom{0}3.19 & \phantom{0}2.80 & \phantom{0}1.70 & \phantom{0}2.38 & \phantom{0}1.56 & \phantom{0}1.75 & \phantom{0}1.79\\
$ (0.75,0.25 )$ & \phantom{0}5.51 & \phantom{0}7.89 & \phantom{0}8.18 & \phantom{0}6.27 & \phantom{0}7.02 & \phantom{0}4.46 & \phantom{0}5.49 & \phantom{0}5.07 & \phantom{0}4.33 & \phantom{0}4.37
\\[3pt]
& \multicolumn{10}{c@{}}{Case WD}\\
$ (0,8 )$ & \phantom{0}7.80 & \phantom{0}8.72 & \phantom{0}7.91 & \phantom{0}8.01 & \phantom{0}7.28 & \phantom{0}5.34 & \phantom{0}7.23 & \phantom{0}4.60 & \phantom{0}4.53 & \phantom{0}4.43\\
$ (0,0.20 )$ & 11.27 & 12.83 & 11.96 & 13.43 & 11.79 & 12.31 & 15.83 & 11.55 & 12.15 & 12.54\\
$ (0.75,8 )$ & 13.01 & 17.84 & 12.66 & 12.10 & 10.22 & \phantom{0}6.65 & \phantom{0}8.78 & \phantom{0}5.86 & \phantom{0}6.30 & \phantom{0}6.62\\
$ (0.75,0.25 )$ & 26.36 & 28.49 & 27.94 & 32.07 & 31.17 & 26.81 & 33.34 & 25.27 & 27.41 & 28.33
\\[3pt]
& \multicolumn{10}{c@{}}{Case SD}\\
$ (0,8 )$ & \phantom{0}6.19 & \phantom{0}6.74 & \phantom{0}6.35 & \phantom{0}6.38 & \phantom{0}5.92 & \phantom{0}4.20 & \phantom{0}5.69 & \phantom{0}3.98 & \phantom{0}4.00 & \phantom{0}3.96\\
$ (0,0.20 )$ & \phantom{0}8.95 & \phantom{0}8.86 & 10.40 & 10.45 & 10.04 & \phantom{0}8.81 & 10.91 & \phantom{0}9.14 & \phantom{0}9.06 & \phantom{0}8.91\\
$ (0.75,8 )$ & 10.50 & 14.23 & 10.34 & \phantom{0}9.81 & \phantom{0}8.72 & \phantom{0}5.19 & \phantom{0}7.31 & \phantom{0}4.74 & \phantom{0}4.99 & \phantom{0}5.13\\
$ (0.75,0.25 )$ & 19.90 & 21.25 & 22.46 & 23.48 & 24.58 & 19.10 & 24.36 & 18.51 & 19.45 & 19.07\\
\hline
\end{tabular*}
\end{table}

\begin{table}
\tabcolsep=0pt
\caption{MISE: high-dimensional slow decay, $K=100$}\label{tab5}
\begin{tabular*}{\tablewidth}{@{\extracolsep{\fill}}@{}lllllllllll@{}}
\hline
& \multicolumn{5}{l}{$n=20$} & \multicolumn{5}{l@{}}{$n=100$}\\[-6pt]
& \multicolumn{5}{l}{\hrulefill} & \multicolumn{5}{l@{}}{\hrulefill}\\
$(\omega,\sigma^{2} )$ & PGA & OGA & RGA & CGA & FWA & PGA & OGA & RGA & CGA & FWA\\
\hline
& \multicolumn{10}{c@{}}{Case ID}\\
$ (0,8 )$ & \phantom{0}0.97 & \phantom{0}0.94 & \phantom{0}0.92 & \phantom{0}0.93 & \phantom{0}1.00 & \phantom{0}0.42 & \phantom{0}0.51 & \phantom{0}0.46 & \phantom{0}0.42 & \phantom{0}0.42\\
$ (0,0.20 )$ & \phantom{0}1.34 & \phantom{0}1.95 & \phantom{0}2.05 & \phantom{0}1.67 & \phantom{0}1.69 & \phantom{0}1.01 & \phantom{0}0.95 & \phantom{0}1.05 & \phantom{0}0.97 & \phantom{0}0.99\\
$ (0.75,8 )$ & \phantom{0}1.10 & \phantom{0}1.56 & \phantom{0}1.08 & \phantom{0}1.07 & \phantom{0}1.08 & \phantom{0}0.51 & \phantom{0}0.77 & \phantom{0}0.53 & \phantom{0}0.56 & \phantom{0}0.56\\
$ (0.75,0.25 )$ & \phantom{0}2.28 & \phantom{0}3.22 & \phantom{0}3.50 & \phantom{0}2.70 & \phantom{0}3.03 & \phantom{0}1.54 & \phantom{0}1.71 & \phantom{0}1.68 & \phantom{0}1.47 & \phantom{0}1.50
\\[3pt]
& \multicolumn{10}{c@{}}{Case WD}\\
$ (0,8 )$ & \phantom{0}3.54 & \phantom{0}4.31 & \phantom{0}3.49 & \phantom{0}3.52 & \phantom{0}3.41 & \phantom{0}1.88 & \phantom{0}2.22 & \phantom{0}1.62 & \phantom{0}1.70 & \phantom{0}1.68\\
$ (0,0.20 )$ & \phantom{0}6.16 & \phantom{0}7.68 & \phantom{0}6.59 & \phantom{0}7.51 & \phantom{0}6.58 & \phantom{0}6.73 & \phantom{0}8.38 & \phantom{0}6.11 & \phantom{0}6.40 & \phantom{0}6.71\\
$ (0.75,8 )$ & \phantom{0}4.48 & \phantom{0}6.73 & \phantom{0}3.99 & \phantom{0}4.03 & \phantom{0}3.89 & \phantom{0}2.46 & \phantom{0}3.33 & \phantom{0}2.21 & \phantom{0}2.44 & \phantom{0}2.55\\
$ (0.75,0.25 )$ & 11.53 & 12.31 & 11.96 & 13.67 & 13.55 & 11.55 & 14.42 & 10.99 & 11.66 & 12.09
\\[3pt]
& \multicolumn{10}{c@{}}{Case SD}\\
$ (0,8 )$ & \phantom{0}2.83 & \phantom{0}3.40 & \phantom{0}2.72 & \phantom{0}2.76 & \phantom{0}2.74 & \phantom{0}1.45 & \phantom{0}1.81 & \phantom{0}1.36 & \phantom{0}1.37 & \phantom{0}1.35\\
$ (0,0.20 )$ & \phantom{0}5.04 & \phantom{0}4.81 & \phantom{0}5.82 & \phantom{0}5.89 & \phantom{0}5.60 & \phantom{0}4.81 & \phantom{0}5.90 & \phantom{0}4.85 & \phantom{0}4.84 & \phantom{0}4.69\\
$ (0.75,8 )$ & \phantom{0}3.37 & \phantom{0}5.08 & \phantom{0}3.24 & \phantom{0}3.41 & \phantom{0}3.22 & \phantom{0}1.97 & \phantom{0}2.60 & \phantom{0}1.75 & \phantom{0}1.80 & \phantom{0}1.82\\
$ (0.75,0.25 )$ & \phantom{0}8.51 & \phantom{0}8.82 & \phantom{0}9.57 & 10.07 & 10.51 & \phantom{0}8.12 & 10.00 & \phantom{0}7.98 & \phantom{0}8.11 & \phantom{0}8.16\\
\hline
\end{tabular*}
\end{table}

While somehow equivalent, the FWA updates the coefficients in a slightly
cruder way than the CGA. This seems to lead the FWA to have slightly
different performance than the CGA in some cases, with no definite
conclusion on which one is best. No attempt was made to use a line
search for $w_{j}$ (e.g., (\ref{EQRGALineSearch})) instead of the
deterministic weights.

\section{Proofs}\label{SectionProofs}

The proof for the results requires first to show that the estimators
nearly minimize the objective function $\llvert Y-\mu(X
)\rrvert _{n}^{2}$
for $\mu\in\mathcal{L} (B )$. Then uniform law of large
numbers for $\llvert Y-\mu(X )\rrvert _{n}^{2}$ with $\mu
\in\mathcal{L} (B )$
or related quantities are established.

To avoid cumbersome notation, for any functions of $ (Y,X )$,
say $f$ and $g$, write $ \langle f,g \rangle_{P}:=\int
f (y,x )g (y,x )\,\mathrm{d}P (y,x )$
where $P$ is the marginal distribution of $ (Y,X )$; moreover,
$\llvert f\rrvert _{P,2}^{2}:= \langle f,f \rangle_{P}$. In
the context of the paper, this means\vspace*{1pt} that $\llvert Y-\mu_{n}\rrvert
_{P,2}^{2}=\int\llvert y-\mu_{n} (x )\rrvert ^{2}\,\mathrm{d}P
(y,x )$
for a possibly random function $\mu_{n} (x )$ (e.g.,
a sample estimator). Clearly, if $\mu_{n}=\mu$ is not random, $\llvert
Y-\mu\rrvert _{P,2}^{2}=\llvert Y-\mu\rrvert _{2}^{2}$.
Consequently, the norm $\llvert \cdot\rrvert _{P,2}^{2}$ means that
$\llvert \mu_{n}-\mu\rrvert _{P,2}^{2}:=\mathbb{E}'\llvert \mu
_{n} (X' )-\mu(X' )\rrvert ^{2}$,
where $X'$ and $\mathbb{E}'$ are as defined just before (\ref
{EQerrorDecomposition}).

For any $\mu(X ):=\sum_{k=1}^{K}b_{k}X^{ (k
)}\in\mathcal{L}$,
$\llvert \mu\rrvert _{\mathcal{L}}=\sum_{k=1}^{K}\llvert
b_{k}\rrvert $
denotes the $l_{1}$ norm of the linear coefficients. Throughout,
$R_{m}:= (Y-F_{m} )$ denotes the residual in the approximation.

\subsection{Approximation rates for the algorithms}

The following provide approximation rates of the algorithms and show
that the resulting minimum converges to the global minimum, which
might not be unique, as the number of iterations $m$ goes to infinity.

%
%
\begin{lemma}\label{LemmaL2Boosting}For the PGA, for any $\mu\in
\mathcal{L} (B )$,
\[
\llvert R_{m}\rrvert_{n}^{2}\leq\bigl\llvert Y-
\mu(X )\bigr\rrvert_{n}^{2}+ \biggl(\frac{4\llvert Y\rrvert
_{n}^{4}B^{2}}{\nu
(2-\nu)m}
\biggr)^{1/3}.
\]
\end{lemma}

\begin{pf}Let $\tilde{R}_{0}=\mu\in\mathcal{L} (B )$,
and
\[
\tilde{R}_{m}=\tilde{R}_{m-1}-\nu\bigl\langle
X^{s (m
)},Y-F_{m-1} \bigr\rangle_{n}X^{s (m )}
\]
so that $\tilde{R}_{m}\in\mathcal{L} (B_{m} )$, where $B_{0}:=B$,
%
%
\begin{eqnarray}
B_{m} &:= & B_{m-1}+\nu\bigl\llvert\bigl\langle
X^{s (m
)},Y-F_{m-1} \bigr\rangle_{n}\bigr\rrvert
.\label{EQbmGrowth}
\end{eqnarray}
Also note that $\tilde{R}_{m}=R_{m}- (Y-\mu)$, where
$R_{m}=Y-F_{m}$,
$F_{0}=0$. Unlike $R_{0}$, $\tilde{R}_{0}$ has coefficients that
are controlled in terms of $B_{m}$, hence, it will be used to derive
a recursion for the gain at each greedy step. Hence, using these
remarks,
\[
\llvert\tilde{R}_{m}\rrvert_{n}^{2}= \langle
\tilde{R}_{m},\tilde{R}_{m} \rangle_{n}= \langle
\tilde{R}_{m},R_{m} \rangle_{n}- \langle
\tilde{R}_{m},Y-\mu\rangle_{n}\leq B_{m}\max
_{k}\bigl\llvert\bigl\langle X^{
(k )},R_{m}
\bigr\rangle_{n}\bigr\rrvert- \langle\tilde{R}_{m},Y-\mu
\rangle_{n}
\]
because $\tilde{R}_{m}\in\mathcal{L} (B_{m} )$, which, by
definition of $X^{s (m+1 )}$ implies
\begin{eqnarray*}
\bigl\llvert\bigl\langle X^{ (m+1 )},R_{m} \bigr\rangle
_{n}\bigr\rrvert& \geq& \frac{ \langle\tilde{R}_{m},\tilde
{R}_{m}+Y-\mu\rangle_{n}}{B_{m}}=\frac{ \langle\tilde
{R}_{m},R_{m} \rangle_{n}}{B_{m}}
\\
& = & \frac{ \langle R_{m},R_{m} \rangle_{n}- \langle
R_{m},Y-\mu\rangle_{n}}{B_{m}}
\end{eqnarray*}
using the definition of $\tilde{R}_{m}$ in the last equality. Then,
by the scalar inequality $ab\leq(a^{2}+b^{2} )/2$ the above
becomes
%
%
\begin{equation}
\bigl\llvert\bigl\langle X^{ (m+1 )},R_{m} \bigr\rangle
_{n}\bigr\rrvert\geq\frac{\llvert R_{m}\rrvert _{n}^{2}-\llvert
Y-\mu
\rrvert _{n}^{2}}{2B_{m}}.\label{EQgainLowerBound}
\end{equation}
Note that the right-hand side is positive, if not, $\llvert
Y-F_{m}\rrvert _{n}^{2}\leq\llvert Y-\mu\rrvert _{n}^{2}$
and the lemma is proved (recall that $R_{m}=Y-F_{m}$). Now, note
that $R_{m}=R_{m-1}-\nu\langle X^{s (m
)},R_{m-1} \rangle_{n}X^{s (m )}$,
so that
\begin{eqnarray*}
\llvert R_{m}\rrvert_{n}^{2} & = & \llvert
R_{m-1}\rrvert_{n}^{2}+\nu^{2}\bigl
\llvert\bigl\langle X^{s (m )},R_{m-1} \bigr\rangle
_{n}\bigr\rrvert^{2}-2\nu\bigl\llvert\bigl\langle
X^{s (m
)},R_{m-1} \bigr\rangle_{n}\bigr\rrvert
^{2}
\\
& = & \llvert R_{m-1}\rrvert_{n}^{2}-\nu(2-\nu)
\bigl\llvert\bigl\langle X^{s (m )},R_{m-1} \bigr
\rangle_{n}\bigr\rrvert^{2}.
\end{eqnarray*}
The above two displays imply
\begin{eqnarray*}
\llvert R_{m}\rrvert_{n}^{2} & \leq& \llvert
R_{m-1}\rrvert_{n}^{2}-\frac{\nu(2-\nu)}{4B_{m-1}^{2}} \bigl(
\llvert R_{m-1}\rrvert_{n}^{2}-\llvert Y-\mu\rrvert
_{n}^{2} \bigr)^{2}.
\end{eqnarray*}
Subtracting $\llvert Y-\mu\rrvert _{n}^{2}$ on both sides, and defining
$a_{m}:=\llvert R_{m}\rrvert _{n}^{2}-\llvert Y-\mu\rrvert
_{n}^{2}$, and
$\tau:=\nu(2-\nu)/4$, the above display is
%
%
\begin{equation}
a_{m}\leq a_{m-1} \bigl(1-\tau a_{m-1}B_{m-1}^{-2}
\bigr).\label
{EQamRecursion}
\end{equation}
The proof then exactly follows the proof of Theorem 3.6 in DeVore and
Temlyakov \cite{DeVTem96}. For completeness, the details are provided.
Define
%
%
\begin{eqnarray}
\rho(R_{m} ) &:= & a_{m}^{-1/2}\bigl\llvert
\bigl\langle X^{s (m+1 )},R_{m} \bigr\rangle_{n}\bigr
\rrvert\geq a_{m}^{1/2}B_{m}^{-1}.\label{EQlowerGainRhoRepresentation}
\end{eqnarray}
Since $B_{m}\geq B_{m-1}$,
%
%
\begin{eqnarray}
a_{m}B_{m}^{-2} & \leq& a_{m-1}B_{m-1}^{-2}
\bigl(1-\tau a_{m-1}B_{m-1}^{-2} \bigr)\leq
\frac{1}{\tau m}\label{EQamBmDivisionBound}
\end{eqnarray}
using Lemma 3.4 in DeVore and Temlyakov \cite{DeVTem96} in the
second step
in order to bound the recursion. Then (\ref{EQbmGrowth}) and (\ref
{EQlowerGainRhoRepresentation})
give
\begin{eqnarray*}
B_{m} & = & B_{m-1} \bigl(1+\nu\rho(R_{m-1}
)a_{m-1}^{1/2}B_{m-1}^{-1} \bigr)
\\
& \leq& B_{m-1} \bigl(1+\nu\rho(R_{m-1} )^{2}
\bigr).
\end{eqnarray*}
Multiply both sides of (\ref{EQamRecursion}) by $B_{m}$, and substitute
the lower bound (\ref{EQlowerGainRhoRepresentation}) into (\ref
{EQamRecursion}),
so that using the above display,
\begin{eqnarray*}
a_{m}B_{m} & \leq& a_{m-1}B_{m-1}
\bigl(1+\nu\rho(R_{m-1} )^{2} \bigr) \bigl(1-\tau\rho
(R_{m-1} )^{2} \bigr)
\\
& = & a_{m-1}B_{m-1} \bigl(1-\nu\tau\rho(R_{m-1}
)^{4} \bigr)\leq\llvert Y\rrvert_{n}^{2}B,
\end{eqnarray*}
where the last inequality follows after iterating because $1-\nu\tau
\rho(R_{m-1} )^{4}\in(0,1 )$
and substituting $B_{0}=B$ and $a_{0}=\llvert Y\rrvert _{n}^{2}$. If
$a_{m}>0$, it is obvious that $1-\nu\tau\rho(R_{m-1}
)^{4}\in(0,1 )$.
If this were not the case, the lemma would hold automatically at step
$m$, by definition of $a_{m}$. Hence, by the above display together
with (\ref{EQamBmDivisionBound}),
\[
a_{m}^{3}= (a_{m}B_{m}
)^{2}a_{m}B_{m}^{-2}\leq
\frac
{4\llvert Y\rrvert _{n}^{4}B^{2}}{\nu(2-\nu)m}
\]
using the definition of $\tau=\nu(2-\nu)/4$, so that
$a_{m}\leq[4\llvert Y\rrvert _{n}^{4}B/ (\nu(2-\nu
)m ) ]^{1/3}$.
\end{pf}

The following bound for the OGA is Theorem 2.3 in Barron \textit
{et~al.} \cite{Baretal08}.

%
%
\begin{lemma}\label{LemmaOGA}For the OGA, for any $\mu\in\mathcal
{L} (B )$,
\[
\llvert R_{m}\rrvert_{n}^{2}\leq\bigl\llvert Y-
\mu(X )\bigr\rrvert_{n}^{2}+4\frac{B^{2}}{m}.
\]
\end{lemma}

The following Lemma~\ref{Lemmarelaxedgreedy} is Theorem 2.4 in
Barron \textit{et~al.} (\cite{Baretal08}, equation (2.41)),
where the CGA bound is inferred from
their proof (in their proof set their $\beta$ on page 78 equal to
$w_{k}\bar{B}$ to satisfy the CGA constraint).

%
%
\begin{lemma}\label{Lemmarelaxedgreedy}For the RGA, for any $\mu\in
\mathcal{L} (B )$,
\[
\llvert R_{m}\rrvert_{n}^{2}\leq\bigl\llvert Y-
\mu(X )\bigr\rrvert_{n}^{2}+\frac{B^{2}}{m}.
\]
For the CGA the above holds with $B$ replaced by $\bar{B}$ in the
above display and any $\mu\in\mathcal{L} (\bar{B} )$.
\end{lemma}

%
%
\begin{lemma}\label{LemmafrankWolfeApproximation}For the FWA, for
any $\mu\in\mathcal{L} (\bar{B} )$, and $m>0$,
\[
\llvert R_{m}\rrvert_{n}^{2}\leq\bigl\llvert Y-
\mu(X )\bigr\rrvert_{n}^{2}+\frac{4\bar{B}^{2}}{m},
\]
when $w_{m}=2/ (1+m )$.
\end{lemma}

\begin{pf}From Jaggi (\cite{Jag13}, equations (3)--(4),
see also Frank and Wolfe \cite{FraWol56}), for every
$m=1,2,3,\ldots,$ infer
the first inequality in the
following display:
\begin{eqnarray*}
\llvert R_{m}\rrvert_{n}^{2}-\bigl\llvert Y-\mu
(X )\bigr\rrvert_{n}^{2} & \leq& (1-w_{m} )
\bigl(\llvert R_{m-1}\rrvert_{n}^{2}-\bigl\llvert
Y-\mu(X )\bigr\rrvert_{n}^{2} \bigr)
\\
& &{} +w_{m}^{2}\max_{\sum_{k=1}^{K}\llvert b_{k}\rrvert \leq\bar
{B},\sum_{k=1}^{K}\llvert c_{k}'\rrvert \leq\bar{B}}\Biggl\llvert
\sum_{k=1}^{K} \bigl(b_{k}-b_{k}'
\bigr)X^{ (k )}\Biggr\rrvert_{n}^{2}
\\
& \leq& (1-w_{m} ) \bigl(\llvert R_{m-1}\rrvert
_{n}^{2}-\bigl\llvert Y-\mu(X )\bigr\rrvert
_{n}^{2} \bigr)+w_{m}^{2}4
\bar{B}^{2}\max_{k\leq K}\bigl\llvert X^{ (k
)}
\bigr\rrvert_{n}^{2},
\end{eqnarray*}
where the second inequality follows because the maximum over the simplex
is at one of the edges of the simplex. Moreover, $\max_{k\leq
K}\llvert X^{ (k )}\rrvert _{n}^{2}=1$
by construction. The result then follows by Theorem 1 in Jaggi \cite{Jag13}
when $w_{m}=2/ (1+m )$.
\end{pf}

\subsection{Size of the functions generated by the algorithms}

The following gives a bound for the size of $F_{m}$ in terms of the
norm $\llvert \cdot\rrvert _{\mathcal{L}}$; $F_{m}$ is the function
generated by each algorithm.

%
%
\begin{lemma}\label{LemmaEstimatorL1Bound}As $n\rightarrow\infty$,
$\Pr(F_{m}\in\mathcal{L} (B_{m} )
)\rightarrow1$,
where:

\begin{longlist}[RGA:]
\item[PGA:] $B_{m}\lesssim\llvert Y\rrvert _{2}m^{1/2}$;

\item[OGA:] $B_{m}\lesssim\llvert Y\rrvert _{2} [ (\frac{m}{\rho
_{m,n}} )^{1/2}\wedge m\wedge K ]$
with $\rho_{m,n}$ as in (\ref{EQrestrictedEigenvalueCondition});

\item[RGA:] $B_{m}\lesssim\llvert Y\rrvert _{2} [ (\frac{m}{\rho
_{m,n}} )^{1/2}\wedge m\wedge K ]$
with $\rho_{m,n}$ as in (\ref{EQrestrictedEigenvalueCondition}),
as long as in Lemma~\ref{Lemmarelaxedgreedy} $B^{2}/m=\mathrm{O} (1)$;

\item[CGA and FWA:] $B_{m}\leq\bar{B}$.
\end{longlist}
\end{lemma}

\begin{pf}Note that $F_{m} (X )=\sum_{k=1}^{m}b_{k}X^{s (k )}$,
where to ease notation $b_{k}$ does not make explicit the dependence
on $m$. A loose bound for $\llvert F_{m}\rrvert _{\mathcal{L}}$ is
found by noting that
\[
\bigl\llvert\bigl\langle X^{ (k )},X^{ (l )} \bigr
\rangle_{n}\bigr\rrvert\leq\max_{k}\bigl\llvert
\bigl\langle X^{
(k )},X^{ (k )} \bigr\rangle_{n}\bigr
\rrvert=1,
\]
so that each coefficient is bounded by $\llvert Y\rrvert _{n}$. Since
at the $m$th iteration we have at most $m$ different terms and
no more than $K$, $\llvert F_{m}\rrvert _{\mathcal{L}}\leq
(m\wedge K )\llvert Y\rrvert _{n}$.
Given that $\llvert Y\rrvert _{n}^{2}=\mathrm{O}_{p} (1 )$, one can
infer the crude bound $\llvert F_{m}\rrvert _{\mathcal{L}}=\mathrm{O}_{p}
(m\wedge K )$.
This is the worse case scenario, and can be improved for all the algorithms.

For the PGA, at the first iteration, $\llvert b_{1}\rrvert:=\max
_{k}\llvert \langle X^{ (k )},Y \rangle
_{n}\rrvert \leq\llvert Y\rrvert _{n}$,
hence there is an $\alpha_{1}\in[0,1 ]$ such that $\llvert
b_{1}\rrvert =\alpha_{1}^{1/2}\llvert Y\rrvert _{n}$
(the root exponent is used to ease notation in the following steps).
Then, by the properties of projections
\[
\llvert R_{1}\rrvert_{n}^{2}=\bigl\llvert
Y-X^{s (1 )}b_{1}\bigr\rrvert_{n}^{2}=
\llvert Y\rrvert_{n}^{2}-\llvert b_{1}\rrvert
^{2}=\llvert Y\rrvert_{n}^{2} (1-
\alpha_{1} ),
\]
where the second inequality follows from $\llvert X^{ (k
)}\rrvert _{n}^{2}=1$
for any $k$. By similar arguments, there is an $\alpha_{2}\in
[0,1 ]$
such that $\llvert b_{2}\rrvert =\alpha_{2}^{1/2}\llvert
R_{1}\rrvert _{n}$
and $\llvert R_{2}\rrvert _{n}^{2}=\llvert R_{1}\rrvert _{n}^{2}
(1-\alpha_{2} )$.
So by induction $\llvert b_{m}\rrvert =\alpha_{m}^{1/2}\llvert
R_{m-1}\rrvert _{n}$
and $\llvert R_{m}\rrvert _{n}^{2}=\llvert R_{m-1}\rrvert _{n}^{2}
(1-\alpha_{m} )$.
By recursion, this implies that
\[
\llvert b_{m}\rrvert^{2}=\alpha_{m} (1-
\alpha_{m-1} )\cdots(1-\alpha_{1} )\llvert Y\rrvert
_{n}^{2}
\]
and in consequence that
\begin{eqnarray*}
\sum_{k=1}^{m}\llvert b_{k}
\rrvert& = & \sum_{k=1}^{m}\alpha
_{k}^{1/2}\prod_{l<k} (1-
\alpha_{l} )^{1/2}\llvert Y\rrvert_{n},
\end{eqnarray*}
where the empty product is $1$. It is clear that if any $\alpha_{k}\in
\{ 0,1 \} $
for $k<m$ then $b_{m}=0$, hence one can assume that all the $\alpha_{k}$'s
are in $ (0,1 )$. The above display is maximized if $\alpha
_{l}\rightarrow0$
fast enough, as otherwise, the product converges to zero exponentially
fast and the result follows immediately. Suppose that $\sum
_{l=1}^{\infty}\alpha_{l}^{2}<\infty$.
Then, using the fact that $\ln(1-\alpha_{l} )=-\alpha
_{l}+\mathrm{O} (\alpha_{l}^{2} )$,
\begin{eqnarray*}
\prod_{l<k} (1-\alpha_{l} ) & = &\cdots  =\exp
\Biggl\{ \sum_{l=1}^{k-1}\ln(1-
\alpha_{l} ) \Biggr\} =\exp\Biggl\{ -\sum
_{l=1}^{k-1}\alpha_{l}+\mathrm{O} \Biggl(\sum
_{l=1}^{k-1}\alpha_{l}^{2}
\Biggr) \Biggr\}
\\
& \asymp& \exp\Biggl\{ -\sum_{l=1}^{k-1}
\alpha_{l} \Biggr\}.
\end{eqnarray*}
The above converges exponentially fast to $0$ if $\alpha_{l}\asymp
l^{-\alpha}$
for $\alpha\in(0.5,1 )$. While the argument is not valid
for $\alpha\in(0,0.5]$, it is clear, that the convergence is even
faster in this case. Hence, restrict attention to $\alpha=1$, in
which case, $\prod_{l<k} (1-\alpha_{l} )\asymp k^{-c}$ for
some $c>0$, that is, polynomial decay. On the other hand for $\alpha>1$,
the product converges. Hence, it must be the case that the maximum
is achieved by setting $\alpha_{l}\asymp l^{-1}$ and assuming that
the product converges. This implies that for the PGA,
\[
\sum_{k=1}^{m}\llvert b_{k}
\rrvert\lesssim\llvert Y\rrvert_{n}\sum_{k=1}^{m}
\bigl(k^{-1} \bigr)^{1/2}\lesssim\llvert Y\rrvert
_{n}m^{1/2}.
\]
Now, consider the OGA and the RGA. The following just follows by standard
inequalities:
%
%
\begin{equation}
(\rho_{m,n}/m )^{1/2}\sum_{k=1}^{m}
\llvert b_{k}\rrvert\leq\rho_{m,n}^{1/2} \Biggl(
\sum_{k=1}^{m}\llvert b_{k}
\rrvert^{2} \Biggr)^{1/2}\leq\llvert F_{m}\rrvert
_{n}.\label{EQstandardInequality}
\end{equation}
For the OGA, by definition of the OLS estimator, $\llvert
F_{m}\rrvert _{n}\leq\llvert Y\rrvert _{n}$
implying the result for the OGA using the above display and the crude
bound. For the RGA, consider the case when $\llvert F_{m}\rrvert _{n}$
is small and large, separately. If $\llvert F_{m}\rrvert _{n}=\mathrm{o}_{p}
(1 )$,
then clearly, $\llvert F_{m}\rrvert _{n}=\mathrm{o} (\llvert Y\rrvert _{n} )$,
because $Y$ is not degenerate. By this remark, the above display
implies that
\[
\llvert F_{m}\rrvert_{\mathcal{L}}:=\sum
_{k=1}^{m}\llvert b_{k}\rrvert
=\mathrm{o}_{p} \bigl(\sqrt{m/\rho_{m,n}}\llvert Y\rrvert
_{n} \bigr)
\]
and the result for the RGA would follow. Hence, one can assume that
$\llvert F_{m}\rrvert _{n}\gtrsim1$ in probability, eventually as
$m\rightarrow\infty$.
In this case, by the approximating Lemma~\ref{Lemmarelaxedgreedy},
if $B^{2}/m=\mathrm{O} (1 )$,
\[
\llvert Y-F_{m}\rrvert_{n}^{2}\leq\llvert Y
\rrvert_{n}^{2}+\mathrm{O} (1 )
\]
which implies
\[
\llvert F_{m}\rrvert_{n}^{2}\leq2 \langle
Y,F_{m} \rangle_{n}+\mathrm{O} (1 )\leq2\llvert Y\rrvert
_{n}\llvert F_{m}\rrvert_{n}+\mathrm{O} (1 )
\]
and in consequence
\[
\llvert F_{m}\rrvert_{n}\leq2\llvert Y\rrvert
_{n}+\mathrm{O} \bigl(\llvert F_{m}\rrvert_{n}^{-1}
\bigr)=2\llvert Y\rrvert_{n}+\mathrm{O}_{p} (1 )
\]
by the fact that $\llvert F_{m}\rrvert _{n}\gtrsim1$, in probablity.
Hence, using the above display together with (\ref{EQstandardInequality}),
the result follows for the RGA as well.

For the CGA, the $b_{k}$'s are all bounded in absolute value by $\bar{B}$.
Since by construction, $F_{m} (X )=m^{-1}\sum_{k=1}^{m}b_{k}X^{s (k
)}$, $\llvert F_{m}\rrvert _{\mathcal
{L}}\leq\bar{B}$.
A similar argument holds for the FWA.
\end{pf}

It is natural to replace the random eigenvalue $\rho_{m,n}$ with
the population one. This is achieved next.

%
%
\begin{lemma}\label{LemmaeigenValue}Suppose\vspace*{1pt} Conditions~\ref{ConditionEY|X}
and~\ref{Conditiondependence} hold. Then $\rho_{m,n}\geq\rho
_{m}-\mathrm{O}_{p} (d_{n,p}mK^{2/p}n^{-1/2} )$
implying that if $d_{n,p}mK^{2/p}n^{-1/2}=\mathrm{o} (\rho_{m} )$,
then $\rho_{m,n}^{-1}=\mathrm{O}_{p} (\rho_{m}^{-1} )$.
\end{lemma}

\begin{pf}Note that
\[
\rho_{m,n}=\inf_{\llvert b\rrvert _{0}\leq m,\llvert b\rrvert
_{2}\leq
1}\frac{1}{n}\sum
_{i=1}^{n} \Biggl(\sum
_{k=1}^{K}b_{k}X_{i}^{
(k )}
\Biggr)^{2},\qquad\rho_{m}=\inf_{\llvert b\rrvert _{0}\leq
m,\llvert b\rrvert _{2}\leq1}
\frac{1}{n}\sum_{i=1}^{n}\mathbb{E}
\Biggl(\sum_{k=1}^{K}b_{k}X_{i}^{ (k )}
\Biggr)^{2},
\]
where $\llvert b\rrvert _{0}=\sum_{k=1}^{K} \{ b_{k}\neq0
\} $
and $\llvert b\rrvert _{2}^{2}=\sum_{k=1}^{K}\llvert b_{k}\rrvert ^{2}$,
that is, the number of non-zero $b_{k}$'s and their squared $l_{2}$
norm, respectively. By obvious manipulations, using the above display,
and the definition of $\rho_{m}$,
\begin{eqnarray*}
\rho_{m,n} & \geq& \rho_{m}-\sup_{\llvert b\rrvert _{0}\leq
m,\llvert b\rrvert _{2}\leq1}
\Biggl\llvert\frac{1}{n}\sum_{i=1}^{n}
(1-\mathbb{E} ) \Biggl(\sum_{k=1}^{K}b_{k}X_{i}^{ (k )}
\Biggr)^{2}\Biggr\rrvert,
\end{eqnarray*}
hence it is sufficient to bound the r.h.s. of the above display. Using
similar arguments as in the control of $\mathit{II}$ in the proof of Lemma
\ref{Lemmaulln} in Section~\ref{SectionuniformControl},
\begin{eqnarray*}
& & \mathbb{E}\sup_{\llvert b\rrvert _{0}\leq m,\llvert b\rrvert
_{2}\leq1}\Biggl\llvert\frac{1}{n}\sum
_{i=1}^{n} (1-\mathbb{E} ) \Biggl(\sum
_{k=1}^{K}b_{k}X_{i}^{ (k )}
\Biggr)^{2}\Biggr\rrvert
\\
&&\quad \leq m\mathbb{E}\max_{k,l\leq K}\Biggl\llvert\frac{1}{n}
\sum_{i=1}^{n} (1-\mathbb{E}
)X_{i}^{ (k
)}X_{i}^{ (l )}\Biggr\rrvert
\lesssim\frac
{d_{n,p}mK^{2/p}}{\sqrt{n}},
\end{eqnarray*}
and the first result follows. The second part is directly inferred
from the first.
\end{pf}

\subsection{Inequalities for dependent random variables}

Two different inequalities will be needed depending on whether one
assumes absolute regularity or mixingales. The following is suitable
for beta mixing random variables. It is somewhat standard, but proved
for completeness due to some adjustments to the present context.

%
%
\begin{lemma}\label{LemmaMaximalInequ}Suppose that $\mathfrak{F}$
is a measurable class of functions with cardinality $K$. Let $
(W_{i} )_{i\in\mathbb{Z}}$
be strictly stationary and beta mixing with mixing coefficients $\beta
(i )\lesssim\beta^{i}$,
$\beta\in[0,1)$. Suppose that for all $f\in\mathfrak{F}$, $\mathbb
{E}\llvert f (W_{1} )\rrvert ^{p}<\infty$
for some $p>2$. Then
\begin{eqnarray*}
\mathbb{E}\max_{f\in\mathfrak{F}}\frac{1}{\sqrt{n}}\Biggl\llvert\sum
_{i=1}^{n} (1-\mathbb{E} )f (W_{i}
)\Biggr\rrvert& \lesssim& \sqrt{\ln K}
\end{eqnarray*}
if $K\lesssim n^{\alpha}$ for some $\alpha< (p-2 )/2$.
If $\max_{f\in\mathfrak{F}}\llvert f\rrvert $ is bounded, the result
holds for $K\lesssim\exp\{ n^{\alpha} \} $, $\alpha\in[0,1)$.
\end{lemma}

\begin{pf}Note that
\begin{eqnarray*}
& & \mathbb{E}\max_{f\in\mathfrak{F}}\frac{1}{\sqrt{n}}\Biggl\llvert
\sum
_{i=1}^{n} (1-\mathbb{E} )f (W_{i}
)\Biggr\rrvert
\\[-2pt]
&&\quad \leq \mathbb{E}\max_{f\in\mathfrak{F}}\frac{1}{\sqrt{n}}\Biggl\llvert
\sum_{i=1}^{n} (1-\mathbb{E} )f
(W_{i} ) \Bigl\{ \max_{f\in\mathfrak{F}}\bigl\llvert f
(W_{i} )\bigr\rrvert\leq M \Bigr\} \Biggr\rrvert
\\[-2pt]
&&\qquad{} +\mathbb{E}\max_{f\in\mathfrak{F}}\frac{1}{\sqrt{n}}\Biggl\llvert
\sum
_{i=1}^{n} (1-\mathbb{E} )f (W_{i}
) \Bigl\{ \max_{f\in\mathfrak{F}}\bigl\llvert f (W_{i} )\bigr
\rrvert>M \Bigr\} \Biggr\rrvert
\\[-2pt]
&&\quad \leq \mathbb{E}\max_{f\in\mathfrak{F}}\frac{1}{\sqrt{n}}\Biggl\llvert
\sum_{i=1}^{n} (1-\mathbb{E} )f
(W_{i} ) \Bigl\{ \max_{f\in\mathfrak{F}}\bigl\llvert f
(W_{i} )\bigr\rrvert\leq M \Bigr\} \Biggr\rrvert
\\[-2pt]
&&\qquad{} +2\sqrt{n}\mathbb{E}\max_{f\in\mathfrak{F}}\bigl\llvert f
(W_{i} )\bigr\rrvert\Bigl\{ \max_{f\in\mathfrak{F}}\bigl\llvert
f (W_{i} )\bigr\rrvert>M \Bigr\}
\\[-2pt]
&&\quad =:  I+\mathit{II},
\end{eqnarray*}
where in the last inequality one uses Minkowski's inequality. (Here,
$ \{ \cdot\} $ is the indicator of a set.) By H\"older's
inequality,
\begin{eqnarray*}
\mathit{II} & \leq& 2\sqrt{n} \Bigl(\mathbb{E}\max_{f\in\mathfrak{F}}\bigl
\llvert f
(W_{i} )\bigr\rrvert^{p} \Bigr)^{1/p}\Pr\Bigl(
\max_{f\in
\mathfrak{F}}\bigl\llvert f (W_{i} )\bigr\rrvert>M
\Bigr)^{
(p-1 )/p}
\\
& \leq& 2\sqrt{n} \Bigl(\mathbb{E}\max_{f\in\mathfrak{F}}\bigl\llvert f
(W_{i} )\bigr\rrvert^{p} \Bigr)^{1/p}K^{ (p-1
)/p}M^{- (p-1 )}
\\
& \lesssim& \sqrt{n}KM^{- (p-1 )}
\end{eqnarray*}
because by Markov inequality and the union bound,
\[
\Pr\Bigl(\max_{f\in\mathfrak{F}}\bigl\llvert f (W_{i} )\bigr
\rrvert>M \Bigr)\lesssim KM^{-p},
\]
while $\mathbb{E}\max_{f\in\mathfrak{F}}\llvert f (W_{i}
)\rrvert ^{p}\lesssim K$
(e.g., Lemma 2.2.2 in van der Vaart and Wellner \cite{vanWel00}). Hence, set
$M= (\sqrt{n}K/\sqrt{\ln K} ){}^{1/ (p-1 )}$ to
ensure that $\mathit{II}=\mathrm{O} (\sqrt{\ln K} )$. Pollard (\cite{Pol02}, equation~(8))
shows that if the $W_{i}$'s are beta mixing, for any integer sequence
$a_{n}=\mathrm{o} (n )$,
%
%
\begin{equation}
I\lesssim\sqrt{\ln K}\llvert f\rrvert_{2\beta}\mathcal{E} \biggl(
\frac{Ma_{n}\sqrt{2\ln K}}{\llvert f\rrvert _{2\beta}\sqrt
{n}} \biggr)+M\beta(a_{n} )\sqrt{n},\label{EQmaximalInequalityBennetStep}
\end{equation}
where $\mathcal{E}$ is some positive increasing function such that
$\lim_{x\rightarrow\infty}\mathcal{E} (x )=\infty$ and
$\llvert \cdot\rrvert _{2\beta}$
is the beta mixing norm introduced by Doukhan \textit{et~al.} \cite{DouMasRio95}
(see also
Rio \cite{Rio00}, equation (8.21)). The exact form
of the norm is
irrelevant for
the present purposes, however, $\llvert f\rrvert _{2\beta}\leq
c_{1}<\infty$
for some constant $c_{1}$ under the condition on the mixing coefficients
(e.g., Rio \cite{Rio00}, p. 15). Since $\beta
(a_{n}
)\lesssim\beta^{a_{n}}$,
for $a_{n}\asymp\ln n/\ln(1/\beta)$, and using the value
for $M$ set in $\mathit{II}$, deduce
\[
I+\mathit{II}\lesssim\sqrt{\ln K}\mathcal{E} \biggl(\frac{c_{2}\ln n
(\sqrt{n}K/\sqrt{\ln K} ){}^{1/ (p-1 )}\sqrt{\ln
K}}{\sqrt{n}} \biggr)+\sqrt{\ln
K},
\]
for some finite positive constant $c_{2}$. Substituting $K\asymp
n^{\alpha}$
for any positive $\alpha< (p-2 )/2$, the argument in the
continuous increasing function $\mathcal{E} (\cdot)$
is bounded and the result follows. Notice that this choice of $K$
also makes $M\beta(a_{n} )\sqrt{n}\lesssim1$.

For the case of bounded $\max_{f\in\mathfrak{F}}\llvert f\rrvert $,
one can take $M$ large enough, but finite so that $\mathit{II}=0$. Given
that $M$ is finite, $K\lesssim\exp\{ n^{\alpha} \} $,
and with $a_{n}$ as before, (\ref{EQmaximalInequalityBennetStep})
becomes
\[
I\lesssim\sqrt{\ln K}\llvert f\rrvert_{2\beta}\mathcal{E} \biggl(
\frac{c_{3}\ln n\sqrt{n^{\alpha}}}{\sqrt{n}} \biggr)
\]
for some finite constant $c_{3}$, and the argument of $\mathcal
{E} (\cdot)$
is bounded because $\alpha<1$. Some tidying up gives the last result.
\end{pf}

The following is an extension of Burkh\"older inequality to mixingales
(see Peligrad, Utev and Wu \cite{PelUteWu07}, Corollary~1 for the exact constants).

%
%
\begin{lemma}\label{LemmaBurkholderInequality}Suppose that $
(W_{i} )_{i\in\mathbb{Z}}$
is a mean zero stationary sequence of random variables. Let
\[
d_{n,p} (W ):=\sum_{i=0}^{n}
(i+1 )^{-1/2}\bigl\llvert\mathbb{E} [W_{i}|
\mathcal{F}_{0} ]\bigr\rrvert_{p},
\]
where $\mathcal{F}_{0}:=\sigma(W_{i}\dvt i\leq0 )$ is the sigma
algebra generated by $ (W_{i}\dvt i\leq0 )$. Then, for all
$p\geq2$,
such that $\llvert W_{i}\rrvert _{p}<\infty$,
\[
\Biggl\llvert\sum_{i=1}^{n}W_{i}
\Biggr\rrvert_{p}\leq C_{p}^{1/p}n^{1/2}d_{n,p}
(W ),
\]
where for $p\in[2,4)$, $C_{p}\lesssim p^{p}$ while for $p\geq4$,
$C_{p}\lesssim(2p )^{p/2}$.
\end{lemma}

\subsection{Uniform control of the estimator}\label{SectionuniformControl}

Next, one needs a uniform control of the objective function. Recall
that $\mu_{B}$ is the best approximation in $\mathcal{L}
(B )$
to $\mu_{0}$ in the $L_{2}$ sense.

Define
\[
\mathcal{L}_{0} (B ):= \Biggl\{ \mu\dvt \mu(X )=\sum
_{k=1}^{K}b_{k}X^{ (k )},\sum
_{k=1}^{K} \{ b_{k}\neq0 \} \leq B
\Biggr\}.
\]
These are linear functions with $l_{0}$ norm less or equal to $B$,
that is, linear functions with at most $B$ non-zero coefficients. The
following is Lemma 5.1 in van~de Geer \cite{van14} with
minor differences.
The proof is given for completeness.

%
%
\begin{lemma}\label{LemmavanDeGeerConsistency}Let $\mu'\in\mathcal
{L} (B )$
be an arbitrary but fixed function and $m$ a positive integer. Suppose
that in probability, for some $\delta_{1}\in(0,1 )$ and
$\delta_{2},\delta_{3}>0$:

\begin{enumerate}
\item $\sup_{\mu\in\mathcal{L}_{0} (2m )\dvt \llvert \mu\rrvert _{2}\leq
1}\llvert
(1-\mathbb{E} )\llvert \mu\rrvert _{n}^{2}\rrvert \leq
\delta_{1}$,

\item $\sup_{\mu\in\mathcal{L}_{0} (2m )\dvt \llvert \mu\rrvert _{2}\leq
1}\llvert 2 (1-\mathbb{E} ) \langle Y-\mu
',\mu\rangle_{n}\rrvert \leq\delta_{2}$,

\item the sequence $\mu_{n}\in\mathcal{L}_{0} (m )$ satisfies
$\llvert Y-\mu_{n}\rrvert _{n}^{2}\leq\llvert Y-\mu'\rrvert
_{n}^{2}+\delta_{3}^{2}$,

\item the moment condition $ \langle Y-\mu',\mu_{n} \rangle_{P}=0$
holds.
\end{enumerate}
Then $\llvert \mu_{n}-\mu'\rrvert _{P,2}\leq(\delta_{2}+\delta
_{3} )/ (1-\delta_{1} )$
in probability (recall the definition of $\llvert \cdot\rrvert _{P,2}$
at the beginning of Section~\ref{SectionProofs}).
\end{lemma}

\begin{pf}Starting from the assumption
\[
\llvert Y-\mu_{n}\rrvert_{n}^{2}\leq\bigl\llvert
Y-\mu'\bigr\rrvert_{n}^{2}+
\delta_{3}^{2},
\]
by algebraic manipulations, $\llvert \mu_{n}-\mu'\rrvert _{n}^{2}\leq
2 \langle Y-\mu',\mu_{n}-\mu' \rangle_{n}+\delta_{3}^{2}$.
Assume that $\llvert \mu_{n}-\mu'\rrvert _{P,2}\geq\delta_{3}$ otherwise,
there is nothing to prove. Hence, $\delta_{3}^{2}\leq\delta_{3}\llvert
\mu_{n}-\mu'\rrvert _{P,2}$.
Also note that $ \langle Y-\mu',\mu_{n}-\mu' \rangle_{P}=0$
by definition of $\mu'$ (point 4 in the statement). Adding and subtracting
$\llvert \mu_{n}-\mu'\rrvert _{n}^{2}$, and using the just derived bounds
\begin{eqnarray*}
\bigl\llvert\mu_{n}-\mu'\bigr\rrvert
_{P,2}^{2}
&\leq& \bigl\llvert\mu_{n}-\mu'\bigr\rrvert
_{P,2}^{2}-\bigl\llvert\mu_{n}-\mu'
\bigr\rrvert_{n}^{2}+2 \bigl( \bigl\langle Y-
\mu',\mu_{n}-\mu' \bigr\rangle_{n}-
\bigl\langle Y-\mu',\mu_{n}-\mu' \bigr\rangle
_{P} \bigr)
\\
&&{} +2 \bigl\langle Y-\mu',\mu_{n}-\mu'
\bigr\rangle_{P}+\delta_{3}^{2}
\\
&\leq& \biggl\llvert\frac{\llvert \mu_{n}-\mu'\rrvert
_{P,2}^{2}-\llvert \mu_{n}-\mu'\rrvert _{n}^{2}}{\llvert \mu_{n}-\mu
'\rrvert _{P,2}^{2}}\biggr\rrvert\bigl\llvert
\mu_{n}-\mu'\bigr\rrvert_{P,2}^{2}
\\
&&{} +2\biggl\llvert\frac{ \langle Y-\mu',\mu_{n}-\mu' \rangle
_{n}- \langle Y-\mu',\mu_{n}-\mu' \rangle_{P}}{\llvert \mu_{n}-\mu
'\rrvert _{P,2}}\biggr\rrvert\bigl\llvert
\mu_{n}-\mu'\bigr\rrvert_{P,2}+
\delta_{3}\bigl\llvert\mu_{n}-\mu'\bigr\rrvert
_{P,2}.
\end{eqnarray*}
Given that $\mu_{n}$ and $\mu'$ are linear with at most $m$ non-zero
coefficients, then $\Delta\mu:= (\mu_{n}-\mu' )/\llvert \mu_{n}-\mu
'\rrvert _{2}$
is linear with at most $2m$-non-zero coefficients and $\llvert \Delta
\mu\rrvert _{2}=1$
by construction. Hence, in probability
\begin{eqnarray*}
\bigl\llvert\mu_{n}-\mu'\bigr\rrvert
_{P,2}^{2} & \leq& \sup_{\Delta\mu\in
\mathcal{L}_{0} (2m )\dvt \llvert \mu\rrvert _{2}\leq1}\bigl\llvert
(1-\mathbb{E} )\llvert\Delta\mu\rrvert_{n}^{2}\bigr\rrvert
\bigl\llvert\mu_{n}-\mu'\bigr\rrvert_{P,2}^{2}
\\
&&{} +\sup_{\Delta\mu\in\mathcal{L}_{0} (2m )\dvt \llvert \mu
\rrvert _{2}\leq1}\bigl\llvert2 (1-\mathbb{E} ) \bigl\langle Y-
\mu',\Delta\mu\bigr\rangle_{n}\bigr\rrvert\bigl\llvert
\mu_{n}-\mu'\bigr\rrvert_{P,2}+
\delta_{3}\bigl\llvert\mu_{n}-\mu'\bigr\rrvert
_{P,2}
\\
& \leq& \delta_{1}\bigl\llvert\mu_{n}-\mu'
\bigr\rrvert_{P,2}^{2}+ (\delta_{2}+
\delta_{3} )\bigl\llvert\mu_{n}-\mu'\bigr
\rrvert_{P,2}.
\end{eqnarray*}
Solving for $\llvert \mu_{n}-\mu'\rrvert _{P,2}^{2}$ gives the result
as long as $\delta_{1}\in[0,1)$.
\end{pf}

The next result is used to verify some of the conditions in the previous
lemma.

%
%
\begin{lemma}\label{LemmacontrolVanDeGeerConditions}Under Condition
\ref{ConditionEY|X} and either Condition~\ref{ConditionabsoluteRegularityBounded}
or~\ref{ConditionabsoluteRegularity}, for any arbitrary but fixed
$\mu'\in\mathcal{L}$, and positive integer $m$, the following hold
with probability going to one:

\begin{enumerate}
\item $\sup_{\mu\in\mathcal{L}_{0} (m )\dvt \llvert \mu\rrvert _{2}\leq
1}\llvert (1-\mathbb{E} )\llvert \mu\rrvert _{n}^{2}\rrvert \lesssim
\sqrt{\frac{m\ln K}{n}}$,

\item $\sup_{\mu\in\mathcal{L}_{0} (m )\dvt \llvert \mu\rrvert _{2}\leq
1}\llvert (1-\mathbb{E} ) \langle Y-\mu
',\mu\rangle_{n}\rrvert \lesssim\sqrt{\frac{m\ln K}{n}}$.
\end{enumerate}
\end{lemma}

\begin{pf}Let $\mathcal{S}$ be an arbitrary but fixed subset
of $ \{ 1,2,\ldots,K \} $ with cardinality $\llvert \mathcal
{S}\rrvert $.
Then, having fixed $\mathcal{S}$, $\mathfrak{F}_{\mathcal{S}}:=
\{ \mu_{\mathcal{S}}:=\sum_{k\in\mathcal{S}}b_{k}X^{ (k
)}\dvt \llvert \mu_{\mathcal{S}}\rrvert _{2}\leq A \} $
is a linear vector space of dimension~$\llvert \mathcal{S}\rrvert $.
In particular let $\Sigma_{\mathcal{S}}$ be the $m\times m$ dimensional
matrix with entries $ \{ \mathbb{E}X^{ (k )}X^{
(l )}\dvt k,l\in\mathcal{S} \} $,
and $b_{\mathcal{S}}$ the $m$ dimensional vector with entries $
\{ b_{k}\dvt k\in\mathcal{S} \} $.
Then\vspace*{1pt} $\llvert \mu_{\mathcal{S}}\rrvert _{2}^{2}=b_{\mathcal
{S}}^{T}\Sigma_{\mathcal{S}}b_{\mathcal{S}}\geq0$,
where the superscript $T$ stands for the transpose. In consequence,
$\Sigma_{\mathcal{S}}=CC^{T}$ for some $m\times m$ matrix $C$.
It follows that there is an isometry between $\mathfrak{F}_{\mathcal{S}}$
and $ \{ a\in\mathbb{R}^{m}\dvt a=C^{T}b_{\mathcal{S}} \} $.
Any vector $a$ in this last set satisfies $a^{T}a=\llvert \mu_{\mathcal
{S}}\rrvert _{2}^{2}$,
hence it is contained into the $m$ dimensional sphere of radius $A$
(under the Euclidean norm). By Lemma 14.27 in B{\"u}hlmann and van~de Geer \cite{PelUteWu07},
such sphere has a $\delta$ cover of cardinality bounded
by $ (\frac{2A+\delta}{\delta} )^{m}$ (under the Euclidean
norm). Then note that the class of functions $\mathcal{L}_{02}
(m,A ):= \{ \mu\in\mathcal{L}_{0} (m )\dvt \llvert \mu\rrvert _{2}\leq A
\} =\bigcup_{\llvert \mathcal{S}\rrvert \leq m}\mathfrak{F}_{\mathcal{S}}$.
Given that the union is over $\sum_{s=1}^{m}{
K\choose s}<mK^{m}$ number\vspace*{1pt} of elements, the covering number of $\mathcal
{L}_{02} (m,A )$
is bounded above by $mK^{m} (\frac{2A+\delta}{\delta} )^{m}$.

An argument in Loh and Wainwright (\cite{LohWai12}, proof of Lemma 15) allows
one to replace the supremum over $\mathcal{L}_{02} (m,A )$
with the maximum over a finite set. Let $ \{ \mu^{ (l
)}\dvt l=1,2,\ldots,N \} $
be an $L_{2}$ $1/3$ cover for $\mathcal{L}_{02} (m,A )$,
that is, for any $\mu\in\mathcal{L}_{02} (m,A )$ there is a
$\mu^{ (l )}$ such that $\llvert \Delta\mu\rrvert _{2}\leq1/3$,
where $\Delta\mu:=\mu-\mu^{ (l )}$. An upper bound for the
cardinality $N$ of such cover has been derived above for arbitrary
$\delta$, so for $\delta=1/3$, $N<mK^{m} (6A+1 )^{m}$.
For a $1/3$ cover, one has that $3\Delta\mu\in\mathcal{L}_{02}
(m,A )$
or equivalently $\Delta\mu\in\mathcal{L}_{02} (m,A/3 )$.
This will be used next. By adding and subtracting quantities such
as $ (1-\mathbb{E} ) \langle\mu^{ (l
)},\mu\rangle_{n}$
and using simple bounds, infer that (e.g., Loh and Wainwright \cite{LohWai12},
proof of Lemma 15),
\begin{eqnarray*}
I &:=&\sup_{\mu\in\mathcal{L}_{02} (m,A )}\bigl\llvert(1-\mathbb{E}
)\llvert\mu
\rrvert_{n}^{2}\bigr\rrvert
\\
& \leq& \max_{l\leq N}\bigl\llvert(1-\mathbb{E} )\bigl\llvert
\mu^{ (l )}\bigr\rrvert_{n}^{2}\bigr\rrvert+2\sup
_{\Delta\mu\in
\mathcal{L}_{02} (m,A/3 )}\max_{l\leq N}\bigl\llvert(1-\mathbb{E} )
\bigl\langle\mu^{ (l )},\Delta\mu\bigr\rangle_{n}\bigr\rrvert
\\
& &{} +\sup_{\Delta\mu\in\mathcal{L}_{02} (m,A/3 )}\bigl\llvert(1-\mathbb
{E} )\llvert\Delta\mu
\rrvert_{n}^{2}\bigr\rrvert
\\
& = & \max_{l\leq N}\bigl\llvert(1-\mathbb{E} )\bigl\llvert\mu
^{ (l )}\bigr\rrvert_{n}^{2}\bigr\rrvert+
\frac{2}{3}\sup_{\Delta
\mu\in\mathcal{L}_{02} (m,A )}\max_{l\leq N}\bigl
\llvert(1-\mathbb{E} ) \bigl\langle\mu^{ (l )},\Delta\mu\bigr
\rangle_{n}\bigr\rrvert
\\
& &{} +\frac{1}{9}\sup_{\Delta\mu\in\mathcal{L}_{02}
(m,A )}\bigl\llvert(1-\mathbb{E} )
\llvert\Delta\mu\rrvert_{n}^{2}\bigr\rrvert
\\
& = & \max_{l\leq N}\bigl\llvert(1-\mathbb{E} )\bigl\llvert\mu
^{ (l )}\bigr\rrvert_{n}^{2}\bigr\rrvert+
\frac{2}{3}\sup_{\mu\in
\mathcal{L}_{02} (m,A )}\bigl\llvert(1-\mathbb{E} )\llvert
\mu\rrvert_{n}^{2}\bigr\rrvert
\\
&&{} +\frac{1}{9}\sup
_{\mu\in
\mathcal{L}_{02} (m,A )}\bigl\llvert(1-\mathbb{E} )\llvert\mu\rrvert
_{n}^{2}\bigr\rrvert.
\end{eqnarray*}
This implies that $I:=\sup_{\mu\in\mathcal{L}_{02} (m,A
)}\llvert (1-\mathbb{E} )\llvert \mu\rrvert _{n}^{2}\rrvert \leq
\frac{9}{2}\max_{l\leq N}\llvert (1-\mathbb{E}
)\llvert \mu^{ (l )}\rrvert _{n}^{2}\rrvert $.
By a similar argument,
\begin{eqnarray*}
 \mathit{II} &:=&\sup_{\mu\in\mathcal{L}_{02} (m,A )}\bigl\llvert(1-\mathbb{E}
) \bigl\langle
Y-\mu',\mu\bigr\rangle_{n}\bigr\rrvert
\\
& \leq& \max_{l\leq N}\bigl\llvert(1-\mathbb{E} ) \bigl\langle
Y-\mu',\mu^{ (l )} \bigr\rangle_{n}\bigr\rrvert+
\sup_{\Delta\mu\in\mathcal{L}_{02} (m,A/3 )}\bigl\llvert(1-\mathbb{E}
) \bigl\langle Y-
\mu',\Delta\mu\bigr\rangle_{n}\bigr\rrvert
\\
& = & \max_{l\leq N}\bigl\llvert(1-\mathbb{E} ) \bigl\langle Y-
\mu',\mu^{ (l )} \bigr\rangle_{n}\bigr\rrvert+
\frac
{1}{3}\sup_{\mu\in\mathcal{L}_{02} (m,A )}\bigl\llvert(1-\mathbb{E} )
\bigl
\langle Y-\mu',\mu\bigr\rangle_{n}\bigr\rrvert
\end{eqnarray*}
implying\vspace*{1pt} $\mathit{II}:=\sup_{\mu\in\mathcal{L}_{02} (m,A )}\llvert (1-\mathbb
{E} ) \langle Y-\mu',\mu\rangle
_{n}\rrvert \leq\frac{3}{2}\max_{l\leq N}\llvert (1-\mathbb
{E} ) \langle Y-\mu',\mu^{ (l )} \rangle
_{n}\rrvert $.
Hence, to bound $I$ and $\mathit{II}$ use the above upper bounds together
with Lemma~\ref{LemmaMaximalInequ} and the upper bound for $N$
($N<mK^{m} (6A+1 )^{m}$
with $A=1$).
\end{pf}

The following is a modification of a standard crude result often used
to derive consistency, but not convergence rates. However, for the
CGA and FWA this will be enough to obtain sharp convergence rates
independently of the number of iterations $m$. Recall $\mu_{0}
(X ):=\mathbb{E} [Y|X ]$.

%
%
\begin{lemma}\label{TheoremcrudeConsistencyRates}Let $\mu'\in
\mathcal{L}$
be arbitrary, but fixed. Suppose that in probability, for some $\delta
_{1}\in(0,1 )$
and $\delta_{2},\delta_{3}>0$, and for a positive $B_{m}$:

\begin{enumerate}
\item $\sup_{\mu\in\mathcal{L} (B_{m} )}\llvert
(1-\mathbb{E} ) (\llvert Y-\mu\rrvert _{n}^{2}-\llvert Y-\mu
'\rrvert _{n}^{2} )\rrvert \leq\delta_{1}$;

\item $\llvert \mu'-\mu_{0}\rrvert _{2}^{2}\leq\delta_{2}$;

\item the sequence $\mu_{n}\in\mathcal{L} (B_{m} )$ satisfies
$\llvert Y-\mu_{n}\rrvert _{n}^{2}-\llvert Y-\mu'\rrvert
_{n}^{2}\leq \delta_{3}$.
\end{enumerate}
Then $\llvert \mu_{n}-\mu_{0}\rrvert _{P,2}\leq\sqrt{\delta
_{1}+\delta_{2}+\delta_{3}}$
in probability.
\end{lemma}

\begin{pf}By simple algebra, $\llvert Y-\mu_{n}\rrvert
_{P,2}^{2}-\llvert Y-\mu'\rrvert _{P,2}^{2}=\llvert \mu_{n}-\mu
_{0}\rrvert _{P,2}^{2}-\llvert \mu'-\mu_{0}\rrvert _{P,2}^{2}$.
Adding and subtracting $\llvert Y-\mu_{n}\rrvert _{n}^{2}-\llvert
Y-\mu
'\rrvert _{n}^{2}$,
\begin{eqnarray*}
\llvert\mu_{n}-\mu_{0}\rrvert_{P,2}^{2}
& \leq& \bigl\llvert\mu'-\mu_{0}\bigr\rrvert
_{P,2}^{2}+ \bigl[\llvert Y-\mu_{n}\rrvert
_{P,2}^{2}-\bigl\llvert Y-\mu'\bigr\rrvert
_{P,2}^{2} \bigr]- \bigl[\llvert Y-\mu_{n}\rrvert
_{n}^{2}-\bigl\llvert Y-\mu'\bigr\rrvert
_{n}^{2} \bigr]
\\
& &{} + \bigl[\llvert Y-\mu_{n}\rrvert_{n}^{2}-
\bigl\llvert Y-\mu'\bigr\rrvert_{n}^{2} \bigr]
\\
& \leq& \delta_{2}+ \bigl[\llvert Y-\mu_{n}\rrvert
_{P,2}^{2}-\bigl\llvert Y-\mu'\bigr\rrvert
_{P,2}^{2} \bigr]- \bigl[\llvert Y-\mu_{n}\rrvert
_{n}^{2}-\bigl\llvert Y-\mu'\bigr\rrvert
_{n}^{2} \bigr]+\delta_{3},
\end{eqnarray*}
where the last step follows by points 2 and 3 in the lemma. However,
\begin{eqnarray*}
& & \bigl(\llvert Y-\mu_{n}\rrvert_{P,2}^{2}-
\bigl\llvert Y-\mu'\bigr\rrvert_{P,2}^{2}
\bigr)- \bigl(\llvert Y-\mu_{n}\rrvert_{n}^{2}-
\bigl\llvert Y-\mu'\bigr\rrvert_{n}^{2} \bigr)
\\
&&\quad \leq \sup_{\mu\in\mathcal{L} (B_{m} )}\bigl\llvert\llvert Y-\mu
\rrvert
_{P,2}^{2}-\bigl\llvert Y-\mu'\bigr\rrvert
_{P,2}^{2}- \bigl(\llvert Y-\mu\rrvert_{n}^{2}-
\bigl\llvert Y-\mu'\bigr\rrvert_{n}^{2} \bigr)
\bigr\rrvert
\\
&&\quad  =  \sup_{\mu\in\mathcal{L} (B_{m} )}\bigl\llvert(1-\mathbb{E} )
\bigl(\llvert Y-
\mu\rrvert_{n}^{2}-\bigl\llvert Y-\mu'\bigr
\rrvert_{n}^{2} \bigr)\bigr\rrvert\leq\delta_{1},
\end{eqnarray*}
where the last inequality follows by assumption. Putting everything
together the result follows.
\end{pf}

In what follows, define $\mathcal{L}_{01} (m,B ):=\mathcal
{L}_{0} (m )\cap\mathcal{L}_{1} (B )$,
where $\mathcal{L}_{1} (B )=\mathcal{L} (B )$ the
usual linear space of functions with absolute sum of coefficients
bounded by $B$. The next result will be used to verify the conditions
of the previous lemma in the case of the CGA and FWA but also as main
ingredient to derive consistency rates for non-mixing data in a variety
of situations.

%
%
\begin{lemma}\label{Lemmaulln}Suppose Condition~\ref{ConditionEY|X}.
For any arbitrary, but fixed $\mu'\in\mathcal{L}_{01}
(m,B )$,
and $B_{m}<\infty$,
\begin{eqnarray*}
\mathbb{E}\sup_{\mu\in\mathcal{L}_{01} (m,B_{m} )\dvt \llvert \mu-\mu
'\rrvert _{2}\leq\delta}\bigl\llvert(1-\mathbb{E} ) \bigl(\bigl
\llvert Y-\mu(X )\bigr\rrvert_{n}^{2}-\bigl\llvert Y-\mu
' (X )\bigr\rrvert_{n}^{2} \bigr)\bigr\rrvert&
\lesssim& \error (\delta),
\end{eqnarray*}
where, under either Condition~\ref{ConditionabsoluteRegularityBounded}
or~\ref{ConditionabsoluteRegularity},
\[
\error (\delta)=\min\biggl\{ \delta\sqrt{\frac{m}{\rho
_{2m}}},B+B_{m}
\biggr\} \biggl(1+\min\biggl\{ \delta\sqrt{\frac
{m}{\rho_{2m}}},B+B_{m}
\biggr\} \biggr) \biggl(\sqrt{\frac{\ln
K}{n}} \biggr)
\]
while under Condition~\ref{Conditiondependence},
\[
\error (\delta)=\min\biggl\{ \delta\sqrt{\frac{m}{\rho
_{2m}}},B+B_{m}
\biggr\} \biggl(1+K^{1/p}\min\biggl\{ \delta\sqrt{
\frac{m}{\rho_{2m}}},B+B_{m} \biggr\} \biggr) \biggl(\frac
{d_{n,p}K^{1/p}}{\sqrt{n}}
\biggr).
\]
\end{lemma}

\begin{pf}Note that $Y=\mu_{0}+Z$, where $Z$ is mean zero conditionally
on $X$. Then, by standard algebra
\begin{eqnarray*}
& & (1-\mathbb{E} )\llvert Y-\mu\rrvert_{n}^{2}- (1-
\mathbb{E} )\bigl\llvert Y-\mu'\bigr\rrvert_{n}^{2}
\\
&&\quad =  \frac{1}{n}\sum_{i=1}^{n}2Z_{i}
\bigl(\mu' (X_{i} )-\mu(X_{i} ) \bigr)
\\
&&\qquad{} +
\frac{1}{n}\sum_{i=1}^{n} (1-
\mathbb{E} ) \bigl(\mu(X_{i} )-\mu' (X_{i} )
\bigr) \bigl(\mu(X_{i} )+\mu' (X_{i} )-2
\mu_{0} (X_{i} ) \bigr)
\\
&&\quad =:  I+\mathit{II},
\end{eqnarray*}
using the fact that $\mathbb{E} [Z|X ]=0$ in the equality.
The two terms above can be bounded separately, uniformly in $\mu$
such that $\llvert \mu-\mu'\rrvert _{2}\leq\delta$. First, let $\mu
' (X )=\sum_{k=1}^{K}b{}_{k}'X_{i}^{ (k )}$,
where by definition of $\mathcal{L}_{01} (m,B )$, only $m$
coefficients are non-zero. Note that for $\mu(X )=\sum
_{k=1}^{K}b_{k}X^{ (k )}$
in $\mathcal{L}_{01} (m,B )$, $ (\mu'-\mu)\in
\mathcal{L}_{01} (2m,B+B_{m} )$,
because $\mu$ and $\mu'$ are arbitrary, hence do not need to have
any variables in common for $2m\leq K$ (recall that there are $K$
variables $X^{ (k )}$, $k\leq K$). Define $c_{k}:=\operatorname{sign}
(b{}_{k}'-b{}_{k} )\sum_{k=1}^{K}\llvert b{}_{k}'-b{}_{k}\rrvert $,
$\lambda_{k}:= (b{}_{k}'-b{}_{k} )/\sum_{k=1}^{K}\llvert
b{}_{k}'-b{}_{k}\rrvert $,
where there are at most $2m$ non-zero $\lambda_{k}$'s by the restriction
imposed by $\mathcal{L}_{01} (2m,B+B_{m} )$. Hence,
\begin{eqnarray*}
\mu' (X )-\mu(X ) & = & \sum_{k=1}^{K}
\bigl(b{}_{k}'-b{}_{k} \bigr)X^{ (k )}=
\sum_{k=1}^{K}\lambda_{k}c_{k}X^{ (k )},
\end{eqnarray*}
with $\llvert c_{k}\rrvert \leq\llvert \mu'-\mu\rrvert _{\mathcal{L}}$,
and $\lambda_{k}$'s in the $2m$ dimensional unit simplex. Given
this restrictions, also note that
\[
\sqrt{\frac{\rho_{2m}}{2m}}\sum_{k=1}^{K}
\bigl\llvert b{}_{k}'-b{}_{k}\bigr\rrvert\leq
\sqrt{\rho_{2m}\sum_{k=1}^{K}
\bigl(b{}_{k}'-b{}_{k} \bigr)^{2}}
\leq\bigl\llvert\mu'-\mu\bigr\rrvert_{2}
\]
so that for any $\delta>0$, $\llvert \mu-\mu'\rrvert _{2}\leq\delta$
implies $\llvert \mu'-\mu\rrvert _{\mathcal{L}}\leq\delta\sqrt
{2m/\rho_{2m}}$
or equivalently $\llvert c_{k}\rrvert \leq\min\{ \delta\sqrt
{2m/\rho_{2m}},B+B_{m} \} $.
Going from right to left, the above inequality is obtained from the
Rayleigh quotient, and by bounding the $l_{1}$ norm by $\sqrt{2m}$
times the $l_{2}$ norm (e.g., use Jensen inequality of Cauchy--Schwarz).
To ease notation, write $\sup_{\llvert \mu-\mu'\rrvert _{2}\leq
\delta}$
for $\sup_{\mu\in\mathcal{L}_{01} (m,B )\dvt \llvert \mu-\mu
'\rrvert _{2}\leq\delta}$.
Then, using the previous remarks, and also noting that the supremum
over the unit simplex is achieved at one of the edges of the simplex,
\begin{eqnarray*}
\mathbb{E}\sup_{\llvert \mu-\mu'\rrvert _{2}\leq\delta}\llvert
I\rrvert& = & 2\mathbb{E}\sup
_{\llvert \mu-\mu'\rrvert _{2}\leq
\delta}\Biggl\llvert\frac{1}{n}\sum
_{i=1}^{n}Z_{i} \Biggl(\sum
_{k=1}^{K} \bigl(b{}_{k}'-b{}_{k}
\bigr)X_{i}^{ (k )} \Biggr)\Biggr\rrvert
\\[-2pt]
& = & 2\mathbb{E}\sup_{\llvert \sum_{k=1}^{K}\lambda
_{k}c_{k}X^{ (k )}\rrvert _{2}\leq\delta}\Biggl\llvert\sum
_{k=1}^{K}\lambda_{k}c_{k}
\frac{1}{n}\sum_{i=1}^{n}Z_{i}X_{i}^{
(k )}
\Biggr\rrvert
\\[-2pt]
& = & 2\mathbb{E}\max_{k\leq K}\sup_{\llvert c_{k}\rrvert \leq\min
\{ \delta\sqrt{2m/\rho_{2m}},B+B_{m} \} }
\Biggl\llvert c_{k}\frac{1}{n}\sum_{i=1}^{n}Z_{i}X_{i}^{ (k )}
\Biggr\rrvert
\\[-2pt]
& = & 2\min\biggl\{ \delta\sqrt{\frac{2m}{\rho
_{2m}}},B+B_{m} \biggr\}
\mathbb{E}\max_{k\leq K}\Biggl\llvert\frac
{1}{n}\sum
_{i=1}^{n}Z_{i}X_{i}^{ (k )}
\Biggr\rrvert.
\end{eqnarray*}
Hence, it is sufficient to bound the expectation of the sequence $
(Z_{i}X_{i}^{ (k )} )_{i\geq1}$,
which is mean zero by construction. Under Conditions~\ref{ConditionabsoluteRegularityBounded}
or~\ref{ConditionabsoluteRegularity}, $\mathbb{E}\max_{k\leq K}\llvert
\frac{1}{n}\sum_{i=1}^{n}Z_{i}X_{i}^{ (k )}\rrvert \lesssim\sqrt{\frac
{\ln K}{n}}$,
by Lemma~\ref{LemmaMaximalInequ}, while under Condition~\ref{Conditiondependence},
\[
\mathbb{E}\max_{k\leq K}\Biggl\llvert\frac{1}{n}\sum
_{i=1}^{n}Z_{i}X_{i}^{ (k )}
\Biggr\rrvert\lesssim K^{1/p}\max_{k\leq K} \Biggl(
\mathbb{E}\Biggl\llvert\frac{1}{n}\sum_{i=1}^{n}Z_{i}X_{i}^{ (k )}
\Biggr\rrvert^{p} \Biggr)^{1/p}\lesssim\frac{d_{n,p}K^{1/p}}{\sqrt{n}}
\]
by Lemma~\ref{LemmaBurkholderInequality}. To bound the terms in
$\mathit{II}$, note that
\[
\mu+\mu'-2\mu_{0}=\mu-\mu'+2 \bigl(
\mu'-\mu_{0} \bigr).
\]
Then, recalling $\Delta(X ):= (\mu' (X
)-\mu_{0} (X ) )$,
\begin{eqnarray*}
\mathbb{E}\sup_{\llvert \mu-\mu'\rrvert _{2}\leq\delta}\llvert
\mathit{II}\rrvert& \leq& \mathbb{E}\sup
_{\llvert \sum_{k=1}^{K}\lambda
_{k}c_{k}X^{ (k )}\rrvert _{2}\leq\delta}\Biggl\llvert\frac
{1}{n}\sum
_{i=1}^{n} (1-\mathbb{E} ) \Biggl(\sum
_{k=1}^{K}\lambda_{k}c_{k}X_{i}^{ (k )}
\Biggr)^{2}\Biggr\rrvert
\\[-2pt]
& &{} +\mathbb{E}\sup_{\llvert \sum_{k=1}^{K}\lambda_{k}c_{k}X^{
(k )}\rrvert _{2}\leq\delta}\Biggl\llvert\frac{2}{n}\sum
_{i=1}^{n} (1-\mathbb{E} ) \Biggl(\sum
_{k=1}^{K}\lambda_{k}c_{k}X_{i}^{ (k )}
\Biggr)\Delta(X_{i} )\Biggr\rrvert
\\[-2pt]
& =: & \mathit{III}+\mathit{IV}.
\end{eqnarray*}
Using arguments similar for the bound of $I$,
\begin{eqnarray*}
\mathit{III} & \leq& \mathbb{E}\sup_{\llvert \sum_{k=1}^{K}\lambda
_{k}c_{k}X^{ (k )}\rrvert _{2}\leq\delta,\llvert \sum
_{l=1}^{K}\lambda_{l}c_{l}X^{ (l )}\rrvert _{2}\leq\delta
}\Biggl\llvert\sum
_{k=1}^{K}\lambda_{k}c_{k}\sum
_{l=1}^{K}\lambda_{l}c_{l}
\frac{1}{n}\sum_{i=1}^{n} (1-
\mathbb{E} )X_{i}^{ (k )}X_{i}^{ (l )}\Biggr
\rrvert
\\[-2pt]
& = & \mathbb{E}\max_{k,l\leq K}\sup_{\llvert c_{k}\rrvert,\llvert
c_{l}\rrvert \leq\min\{ \delta\sqrt{2m/\rho
_{2m}},B+B_{m} \} }\Biggl
\llvert c_{k}c_{l}\frac{1}{n}\sum
_{i=1}^{n} (1-\mathbb{E} )X_{i}^{ (k
)}X_{i}^{ (l )}
\Biggr\rrvert
\\[-2pt]
& \leq& \biggl(\min\biggl\{ \delta\sqrt{\frac{2m}{\rho
_{2m}}},B+B_{m}
\biggr\} \biggr)^{2}\mathbb{E}\max_{k,l\leq K}\Biggl\llvert
\frac{1}{n}\sum_{i=1}^{n} (1-
\mathbb{E} )X_{i}^{
(k )}X_{i}^{ (l )}\Biggr
\rrvert.
\end{eqnarray*}
To finish the control of $\mathit{III}$, one can then proceed along the lines
of the control of the $I$ term:
\[
\mathbb{E}\max_{k,l\leq K}\Biggl\llvert\frac{1}{n}\sum
_{i=1}^{n} (1-\mathbb{E} )X_{i}^{ (k )}X_{i}^{ (l
)}
\Biggr\rrvert\lesssim\cases{ \displaystyle\sqrt{\ln K^{2}}, &\quad
under Condition~\ref{ConditionabsoluteRegularityBounded} or~\ref{ConditionabsoluteRegularity},
\vspace*{3pt}\cr
\displaystyle
\frac{d_{n,p}K^{2/p}}{\sqrt{n}}, &\quad under Condition~\ref{Conditiondependence}.}
\]
Similar arguments are used to bound $\mathit{IV}$. Putting these bounds together,
and disregarding irrelevant constants, the result follows.
\end{pf}

\subsection{Proof of theorems}

\begin{pf*}{Proof of Theorem~\ref{TheoremabsoluteRegularity}} At first,
prove the result for the PGA, OGA and RGA. The estimators satisfy
$F_{m}\in\mathcal{L}_{0} (m )$. Hence, apply Lemma~\ref{LemmavanDeGeerConsistency}.
Verify points 1--2 in Lemma~\ref{LemmavanDeGeerConsistency}, using
Lemma~\ref{LemmacontrolVanDeGeerConditions}, so that $\delta
_{1},\delta_{2}\lesssim\sqrt{\frac{m\ln K}{n}}$
in Lemma~\ref{LemmavanDeGeerConsistency}. By Lemmas~\ref{LemmaL2Boosting},
\ref{LemmaOGA} and~\ref{Lemmarelaxedgreedy}, point 3 in Lemma
\ref{LemmavanDeGeerConsistency} is verified with $\delta_{3}$ proportional
to $B^{1/3}m^{-1/6}$ for the PGA, $Bm^{-1/2}$ for the OGA and RGA
with $\mu'=\mu_{B}$. Point 4 is satisfied by the remark around (\ref{EQOLSB})
for $B\geq B_{0}$ as required in (\ref{EQsizeTrueParameter}). Hence,
in probability, by the triangle inequality,
\[
\llvert\mu_{0}-F_{m}\rrvert_{P,2}\lesssim
\sqrt{\frac{m\ln
K}{n}}+\llvert\mu_{0}-\mu_{B}\rrvert
_{2}+\algo (B,m ),
\]
where $\algo (B,m )$ is the appropriate error term in Lemmas
\ref{LemmaL2Boosting},~\ref{LemmaOGA},~\ref{Lemmarelaxedgreedy}.

For the CGA and FWA use Lemma~\ref{TheoremcrudeConsistencyRates}
with $\mu'=\mu_{B}\in\mathcal{L} (B )$, $B=B_{m}=\bar{B}$,
and $\mu_{n}=F_{m}$; recall $\mu_{B}$ is the minimizer in (\ref
{EQapproximationError}).
In Lemma~\ref{TheoremcrudeConsistencyRates}, $\delta_{1}\lesssim\bar
{B}\sqrt{\frac{\ln K}{n}}$
by Lemma~\ref{Lemmaulln} with $m=K$, so that $\mathcal{L}_{01}
(m,B )=\mathcal{L} (B )$.
By definition of $\mu_{B}$, in Lemma~\ref{TheoremcrudeConsistencyRates},
$\delta_{2}=\gamma^{2} (\bar{B} )$. Moreover, $\delta
_{3}\lesssim\bar{B}^{2}m^{-1}$
by Lemmas~\ref{Lemmarelaxedgreedy} and~\ref{LemmafrankWolfeApproximation}.
Hence, Lemma~\ref{TheoremcrudeConsistencyRates} is verified.
\end{pf*}

The proof of Theorem~\ref{TheoremDependenceCrude} is next.

\begin{pf*}{Proof of Theorem~\ref{TheoremDependenceCrude}}
By Lemma
\ref{LemmaEstimatorL1Bound}, $F_{m}\in\mathcal{L} (B_{m} )$
in probability, for some suitable $B_{m}$ depending on the algorithm.
The theorem then follows by an application of Lemma~\ref{TheoremcrudeConsistencyRates}
with $\mu'=\mu_{B}\in\mathcal{L} (B )$ for arbitrary $B$,
and $\mu_{n}=F_{m}$. In Lemma~\ref{TheoremcrudeConsistencyRates},
$\delta_{1}\lesssim(B+B_{m} )^{2} (\frac
{d_{n,p}K^{2/p}}{\sqrt{n}} )$
by Lemma~\ref{Lemmaulln}. Then substitute $B_{m}$ with the upper
bounds given in Lemma~\ref{LemmaEstimatorL1Bound}. Finally, in Lemma
\ref{TheoremcrudeConsistencyRates}, $\delta_{2}=\gamma(B )$
and $\delta_{3}=\algo (B,m )$ by Lemmas~\ref{LemmaL2Boosting},
\ref{LemmaOGA}, or~\ref{Lemmarelaxedgreedy}. Hence, Lemma~\ref{TheoremcrudeConsistencyRates}
and the fact that $\sqrt{\delta_{1}+\delta_{2}+\delta_{3}}\leq
\sqrt{\delta_{1}}+\sqrt{\delta_{2}}+\sqrt{\delta_{3}}$
imply the result.
\end{pf*}

Theorem~\ref{TheoremDependenceEigenvalue} relies on Theorem 3.4.1
in van der Vaart and Wellner \cite{vanWel00}, which is here recalled as a lemma
for convenience, using the present notation and adapted to the current
purposes.%
%
\begin{lemma}\label{TheoremRatesVW00}Suppose that for
any $\delta>\delta_{n}>0$, and for $B_{m}\geq B$, and fixed function
$\mu'\in\mathcal{L}_{0,1} (m,B )$:

\begin{enumerate}
\item  $\mathbb{E}\llvert Y-\mu(X )\rrvert _{n}^{2}-\mathbb
{E}\llvert Y-\mu' (X )\rrvert _{n}^{2}\gtrsim\mathbb
{E}\llvert \mu(X )-\mu' (X )\rrvert _{n}^{2}$
for any $\mu\in\mathcal{L}_{0,1} (m,B_{m} )$ such that
$\llvert \mu-\mu'\rrvert _{2}\leq\delta$;

\item $\mathbb{E}\sup_{\mu\in\mathcal{L}_{0,1} (m,B_{m}
)\dvt \llvert \mu-\mu'\rrvert _{2}\leq\delta}\llvert (1-\mathbb
{E} )\llvert Y-\mu(X )\rrvert _{n}^{2}-
(1-\mathbb{E} )\llvert Y-\mu' (X )\rrvert _{n}^{2}\rrvert \lesssim
\delta\frac{a_{n}}{n^{1/2}}$
for some sequence $a_{n}=\mathrm{o} (n^{1/2} )$;

\item there is a sequence $r_{n}$ such that $r_{n}\lesssim\delta_{n}^{-1}$
and $r_{n}\lesssim\frac{n^{1/2}}{a_{n}}$;

\item $\Pr(F_{m}\in\mathcal{L}_{0,1} (m,B_{m} )
)\rightarrow1$,
and $\llvert Y-F_{m}\rrvert _{n}^{2}\leq\llvert Y-\mu' (X
)\rrvert _{n}^{2}+\mathrm{O}_{P} (r_{n}^{-2} )$.
\end{enumerate}
Then $ (\mathbb{E}\llvert \mu_{0} (X' )-F_{m}
(X' )\rrvert _{n}^{2} )^{1/2}\lesssim\llvert \mu_{0}-\mu
'\rrvert _{2}+r_{n}^{-1}$
in probability.
\end{lemma}

Here is the proof of Theorem~\ref{TheoremDependenceEigenvalue}.

\begin{pf*}{Proof of Theorem~\ref{TheoremDependenceEigenvalue}} It
is enough to verify the conditions in Lemma~\ref{TheoremRatesVW00}
and then show that one can replace the approximation error w.r.t.
$\mu'\in\mathcal{L}_{0,1} (m,B )$ with the one w.r.t.
$\mu_{B}$.
To verify point 1 in Lemma~\ref{TheoremRatesVW00}, restrict attention
to $\mu$ such that $\mathbb{E}\llvert \mu(X )-\mu'
(X )\rrvert _{n}^{2}\geq4\mathbb{E}\llvert \mu_{0}
(X )-\mu' (X )\rrvert _{n}^{2}$.
If this is not the case, the convergence rate (error) is proportional
to $\mathbb{E}\llvert \mu_{0} (X )-\mu' (X
)\rrvert _{n}^{2}$
and Lemma~\ref{TheoremRatesVW00} would apply trivially. Hence, suppose
this is not the case. By standard algebra,
\begin{eqnarray*}
\mathbb{E}\bigl\llvert Y-\mu(X )\bigr\rrvert_{n}^{2}-
\mathbb{E}\bigl\llvert Y-\mu' (X )\bigr\rrvert_{n}^{2}
& = & \mathbb{E}\bigl\llvert\mu_{0} (X )-\mu(X )\bigr\rrvert
_{n}^{2}-\mathbb{E}\bigl\llvert\mu_{0} (X )-
\mu' (X )\bigr\rrvert_{n}^{2}
\\
& \geq& \tfrac{1}{4}\mathbb{E}\bigl\llvert\mu(X )-\mu' (X )
\bigr\rrvert_{n}^{2},
\end{eqnarray*}
where the inequality follows by problem 3.4.5 in van der Vaart and
Wellner \cite{vanWel00}. Hence, point 1 in Lemma~\ref{TheoremRatesVW00} is
satisfied. By construction, $F_{m}$ has at most $m$ non-zero coefficients.
By this remark and Lemma~\ref{LemmaEstimatorL1Bound}, $F_{m}\in
\mathcal{L}_{0,1} (m,B_{m} )$
with $B_{m}=\mathrm{O}_{p} (m^{1/2} )$ for the PGA, $B_{m}=\mathrm{O}_{p}
(m^{1/2}/\rho_{m}^{1/2} )$
for the OGA, and $B_{m}=\mathrm{O}_{p} (m^{1/2}/\rho_{m,n}^{1/2} )$
for the RGA if $\algo (B,m )=B^{2}/m=\mathrm{O} (1 )$, which
holds by the conditions in the theorem. The equality $\algo
(B,m )=B^{2}/m$
follows by Lemma~\ref{Lemmarelaxedgreedy}. By Lemma~\ref{LemmaeigenValue},
if $d_{n,p}mK^{2/p}n^{-1/2}=\mathrm{o} (\rho_{m} )$, then $\rho
_{m,n}^{-1}=\mathrm{O}_{p} (\rho_{m}^{-1} )$.
By the conditions in the theorem, $\rho_{m}>0$, and
$d_{n,p}mK^{2/p}n^{-1/2}\lesssim \error (B,K,n,m )=\mathrm{o}
(1 )$.
Hence, infer that $B_{m}=\mathrm{O}_{p} (m^{1/2} )$ for the OGA and
RGA. Hence, point 2 in Lemma~\ref{TheoremRatesVW00}, is satisfied
for any $\delta$ and $a_{n}=m^{1/2} (B+m^{1/2} )d_{n,p}K^{2/p}$
by Lemma~\ref{Lemmaulln}, where
\[
\error (\delta)\lesssim\delta m^{1/2} \bigl(B+m^{1/2} \bigr)
\biggl(\frac{d_{n,p}K^{2/p}}{\sqrt{n}} \biggr)
\]
using the fact that $\rho_{2m}>0$ and $B_{m}=\mathrm{O}_{p} (m^{1/2} )$.

It follows that point 3 in Lemma~\ref{TheoremRatesVW00} is satisfied
by $r_{n}=n^{1/2}/ [m^{1/2} (B+m^{1/2}
)d_{n,p}K^{2/p} ]$.
Moreover, by Lemma~\ref{LemmaL2Boosting},~\ref{LemmaOGA} and
\ref{Lemmarelaxedgreedy}
\[
\llvert Y-F_{m}\rrvert_{n}^{2}\leq\bigl\llvert
Y-\mu' (X )\bigr\rrvert_{n}^{2}+\mathrm{O}_{p}
\bigl(u_{n}^{-2} \bigr),
\]
with $u_{n}^{-2}$ as given in the aforementioned lemmas because $\mu
'\in\mathcal{L}_{0,1} (m,B )\subseteq\mathcal{L}
(B )$.
Since point~4 in Lemma~\ref{TheoremRatesVW00} requires $u_{n}=\mathrm{O}
(r_{n} )$,
the actual rate of convergence is $u_{n}^{-1}\vee r_{n}^{-1}\leq
u_{n}^{-1}+r_{n}^{-1}$
as stated in the theorem.

It is now necessary to replace the approximation error $\mathbb
{E}\llvert \mu_{0} (X )-\mu' (X )\rrvert _{n}^{2}$
with $\gamma(B ):=\mathbb{E}\llvert \mu_{0} (X
)-\mu_{B} (X )\rrvert _{n}^{2}$.
To this end, consider Lemmas~\ref{LemmaL2Boosting},~\ref{LemmaOGA}
and~\ref{Lemmarelaxedgreedy} with the empirical norm $\llvert \cdot
\rrvert _{n}$
replaced by $\llvert \cdot\rrvert _{P,2}$. Going through the proof,
the results are seen to hold as well with the same error rate (implicitly
using Condition~\ref{ConditionEY|X}). Hence, note that, by standard
algebra,
\[
\bigl\llvert Y-\mu' (X )\bigr\rrvert_{P,2}^{2}-
\mathbb{E}\bigl\llvert Y-\mu_{B} (X )\bigr\rrvert_{P,2}^{2}=
\mathbb{E}\bigl\llvert\mu_{0} (X )-\mu' (X )\bigr\rrvert
_{P,2}^{2}-\mathbb{E}\bigl\llvert\mu_{0} (X )-
\mu_{B} (X )\bigr\rrvert_{P,2}^{2}.
\]
The above display together with the previous remark and Lemmas~\ref{LemmaL2Boosting},
\ref{LemmaOGA} and~\ref{Lemmarelaxedgreedy} imply that
\[
\mathbb{E}\bigl\llvert\mu_{0} (X )-\mu' (X )\bigr
\rrvert_{2}^{2}\leq\mathbb{E}\bigl\llvert\mu_{0}
(X )-\mu_{B} (X )\bigr\rrvert_{n}^{2}+\mathrm{O}
\bigl(u_{n}^{-2} \bigr),
\]
with $u_{n}$ as defined above. Hence, Lemma~\ref{TheoremRatesVW00}
together with the above display gives the result which is valid for
any $B$.
\end{pf*}

\subsection{Proof of Lemmas \texorpdfstring{\protect\ref{LemmaBbarapproximation}}{1},
\texorpdfstring{\protect\ref{LemmaconditionAbsoluteRegularity}}{2}
and \texorpdfstring{\protect\ref{LemmaconditionDependence}}{3}}\vspace*{-15pt}\label{SectionproofLemmataConditions}

\begin{pf*}{Proof of Lemma~\ref{LemmaBbarapproximation}}
If $B'\geq B$,
the lemma is clearly true because $\mathcal{L} (B
)\subseteq\mathcal{L} (B' )$.
Hence, assume $B'<B$. W.n.l.g.\vspace*{2pt} assume that $\sum_{k}\llvert
b_{k}\rrvert =B$,
as $\mu\in\mathcal{L} (B )$. Let $\lambda_{k}=
(\llvert b_{k}\rrvert /B )\geq0$,
and $c_{k}=B (b_{k}/\llvert b_{k}\rrvert )$. Then $\mu=\sum_{k}\lambda
_{k}c_{k}X^{ (k )}$.
Define
\[
\mu''=\sum_{k}
\lambda_{k} \biggl(\frac{B'}{B} \biggr)c_{k}X^{
(k )}
\]
and note that $\mu''\in\mathcal{L} (B' )$ by construction
and $B'/B<1$. Then
\begin{eqnarray*}
\inf_{\mu'\in\mathcal{L} (B' )}\bigl\llvert\mu'-\mu\bigr\rrvert
_{2}^{2} & \leq& \biggl\llvert\sum
_{k}\lambda_{k}c_{k}X^{ (k
)}-
\sum_{k}\lambda_{k} \biggl(
\frac{B'}{B} \biggr)c_{k}X^{
(k )}\biggr\rrvert
_{2}^{2}
\\
& = & \biggl[1- \biggl(\frac{B'}{B} \biggr) \biggr]^{2}\sum
_{k,l}\lambda_{k}c_{k}
\lambda_{l}c_{l}\mathbb{E}X^{ (k
)}X^{ (l )}
\\
& \leq& \biggl[1- \biggl(\frac{B'}{B} \biggr) \biggr]^{2}
\biggl(\sum_{k}\llvert\lambda_{k}c_{k}
\rrvert\biggr)^{2}= \biggl[1- \biggl(\frac
{B'}{B} \biggr)
\biggr]^{2}B^{2},
\end{eqnarray*}
where the second inequality follows using the fact that $\llvert
\mathbb
{E}X^{ (k )}X^{ (l )}\rrvert \leq\mathbb
{E}\llvert X^{ (k )}\rrvert ^{2}=1$
and the last equality because $\sum_{k}\llvert \lambda
_{k}c_{k}\rrvert =\sum_{k}\llvert b_{k}\rrvert =B$.
\end{pf*}

\begin{pf*}{Proof of Lemma~\ref{LemmaconditionAbsoluteRegularity}}
At
first, show (\ref{EQboundVariance}). By independence, $\llvert ZX^{ (k
)}\rrvert _{p}=\llvert Z\rrvert _{p}\llvert X^{
(k )}\rrvert _{p}<\infty$.
Let $A_{kl}$ be the $ (k,l )$ entry in $A$ and similarly
for $H_{kl}$. By stationarity, and the fact that $A$ is diagonal,
the $l$th entry in $W_{i}$ is $W_{il}=\sum_{s=0}^{\infty
}A_{ll}^{s}\varepsilon_{i-s,l}$,
and by definition $X_{i}^{ (k )}=\sum_{l=1}^{L}H_{kl}\sum_{s=0}^{\infty
}A_{ll}^{s}\varepsilon_{i-s,l}$.
Hence, by Minkowski inequality, and the fact that $A_{ll}^{s}$ decays
exponentially fast because less than one in absolute value, and the
fact that $ \{ H_{kl}\dvt l=1,2,\ldots,L \} $ is in the unit simplex,
$\llvert X^{ (k )}\rrvert _{2p}\leq\max_{l\leq L}\sum_{s=0}^{\infty
}\llvert A_{ll}^{s}\rrvert \llvert \varepsilon
_{i-s,l}\rrvert _{2p}<\infty$.
Finally, by H\"older's inequality, the Lipschitz condition for $g$,
and Minkowski inequality, for any $\mu_{B} (X )=\sum_{k=1}^{K}X^{ (k )}b_{k}$,
\begin{eqnarray*}
\bigl\llvert\Delta(X )X^{ (k )}\bigr\rrvert_{p} & \leq&
\Biggl\llvert\sum_{l=1}^{K} \bigl(
\lambda_{l}+\llvert b_{l}\rrvert\bigr)\bigl\llvert
X^{ (l )}\bigr\rrvert\Biggr\rrvert_{2p}\bigl\llvert
X^{
(k )}\bigr\rrvert_{2p}\lesssim(1+B )\max
_{l}\bigl\llvert X^{ (l )}\bigr\rrvert_{2p}
\bigl\llvert X^{ (k )}\bigr\rrvert_{2p}
\\
& \leq& (1+B )\max_{k}\bigl\llvert X^{ (k )}\bigr
\rrvert_{2p}^{2}<\infty.
\end{eqnarray*}
This completes the proof of (\ref{EQboundVariance}). Geometric absolute
regularity follows using the fact that by construction, the mixing
coefficients of $ (X_{i} )_{i\in\mathbb{Z}}$ are equal to
the mixing coefficients of $ (W_{i} )_{i\in\mathbb{Z}}$
because $X_{i}$ is just a linear transformation of $W_{i}$ (i.e.,
the sigma algebras generated by the two processes are the same). The
process $ (W_{i} )_{i\in\mathbb{Z}}$ follows a $L$ dimensional
stationary AR(1) model with i.i.d. innovations having a density
w.r.t. the Lebesgue measure. Hence, Theorem 1 in Mokkadem
\cite{Mok88} says
that the vector autoregressive process $ (W_{i} )_{i\in
\mathbb{Z}}$
is absolutely regular with geometrically decaying mixing coefficients
as long as $L$ is bounded. By independence, the sigma algebra generated
by $ (W_{i} )_{i\in\mathbb{Z}}$ and $ (Z_{i}
)_{i\in\mathbb{Z}}$
are independent. Then Theorem 5.1 Bradley \cite{Bra05}
says that the mixing
coefficient of $ (W_{i},Z_{i} )_{\in\mathbb{Z}}$ are bounded
by the sum of the mixing coefficients of $ (W_{i} )_{\in
\mathbb{Z}}$
and $ (Z_{i} )_{i\in\mathbb{Z}}$. Since the latter mixing
coefficients are zero at any non-zero lags because of independence,
geometric beta mixing follows, and Condition~\ref{ConditionabsoluteRegularity}
holds.
\end{pf*}

\begin{pf*}{Proof of Lemma~\ref{LemmaconditionDependence}}
By assumption
$X_{i}=W_{i}$. Andrews \cite{And84} and Bradley \cite{Bra86} show that
the AR(1)
model as in the lemma is not strong mixing, hence is not absolutely
regular and Condition~\ref{ConditionabsoluteRegularity} fails. Consider
each term in the sum in Condition~\ref{Conditiondependence} separately.
First, by independence, $\llvert \mathbb{E}_{0}Z_{i}X_{i}^{
(k )}\rrvert _{p}=0$
for any $i>0$. Second, using the infinite MA representation of the
AR(1), the fact that the error terms are i.i.d., and then the triangle
inequality
\begin{eqnarray*}
\bigl\llvert\mathbb{E}_{0} (1-\mathbb{E} )\bigl\llvert
X_{i}^{
(k )}\bigr\rrvert^{2}\bigr\rrvert
_{p} & = & \Biggl\llvert\sum_{s,r=i}^{\infty
}A_{kk}^{s}A_{kk}^{r}
(1-\mathbb{E} )\varepsilon_{i-s,k}\varepsilon_{i-r,k}\Biggr
\rrvert_{p}
\\
&  \leq& 2\sum_{s,r=i}^{\infty}\bigl\llvert
A_{kk}^{s}\bigr\rrvert\bigl\llvert A_{kk}^{r}
\bigr\rrvert\llvert\varepsilon_{i-s,k}\varepsilon_{i-r,k}
\rrvert_{p}\lesssim A_{kk}^{2i}
\end{eqnarray*}
which is summable. Third, by the triangle inequality,
\begin{eqnarray*}
\bigl\llvert\mathbb{E}_{0} (1-\mathbb{E} )\Delta(X_{i}
)X_{i}^{ (k )}\bigr\rrvert_{p} & = & \Biggl\llvert
\mathbb{E}_{0} (1-\mathbb{E} ) \Biggl(g \bigl(X_{i}^{ (l
)};l
\leq K \bigr)-\sum_{k=1}^{K}X_{i}^{ (l )}b_{l}
\Biggr)X_{i}^{ (k )}\Biggr\rrvert_{p}
\\
& \leq& \bigl\llvert\mathbb{E}_{0} (1-\mathbb{E} )g
\bigl(X_{i}^{ (l )};l\leq K \bigr)X_{i}^{ (k )}
\bigr\rrvert_{p}+\Biggl\llvert\mathbb{E}_{0} (1-\mathbb{E}
)X_{i}^{
(k )}\sum_{l=1}^{K}X_{i}^{ (l )}b_{l}
\Biggr\rrvert_{p}
\\
& =: & I+\mathit{II}.
\end{eqnarray*}
Consider each term separately. Define $X_{i0}^{ (k )}:=\sum
_{s=0}^{i-1}A_{kk}^{s}\varepsilon_{i-s,k}$
and $X_{i1}^{ (k )}:=\sum_{s=i}^{\infty
}A_{kk}^{s}\varepsilon_{i-s,k}$,
and note that $X_{i}^{ (k )}=X_{i0}^{ (k
)}+X_{i1}^{ (k )}$.
By simple algebraic manipulations and repeated use of Minkowski inequality,
\begin{eqnarray*}
I & \leq& \bigl\llvert\mathbb{E}_{0} (1-\mathbb{E} )g
\bigl(X_{i0}^{ (l )};l\leq K \bigr)X_{i}^{ (k
)}
\bigr\rrvert_{p}+\bigl\llvert\mathbb{E}_{0} (1-\mathbb{E}
) \bigl(g \bigl(X_{i}^{ (l )};l\leq K \bigr)-g
\bigl(X_{i0}^{
(l )};l\leq K \bigr) \bigr)X_{i}^{ (k )}
\bigr\rrvert_{p}
\\
& \leq& \bigl\llvert(\mathbb{E}_{0}-\mathbb{E} )g
\bigl(X_{i0}^{ (l )};l\leq K \bigr)X_{i0}^{ (k
)}
\bigr\rrvert_{p}+\bigl\llvert(\mathbb{E}_{0}-\mathbb{E} )g
\bigl(X_{i0}^{ (l )};l\leq K \bigr)X_{i1}^{ (k
)}
\bigr\rrvert_{p}
\\
& &{} +\bigl\llvert(\mathbb{E}_{0}-\mathbb{E} ) \bigl(g
\bigl(X_{i}^{ (l )};l\leq K \bigr)-g \bigl(X_{i0}^{
(l )};l
\leq K \bigr) \bigr)X_{i0}^{ (k )}\bigr\rrvert_{p}.
\end{eqnarray*}
The first term on the right-hand side of the second inequality is zero
by construction
when taking expectations $ (\mathbb{E}_{0}-\mathbb{E} )$.
To bound the second term, note that by the properties of $g$, and
Minkowski inequality,
\begin{eqnarray*}
\bigl\llvert(\mathbb{E}_{0}-\mathbb{E} )g \bigl(X_{i0}^{
(l )};l
\leq K \bigr)X_{i1}^{ (k )}\bigr\rrvert_{p} & \leq&
\biggl\llvert(\mathbb{E}_{0}+\mathbb{E} )\sum
_{l\leq
K}\lambda_{l}\bigl\llvert X_{i0}^{ (l )}
\bigr\rrvert\bigl\llvert X_{i1}^{ (k )}\bigr\rrvert\biggr
\rrvert_{p}
\\
& \leq& \biggl\llvert\biggl(\sum_{l\leq K}
\lambda_{l}\mathbb{E}_{0}\bigl\llvert X_{i0}^{ (l )}
\bigr\rrvert\biggr)\bigl\llvert X_{i1}^{ (k
)}\bigr\rrvert
\biggr\rrvert_{p}+\biggl\llvert\sum_{l\leq K}
\lambda_{l}\mathbb{E}\bigl\llvert X_{i0}^{ (l )}
\bigr\rrvert\mathbb{E}\bigl\llvert X_{i1}^{ (k )}\bigr\rrvert
\biggr\rrvert_{p}
\\
& \lesssim& \bigl\llvert X_{i1}^{ (k )}\bigr\rrvert
_{p}\lesssim\sum_{s=i}^{\infty}\bigl
\llvert A_{kk}^{s}\bigr\rrvert
\end{eqnarray*}
by the independence of $X_{i0}^{ (k )}$ and $X_{i1}^{
(k )}$
and the existence of $p$ moments. The third term was bounded in a
similar way. Hence, $I\lesssim\sum_{s=i}^{\infty}\llvert
A_{kk}^{s}\rrvert \lesssim\llvert A_{kk}^{i}\rrvert $.
Finally,
\begin{eqnarray*}
\mathit{II} & \leq& B\max_{l}\bigl\llvert\mathbb{E}_{0}
(1-\mathbb{E} )X_{i}^{ (k )}X_{i}^{ (l )}\bigr
\rrvert_{p}
\\
& \leq& \Biggl\llvert\sum_{s,r=0}^{i-1}A_{kk}^{s}A_{ll}^{r}
(\mathbb{E}_{0}-\mathbb{E} )\varepsilon_{i-s,k}\varepsilon
_{i-r,l}\Biggr\rrvert_{p}+\Biggl\llvert\sum
_{s,r=i}^{\infty
}A_{kk}^{s}A_{ll}^{r}
(\mathbb{E}_{0}-\mathbb{E} )\varepsilon_{i-s,k}
\varepsilon_{i-r,l}\Biggr\rrvert_{p}
\\
& \lesssim& A_{kk}^{2i}
\end{eqnarray*}
as the first term is exactly zero. Clearly, both $I$ and $\mathit{II}$ are
summable and the lemma is proved.
\end{pf*}

\subsection{Proof of Example \texorpdfstring{\protect\ref{ExamplelongMemory}}{4}}\label{SectionProofExample}

For simplicity, write $\varepsilon_{i}$ in place of $\varepsilon_{i,k}$.
By independence of the $\varepsilon_{i}$'s and stationarity,
\begin{eqnarray*}
\mathbb{E}_{0} (1-\mathbb{E} )\bigl\llvert X_{i}^{ (k
)}
\bigr\rrvert^{2} & = & \mathbb{E}_{0} (1-\mathbb{E} )\sum
_{s,r=0}^{\infty}a_{s}a_{r}
\varepsilon_{i-s}\varepsilon_{i-r}=\sum
_{s,r\geq i}a_{s}a_{r} \bigl[ (1-\mathbb{E} )
\varepsilon_{i-s}\varepsilon_{i-r} \bigr]
\\
& = & \sum_{s\geq i}a_{s}^{2} (1-
\mathbb{E} )\varepsilon_{i-s}^{2}+2\sum
_{r>s\geq i}a_{s}a_{r}\varepsilon_{i-s}
\varepsilon_{i-r}=:I+\mathit{II}.
\end{eqnarray*}
For $i>0$, define $\bar{a}_{s}=i^{ (1+\epsilon)}a_{s}^{2}$
and $\bar{a}:=\sum_{s\geq i}\bar{a}_{s}$ and note that $\bar{a}$
depends on $i$ but is finite for any $i$ because $i^{ (1+\epsilon
)}a_{s}^{2}=i^{ (1+\epsilon)}s^{- (1+\epsilon
)}\leq1$
for $s\geq i$ (recall the definition of $a_{s}$). Then, by the definition
of $\bar{a}_{s}$ and then by Jensen inequality,
\begin{eqnarray*}
\mathbb{E}\llvert I\rrvert^{p} & \leq& \mathbb{E}\sum
_{s\geq
i} \biggl(\frac{\bar{a}_{s}}{\bar{a}} \biggr)\bigl\llvert
i^{-
(1+\epsilon)}\bar{a} (1-\mathbb{E} )\varepsilon_{i-s}^{2}
\bigr\rrvert^{p}
\\
& = & \sum_{s\geq i} \biggl(\frac{\bar{a}_{s}}{\bar{a}}
\biggr)i^{- (1+\epsilon)p}\bar{a}^{p}\mathbb{E}\bigl\llvert(1-\mathbb{E}
)\varepsilon_{i-s}^{2}\bigr\rrvert^{p}
\\
& \leq& i^{- (1+\epsilon)p}\bar{a}^{p}\max_{s}
\mathbb{E}\bigl\llvert(1-\mathbb{E} )\varepsilon_{s}^{2}
\bigr\rrvert^{p}
\end{eqnarray*}
because $\sum_{s\geq i} (\frac{\bar{a}_{s}}{\bar{a}} )=1$
and the $\bar{a}_{s}\geq0$. The above display implies that $\llvert
I\rrvert _{p}\lesssim i^{- (1+\epsilon)}$.
It remains to bound $\mathit{II}$. For any random variable $W$ such that
$\mathbb{E}\exp\{ \llvert W\rrvert /\tau\} \leq4$ for some
$\tau>0$, it is clear that
\[
\mathbb{E}\llvert W/\tau\rrvert^{p}\leq p! \bigl(\mathbb{E}\exp
\bigl\{ \llvert W\rrvert/\tau\bigr\} -1 \bigr)\leq p!\times3
\]
using Taylor series expansion. This implies that $ (\mathbb
{E}\llvert W\rrvert ^{p} )^{1/p}\leq3p\tau$
for such $\tau$ if it exists. Hence, apply this inequality to bound
$\mathbb{E}\llvert \mathit{II}\rrvert ^{p}$. Noting that $\mathbb{E}\exp
\{ \tau^{-1}\llvert \mathit{II}\rrvert \} \leq\mathbb{E}\exp\{
\tau^{-1}\mathit{II} \} +\mathbb{E}\exp\{ -\tau^{-1}\mathit{II} \} $,
it is enough to bound $\mathbb{E}\exp\{ \tau^{-1}\mathit{II} \} $.
By Gaussianity, independence of the $\varepsilon_{i}$'s, and the
fact that $\exp\{ \cdot\} $ is non-negative, letting
$\mathbb{E}_{i}$
be expectation conditional on $\varepsilon_{i}$ and its past,
\begin{eqnarray*}
\mathbb{E}\exp\bigl\{ \tau^{-1}\mathit{II} \bigr\} & = & \prod
_{r>s\geq
i}\mathbb{E}\exp\bigl\{ \tau^{-1}2a_{s}
\varepsilon_{i-s}a_{r}\varepsilon_{i-r} \bigr\} =
\prod_{r>s\geq i}\mathbb{E}\exp\bigl\{
\mathbb{E}_{i-r} \bigl(\tau^{-1}2a_{s}\varepsilon
_{i-s}a_{r}\varepsilon_{i-r} \bigr)^{2}
\bigr\}
\\
& = & \prod_{r>s\geq i}\mathbb{E}\exp\bigl\{ \bigl(\tau
^{-1}2a_{s}a_{r} \bigr)^{2}
\varepsilon_{i-r}^{2} \bigr\} =\prod
_{r>s\geq i}\mathbb{E}\exp\biggl\{ 4\bar{a}^{2}i^{-2
(1+\epsilon)}
\tau^{-2}\frac{\bar{a}_{s}}{\bar{a}}\frac
{\bar{a}_{r}}{\bar{a}}\varepsilon_{i-r}^{2}
\biggr\},
\end{eqnarray*}
where the last three steps use the properties of the moment generating
function of a Gaussian random variable and the definition of $\bar{a}_{s}$
and $\bar{a}$, as used in the control of $I$. Hence, setting $\tau
=4\bar{a}i^{- (1+\epsilon)}$,
and recalling that $\bar{a}_{s}/\bar{a}\leq1$ by construction, the
above is then bounded by
\[
\max_{r\geq i}\mathbb{E}\exp\biggl\{ \frac{\varepsilon
_{i-r}^{2}}{4} \biggr\}
=\int_{\mathbb{R}}\mathrm{e}^{z^{2}/4}\frac
{\mathrm{e}^{-z^{2}/2}}{\sqrt{2\uppi}}\,\mathrm{d}z=\sqrt{2},
\]
where the two equalities follow from the fact that $\varepsilon_{i-r}^{2}$
is a standard normal random variable, and then performing the integration.
The above two display show that for $\tau=4\bar{a}i^{-
(1+\epsilon)}$,
$\mathbb{E}\exp\{ \tau^{-1}\llvert \mathit{II}\rrvert \} \leq
\exp\{ \tau^{-1}\mathit{II} \} +\exp\{ -\tau^{-1}\mathit{II}
\} \leq2\sqrt{2}<4$,
which implies $\llvert \mathit{II}\rrvert _{p}\lesssim\bar{a}i^{-
(1+\epsilon)}$.
The upper bounds for the $L_{p}$ norms of $I$ and $\mathit{II}$ imply that
$\llvert \mathbb{E}_{0} (1-\mathbb{E} )\llvert X_{i}^{
(k )}\rrvert ^{2}\rrvert _{p}\lesssim i^{- (1+\epsilon
)}$.




%

\printhistory
\end{document}